\documentclass[11pt,a4paper]{amsart}

\usepackage{amsfonts,amssymb,amsmath,amsthm}
\usepackage{mathtools,esint}
\usepackage{lmodern}
\usepackage[noadjust]{cite}
\usepackage{enumitem}
\usepackage{needspace}
\usepackage{xcolor}
\usepackage[plainpages=false,pdfpagelabels,backref=page]{hyperref}
\hypersetup{
 colorlinks=true,
 linkcolor={cyan!90!black},
 citecolor={magenta},
 urlcolor={green!40!black},
 pdftitle={The kinetic Kato square root problem},
 pdfauthor={Kaj Nystr\"om},
 pdfsubject={The kinetic Kato square root problem}
}

\newtheorem{thm}{Theorem}[section]
\newtheorem{lem}[thm]{Lemma}
\newtheorem{prop}[thm]{Proposition}
\newtheorem{cor}[thm]{Corollary}

\theoremstyle{definition}

\newtheorem{rem}[thm]{Remark}
\numberwithin{equation}{section}

\usepackage{tikz}
\usetikzlibrary{arrows.meta,calc,patterns,decorations.pathreplacing}
\definecolor{KPBlue}{RGB}{35,82,111}
\definecolor{KPGreen}{RGB}{45,100,83}
\definecolor{KPGold}{RGB}{163,115,43}
\tikzset{
 kp/past/.style={draw=KPBlue,densely dashed,line width=.8pt,fill=KPBlue!7},
 kp/future/.style={draw=KPGreen,line width=.85pt,fill=KPGreen!6},
 kp/outer/.style={draw=black!55,dash pattern=on 4pt off 2pt,line width=.65pt},
 kp/arrow/.style={-{Stealth[length=2mm]},draw=black!65,line width=.7pt},
 kp/step/.style={-{Stealth[length=2.2mm]},draw=KPBlue,line width=1pt},
 kp/vertex/.style={draw=black!70,rounded corners=1.5pt,fill=white,
 minimum width=2.05cm,minimum height=1.02cm,inner sep=3pt,
 align=center,font=\small},
 kp/terminal/.style={kp/vertex,double,double distance=1.2pt},
 kp/small/.style={font=\footnotesize,align=center},
 kp/heading/.style={font=\small\bfseries,anchor=west},
 kp/leader/.style={draw=black!65,line width=.5pt},
 kp/axis/.style={-{Stealth[length=1.7mm]},line width=.55pt,draw=black!75}
}

\newcommand{\R}{\mathbb R}
\newcommand{\C}{\mathbb C}

\usepackage{mathrsfs}
\renewcommand{\d}{\,\mathrm{d}}
\newcommand{\V}{\mathrm V}

\DeclareMathOperator{\sgn}{sgn}

\DeclareMathOperator{\Ran}{Ran}
\allowdisplaybreaks
\makeatletter
\renewcommand{\l@section}{%
  \@tocline{1}{3pt}{0pt}{}{\normalfont\small}}
\makeatother

\begin{document}

\title[The kinetic Kato square root problem]
{The kinetic Kato square root problem}
\author{Kaj Nystr\"om}
\address{Department of Mathematics, Uppsala University, Box 480\\
SE-751 06 Uppsala, Sweden}
\email{kaj.nystrom@math.uu.se}
\date{}
\renewcommand{\subjclassname}{\textup{2020} Mathematics Subject Classification}

\begin{abstract}
We solve the kinetic Kato square root problem for the whole-space
Kolmogorov operator
$$H_A=\partial_t-X\cdot\nabla_Y-\operatorname{div}_X(A\nabla_X)$$
on \(\mathrm L^2\) with Lebesgue measure, where \(A\) is a bounded
measurable complex uniformly accretive matrix depending on all variables.
We define its realization through a sesquilinear form and identify the
square-root domains of the resulting maximal accretive operator and
its adjoint with the intrinsic Sobolev space defined by one velocity
derivative and half a derivative along the transport flow. The
homogeneous norm equivalences depend only on the dimension and the
coefficient bounds. The energy also controls a position derivative
of order one third.
The proof uses second-order resolvents and localized affine testing.
A regularized transport symmetrizer gives hidden coercivity, while
comparison with a positive Rockland operator supplies intrinsic
smoothing. Commutator and dyadic approximation estimates remove
the smoothing before a Carleson argument is applied. We also obtain
uniform shifted square-root estimates and bounded velocity, transport
and position Riesz transforms.
\end{abstract}

\maketitle
\vspace{-1\baselineskip}
\enlargethispage{2pt}
\begingroup
\hypersetup{linkcolor=black}
\tableofcontents
\endgroup

\vspace{-2\baselineskip}

\section{Introduction and statement of the main result}\label{sec:introduction}

The Kato square root problem asks whether the domain of the square
root of a differential operator coincides with its natural energy
space, with equivalent norms. Such an identification connects the
functional calculus of the operator with the differential regularity
underlying its variational formulation. For operators with rough, nonsymmetric coefficients, this
identification does not follow from abstract accretivity
and requires additional analysis.  In the kinetic setting, establishing square-root estimates
requires accounting for both degenerate diffusion and the
transfer of regularity through transport. To the best of
our knowledge, the corresponding square-root domain identification
has not previously been established for kinetic Kolmogorov operators
with bounded measurable, uniformly accretive complex coefficients.
We fill this gap by identifying the square-root domains of the
operator and its adjoint with the natural intrinsic energy space,
with equivalent norms. By analogy with the elliptic and parabolic theories, we expect
our results to have far-reaching consequences for the analysis
of kinetic boundary value problems with rough coefficients.

We study the Kato square root problem for the Kolmogorov operator
\begin{equation}\label{eq:operator-intro}
 H_A:=T-\operatorname{div}_X(A\nabla_X),
 \qquad T:=\partial_t-X\cdot\nabla_Y,
\end{equation}
on $\mathrm L^2(G)$, where $G:=\R^m\times\R^m\times\R$
and $m\ge1$ is an integer. We write $z=(X,Y,t)\in G$
and use Lebesgue measure $\d z=\d X\d Y\d t$. The measurable matrix
$A\colon G\to\C^{m\times m}$ satisfies
\begin{equation}\label{eq:ellipticity}
 |A(z)\zeta|\le a_1|\zeta|,
 \qquad
 \operatorname{Re}\bigl(A(z)\zeta\cdot\overline\zeta\bigr)
       \ge a_0|\zeta|^2
 \quad(\zeta\in\C^m)
\end{equation}
almost everywhere, for fixed $0<a_0\le a_1<\infty$.
The coefficients may depend measurably on all variables. We use
the accretive sign convention in \eqref{eq:operator-intro}.

For $\lambda>0$, the intrinsic dilations are \[\delta_\lambda(X,Y,t)=(\lambda X,\lambda^3Y,\lambda^2t).\] Consequently, horizontal differentiation has degree one and $T$ has degree two. The natural first-order
energy  combines one derivative in $X$ with half a derivative along the transport flow.  More precisely,  the transport flow associated with $T$ acts on functions by
\[
 (e^{sT}u)(X,Y,t)=u(X,Y-sX,t+s),\qquad s\in\R.
\]
These operators form a strongly continuous unitary group on $\mathrm L^2(G)$.
Its infinitesimal generator is the
skew-adjoint realization of $T$, with $T^*=-T$. The spectral
theorem defines
\[
 |T|=(T^*T)^{1/2}=(-T^2)^{1/2},
\]
as the unique nonnegative self-adjoint square root of
$T^*T=-T^2$. The energy space and its homogeneous energy norm are
\begin{equation}\label{eq:energy-space}
 \V:=\mathrm D(\nabla_X)\cap\mathrm D(|T|^{1/2}),
 \qquad E(u):=\|\nabla_Xu\|_2+\||T|^{1/2}u\|_2.
\end{equation}
Here
\[
 \mathrm D(\nabla_X)
 :=\{u\in\mathrm L^2(G):
       \nabla_Xu\in\mathrm L^2(G;\C^m)\},
\]
with derivatives understood in the distributional sense, and
$\mathrm D(|T|^{1/2})$ denotes the domain of the fractional
power defined by the spectral theorem. We equip $\V$ with the Hilbert norm
\[
 \|u\|_\V
 :=\bigl(\|u\|_2^2+\|\nabla_Xu\|_2^2
                    +\||T|^{1/2}u\|_2^2\bigr)^{1/2}
 \sim\|u\|_2+E(u).
\]
All inner products are linear in the first argument, and all
norms on $G$ use Lebesgue measure unless explicitly stated
otherwise. These operators and the energy space $\V$ are
discussed in detail in Section~\ref{sec:geometry-energy}.

We define the operator $H_A$ through a bounded sesquilinear form on $\V$.
Let $J_0=T|T|^{-1}$ denote the bounded polar factor of $T$.
After partial Fourier transformation in $(Y,t)$ it has multiplier
$i\sgn(\tau-X\cdot\eta)$, with $\sgn(0)=0$. Set
\begin{equation}\label{eq:form-intro}
 \mathfrak a_A(u,v)
 =\langle A\nabla_Xu,\nabla_Xv\rangle
  +\langle J_0|T|^{1/2}u,|T|^{1/2}v\rangle,
 \qquad u,v\in\V.
\end{equation}
The second term is the continuous extension of $\langle Tu,v\rangle$
from smooth functions. The $\mathrm L^2$ realization of
\eqref{eq:operator-intro} is defined by
\begin{equation}\label{eq:realization-intro}
 \begin{split}
 \mathrm D(H_A)&=\{u\in\V:\text{ there is }f\in\mathrm L^2(G)
       \text{ with }\mathfrak a_A(u,v)=\langle f,v\rangle
       \text{ for all }v\in\V\},\quad H_Au=f.
 \end{split}
\end{equation}
The representing function is unique because $\V$ contains
$\mathrm C_0^\infty(G)$, which is dense in $\mathrm L^2(G)$.

The adjoint of $\mathfrak a_A$ is the form
\begin{equation}\label{eq:adjoint-form-intro}
 \begin{split}
 \mathfrak a_A^*(u,v)
 &=\overline{\mathfrak a_A(v,u)}=\langle A^*\nabla_Xu,\nabla_Xv\rangle
   -\langle J_0|T|^{1/2}u,|T|^{1/2}v\rangle,
 \qquad u,v\in\V,
 \end{split}
\end{equation}
where $A^*$ denotes the pointwise conjugate transpose of $A$.
We define its $\mathrm L^2$ realization $H_A^*$ by
\eqref{eq:realization-intro}, with $\mathfrak a_A$ replaced by
$\mathfrak a_A^*$. Its formal differential expression is
\begin{equation}\label{eq:adjoint-operator-intro}
 H_A^*=-T-\operatorname{div}_X(A^*\nabla_X).
\end{equation}
As proved in Proposition~\ref{prop:kinetic-realization}, this
realization is the Hilbert-space adjoint of $H_A$.

Our main result solves the kinetic Kato square root problem.

\begin{thm}\label{thm:main}
Assume \eqref{eq:ellipticity}. The variational realizations $H_A$
and $H_A^*$ defined above are densely defined and maximal accretive
on $\mathrm L^2(G)$ and are Hilbert-space adjoints of one another.
The domains of their maximal accretive square roots coincide with
the common form domain, that is,
\begin{equation}\label{eq:main-domain}
 \mathrm D(H_A^{1/2})=\mathrm D((H_A^*)^{1/2})=\V.
\end{equation}
For every $u\in\V$,
\begin{equation}\label{eq:main-norm}
 \|H_A^{1/2}u\|_2
 \sim \|\nabla_Xu\|_2+\||T|^{1/2}u\|_2
 \sim \|(H_A^*)^{1/2}u\|_2.
\end{equation}
The implicit constants depend only on $m,a_0,a_1$.
\end{thm}

In Theorem~\ref{thm:main}, no symmetry, smallness or continuity assumption is imposed on $A$.
The half-transport derivative $|T|^{1/2}$ in
\eqref{eq:energy-space} measures half-order regularity along
the kinetic transport flow generated by
$T=\partial_t-X\cdot\nabla_Y$.

The real part of $\mathfrak a_A(u,u)$ in \eqref{eq:form-intro}
controls the horizontal energy $\|\nabla_Xu\|_2^2$ but provides
no coercive control of $\||T|^{1/2}u\|_2^2$. It therefore does
not control the full energy norm on $\V$, and the usual
representation theorem for closed sectorial forms is not
directly applicable. Our variational construction of $H_A$
relies on \emph{hidden coercivity}: we modify the test function
so that the skew-adjoint part of the operator contributes
positively to the energy estimate. This principle originates
in Kaplan's variational treatment of parabolic equations
\cite{Kaplan} and is used in the parabolic Kato theory
\cite{Nystrom2016,AuscherEgertNystrom2020,AtaeiEgertNystrom,AtaeiNystromBMO}.

Implementing this principle in the kinetic setting requires an
additional construction. The natural candidate $(I+\delta J_0)u$
is unsuitable because $J_0$ need not preserve horizontal Sobolev
regularity: even if $u\in\V$, the distributional gradient
$\nabla_X(J_0u)$ need not belong to $\mathrm L^2(G;\C^m)$.
Thus $J_0u$ need not be an admissible test function in the
form domain $\V$. One of our contributions is to overcome this obstruction
by constructing a regularized transport symmetrizer.  We use
\[
 \widehat{J_\varepsilon u}(X,\eta,\tau)
 :=\frac{iq}{(q^2+\varepsilon^2|\eta|^{4/3})^{1/2}}
       \widehat u(X,\eta,\tau),
 \qquad q:=\tau-X\cdot\eta,
\]
for $\varepsilon>0$ small, as the regularization of $J_0$. Here the hat denotes partial Fourier transformation in $(Y,t)$,
and the multiplier is understood to be zero when $q=0$. The kinetic uncertainty estimate of Lemma~\ref{lem:uncertainty},
\begin{equation}\label{eq:y-consequence}
 \|D_Y^{1/3}u\|_2\lesssim E(u),
 \qquad
 D_Y^{1/3}:=(-\Delta_Y)^{1/6},
 \qquad u\in\V,
\end{equation}
controls the additional term arising when the multiplier defining
$J_\varepsilon$ is differentiated with respect to $X$. In particular,
we prove that $J_\varepsilon$ is bounded on $\V$ for each fixed $\varepsilon>0$.

Choosing $\varepsilon>0$ and then $\delta>0$ sufficiently small,
with dependence only on $m,a_0,a_1$, gives a bounded invertible
map $S_+:=I+\delta J_\varepsilon$ on $\V$ such that
\begin{equation}\label{eq:hidden-intro}
 \operatorname{Re}\mathfrak a_A(u,S_+u)\gtrsim E(u)^2,
 \qquad
 \operatorname{Re}\langle u,S_+u\rangle=\|u\|_2^2,
 \qquad u\in\V.
\end{equation}
Thus, after this change of test function, the shifted form is
coercive on $\V$. The Lax--Milgram theorem yields the resolvent
and establishes maximal accretivity of $H_A$.
Section~\ref{sec:variational-realization} gives the details and
identifies the adjoint. This construction is independent of the
square-root estimates. 

The estimate \eqref{eq:y-consequence} shows that a position derivative
of order $1/3$ may also be included in the energy without changing
the square-root domain. Corollary~\ref{cor:shifted-roots} gives
uniform inhomogeneous estimates for $(\mu+H_A)^{1/2}$,
$\mu>0$, and the associated Riesz-transform bounds.

\begin{rem}\label{convention}
Throughout the paper, implicit constants depend only on
$m,a_0,a_1$ and fixed geometric choices, unless stated otherwise.
For nonnegative quantities $F$ and $G$, we write $F\lesssim G$
if $F\le CG$ for such a constant $C\ge1$, and $F\gtrsim G$
if $G\lesssim F$. We write $F\sim G$ when both inequalities hold.
Thus, $\sim$ denotes two-sided comparability, not asymptotic
equivalence. A subscript on a comparison symbol records additional
dependence. We write $\mathrm C_0^\infty(G)$ for the space of smooth
compactly supported functions and $\mathscr S(G)$ for
the Schwartz space in the Euclidean coordinates on $G$.
Unless otherwise specified, functions are complex-valued.
\end{rem}

\subsection{The square-root problem and the parabolic theory}

The abstract square-root problem originates in Kato's study
of fractional powers of accretive operators
\cite{KatoFractionalI,KatoFractionalII}, see also
\cite[Chapter~VI, Remark~2.29]{KatoPerturbation}
for its formulation in terms of closed sectorial forms. For a maximal accretive operator $B$ on a Hilbert space,
Kato proved that
\[
 \mathrm D(B^\alpha)=\mathrm D((B^*)^\alpha),
 \qquad 0<\alpha<1/2.
\]
Thus, below the square-root exponent, the fractional powers
of an operator and its adjoint have the same domain.  At the square-root exponent $\alpha=1/2$,
however, this equality can fail. McIntosh's counterexample
\cite{McIntoshComparability} shows that, even for regularly
accretive operators, comparability of the two square roots
does not follow from abstract form theory alone. A positive
answer for a differential operator therefore requires
structure beyond accretivity.

For divergence-form elliptic operators with bounded measurable complex
coefficients, Auscher, Hofmann, Lacey, McIntosh and Tchamitchian
\cite{AHLMcT} proved
\[
 \mathrm D\bigl((-\operatorname{div}(A\nabla))^{1/2}\bigr)
 =\mathrm H^1(\R^n),
 \qquad
 \|(-\operatorname{div}(A\nabla))^{1/2}u\|_2\sim\|\nabla u\|_2.
\]
The development of this result is closely connected with singular
integrals and square functions. The Cauchy-integral theorem of Coifman,
McIntosh and Meyer \cite{CoifmanMcIntoshMeyer}, the second-order theory
of Auscher and Tchamitchian \cite{AuscherTchamitchian}, and the
first-order framework of Axelsson, Keith and McIntosh
\cite{AxelssonKeithMcIntosh} provide several of the analytic methods
that now underlie square-root estimates.
For a recent review of the Kato square root problem, its historical
development and its connections with elliptic and parabolic boundary
value problems, we refer to Auscher \cite{AuscherMeyerKato}.

For parabolic operators, the natural energy includes half a derivative
in time. Kaplan's whole-line variational method \cite{Kaplan} reveals
hidden coercivity through a change of test function involving the
time Hilbert transform. Nystr\"om \cite{Nystrom2016} developed a
parabolic square-function approach and proved the Kato estimate for
time-independent coefficients. Auscher, Egert and Nystr\"om
\cite{AuscherEgertNystrom2020} subsequently solved the whole-space
parabolic Kato problem for bounded complex accretive coefficients
depending measurably on both space and time. Their Theorem~2.6
shows that, for
$\mathcal H=\partial_t-\operatorname{div}_x(A(x,t)\nabla_x)$,
\[
 \mathrm D(\mathcal H^{1/2})
 =\mathrm D(\nabla_x)\cap\mathrm D(D_t^{1/2}),\qquad
 \|\mathcal H^{1/2}u\|_2
 \sim\|\nabla_xu\|_2+\|D_t^{1/2}u\|_2.
\]
Here $D_t^{1/2}$ denotes the half-order time derivative, defined
by the Fourier multiplier $|\tau|^{1/2}$ in the time variable. The work \cite{AuscherEgertNystrom2020} also develops a parabolic Dirac framework for boundary
value problems and layer potentials. Its variational formulation
and hidden-coercivity argument provide a model for
\eqref{eq:form-intro}--\eqref{eq:hidden-intro}. However, replacing
$\partial_t$ by the kinetic transport field $T$ introduces the
noncommutation relation
$[\partial_{X_j},T]=-\partial_{Y_j}$, which must be addressed
both in the coercivity argument and in the intrinsic smoothing
estimates.

Our second-order organization of the proof of Theorem~\ref{thm:main} is inspired by
Ataei, Egert and Nystr\"om \cite{AtaeiEgertNystrom}. Their work
proves the parabolic Kato estimate with a spatial Muckenhoupt
$A_2$ weight and merely measurable space-time coefficients, using
resolvents, principal-part approximation and localized affine
testing.  We extend this second-order strategy to the unweighted kinetic
setting, where transport and horizontal differentiation do not
commute. Our key contributions to this strategy are the regularized
transport symmetrizer that recovers coercivity in the full kinetic
energy space, and commutator square-function estimates that
control the interaction between horizontal differentiation and
intrinsic smoothing. Together with a dyadic approximation estimate,
these commutator bounds allow us to remove the auxiliary smoothing
before the Carleson argument and to complete the square-root
analysis for complex coefficients that depend measurably on
all variables.

\subsection{Kinetic regularity and intrinsic Sobolev spaces}

The kinetic structure goes back to Kolmogorov's explicit fundamental
solution for a degenerate diffusion equation \cite{Kolmogorov1934}.
H\"ormander's bracket criterion \cite{Hormander1967} placed this example
within a general theory of hypoelliptic differential operators.
In the present notation,  $[\partial_{X_j},T]=-\partial_{Y_j}$ expresses how diffusion in the velocity variable communicates
regularity to the position variable. Bouchut \cite{Bouchut2002}
developed quantitative Sobolev estimates for this transfer, including
equations with differentiated source terms. The orders of the
derivatives in \eqref{eq:energy-space} and \eqref{eq:y-consequence}
agree with the kinetic scaling.

The space $\V$ also has a direct interpretation in the intrinsic
Sobolev framework of Pascucci and Pesce \cite{PascucciPesce}.
Their Definitions~2.5--2.6 introduce first-order spaces through
horizontal derivatives and the flow-difference seminorm
\[
 [u]_{T,1/2,p}
 =\left(\int_{|s|\le1}
   \frac{\|e^{sT}u-u\|_p^p}{|s|^{1+p/2}}\d s\right)^{1/p}.
\]
After matching coordinates and signs, their inhomogeneous space
$\mathrm W_B^{1,2}$ is exactly $\V$, with an equivalent norm.
Indeed, the spectral theorem gives
\begin{equation}\label{eq:transport-flow-literature}
 \int_{\R}\frac{\|e^{hT}u-u\|_2^2}{|h|^2}\d h
 =2\pi\,\||T|^{1/2}u\|_2^2.
\end{equation}
The contribution from $|h|>1$ to the integral on the left-hand side in \eqref{eq:transport-flow-literature} is at most $8\|u\|_2^2$.
Restricting the differences to $|h|\le1$, as in their definition,
therefore gives an equivalent inhomogeneous norm, in agreement with
their Remark~2.8. In particular, the half-transport condition itself
has this established intrinsic interpretation. The spectral identity in \eqref{eq:transport-flow-literature} is a Hilbert-space statement. For general $p$, a flow
Slobodeckij norm should not automatically be identified with a
fractional-power norm.

The kinetic spaces of Auscher and his coauthors have a different
definition. In the local case, a basic homogeneous space in the scale
of Auscher--Niebel \cite[Section~2.1.5]{AuscherNiebelSobolev},
denoted there by $\dot{\mathcal F}^{1,2}_1$, requires
\begin{equation}\label{eq:weak-kinetic-graph}
 u\in\mathrm L^2_{t,Y}\dot{\mathrm H}^{1}_X,\qquad
 Tu\in\mathrm L^2_{t,Y}\dot{\mathrm H}^{-1}_X.
\end{equation}
These conditions are understood with their specified homogeneous
Sobolev realizations and mixed-space conventions, rather than solely
through the formal requirement $\nabla_Xu\in\mathrm L^2$.
The norm measures one transport derivative in a negative-order
velocity space.
In contrast, $\V$ requires half a transport derivative in
$\mathrm L^2(G)$ and also includes the $\mathrm L^2$ norm of $u$.
These distinctions must be retained when comparing homogeneous
and inhomogeneous spaces.

Auscher, Imbert and Niebel
\cite[Theorem~1.16 and Section~8]{AuscherImbertNiebelWeak}
prove sharp transfer of regularity and temporal trace estimates.
Under their specified distributional hypotheses, their conclusions
at the basic local energy level include
$D_Y^{1/3}u\in\mathrm L^2$ and
$u\in\mathrm C_0(\R_t;\mathrm L^2_{X,Y})$.
Here $\mathrm C_0$ denotes continuous functions of $t$ with values
in $\mathrm L^2_{X,Y}$ whose norm tends to zero as $t\to\pm\infty$.
Their weak-solution theory permits transport sources in the larger sum
\[
 \mathrm L^2_{t,Y}\dot{\mathrm H}^{-1}_X
 +\mathrm L^2_{t,X}\dot{\mathrm H}^{-1/3}_Y
 +\mathrm L^1_t\mathrm L^2_{X,Y}.
\]
Their companion paper \cite{AuscherImbertNiebelFundamental}
constructs fundamental solution operators and proves
$\mathrm L^2$ off-diagonal estimates for local diffusion.

Auscher and Niebel \cite{AuscherNiebelSobolev} extend the kinetic
Sobolev theory to $1<p<\infty$. They explicitly discuss comparison
with Pascucci--Pesce's spaces and leave half-transport regularity for
their larger class $\mathcal L^{1,p}_1$ as an open comparison problem.
This statement concerns that larger class, not an asserted equality
with $\V$. We use neither such an identification nor a general
embedding of their graph spaces into $\V$ in the present proof.

Garofalo and Tralli \cite{GarofaloTralliExtensions} study fractional
powers of constant-matrix hypoelliptic operators through degenerate
extension problems. Our Theorem~\ref{thm:main} treats measurable
complex diffusion depending on all variables and identifies the
square-root domains of the whole-space operator and its adjoint
with an established intrinsic Sobolev space. This gives an
operator-domain statement beyond regularity estimates for
smooth functions.

\subsection{Variational hidden coercivity in the kinetic setting}

The estimate \eqref{eq:hidden-intro} holds on the whole space--time
energy space. The regularized transport phase recovers the
half-transport energy without differentiating the measurable
coefficients. We here place this variational construction in the context
of other coercivity methods for kinetic equations.

Our use of hidden coercivity should be distinguished from the
modified-energy methods of kinetic hypocoercivity. Villani
\cite{VillaniHypocoercivity} develops a systematic approach
to convergence to equilibrium through suitable Lyapunov
functionals. H\'erau \cite{HerauHypocoercivity} presents
modified-energy estimates for decay and short-time
regularization of kinetic evolution equations and explicitly
describes hypocoercivity in terms of hidden coercivity.

Although these methods share the use of transport--diffusion
interaction, the purpose and form of our estimate are different.
We establish a variational coercivity estimate for the full
space--time form: testing against $(I+\delta J_\varepsilon)u$
recovers control of both $\nabla_Xu$ and $|T|^{1/2}u$.
This estimate provides the variational realization of the
operator and the energy resolvent bounds needed for the
square-root problem, rather than a decay estimate for an
evolving solution. The argument takes place on
$\mathrm L^2(G)$ with Lebesgue measure and requires no invariant
probability measure, confinement assumption or spectral gap.

Related variational approaches to kinetic equations have been
developed in several settings. Albritton, Armstrong, Mourrat
and Novack \cite{AlbrittonArmstrongMourratNovack} establish
a weak-solution theory in kinetic graph spaces using functional
inequalities and, in the conservative-force setting, uniformly
convex minimization. Litsg{\aa}rd and Nystr\"om
\cite{LitsgardNystrom2021} use a uniformly convex minimization
principle to prove existence and uniqueness for
Cauchy--Dirichlet problems in bounded Lipschitz-type cylinders,
with bounded measurable real symmetric uniformly elliptic
diffusion matrices. In these graph-space formulations,
horizontal Sobolev regularity is combined with control of
the transport derivative in a negative-order velocity Sobolev
space. Auscher, Imbert and Niebel
\cite{AuscherImbertNiebelWeak} construct weak solutions using
a Lions variational theorem, with the transport derivative
acting on the test function. Their diffusion forms also
allow complex coefficients.

These results concern variational constructions of weak
solutions. In the present setting, \eqref{eq:hidden-intro}
provides coercivity on the half-transport energy space $\V$.
This yields a direct variational construction of weak solutions
of the shifted equation in this natural energy space,
independently of the square-root estimates. Identifying $\V$
with the square-root domains requires the additional quadratic
estimates and does not follow from variational well-posedness
alone. We expect these results to support kinetic analogues of the
parabolic boundary value theories developed in
\cite{AuscherEgertNystrom2020,AuscherEgertNystrom2018}.

\subsection{The kinetic arguments and the proof strategy}

The hidden-coercivity argument described above provides the operator
realization, its adjoint and the energy resolvent estimates
(Section~\ref{sec:variational-realization}). The
remaining task is to prove two quadratic estimates, one for each
realization. We use the scale-dependent families
\[
 R_\lambda=(I+\lambda^2H_A)^{-1},\qquad
 Q_\lambda=\lambda H_AR_\lambda,
 \qquad \lambda>0.
\]
The desired estimate is
\begin{equation}\label{ThmSq}
 \int_0^\infty\|Q_\lambda u\|_2^2\frac{\d\lambda}{\lambda}
       \lesssim E(u)^2,
 \qquad u\in\V,
\end{equation}
with the corresponding assertion for $H_A^*$.

To prove these estimates, we introduce smoothing adapted to the
kinetic geometry, providing the additional regularity needed
for the principal-part approximation below. We construct this
smoothing using the positive Rockland operator
\begin{equation}\label{Rockland}
 \mathcal R:=(-\Delta_X)^2+T^\ast T=(-\Delta_X)^2-T^2
\end{equation}
on the Kolmogorov group $(G,\circ)$ equipped with the kinetic
dilations $\delta_\lambda$. We use the nonnegative self-adjoint
realization of $\mathcal R$ on $\mathrm L^2(G)$. Both terms have homogeneous degree
four under the kinetic dilations, and $-T^2=T^*T$ is nonnegative
because $T$ is skew-adjoint. Thus $\mathcal R^{1/4}$ has the same
homogeneous degree as $\nabla_X$ and $|T|^{1/2}$.
Rockland Riesz-transform estimates, interpolation and duality
make this correspondence precise
(Proposition~\ref{wr:energy-comparison}):
\[
 \mathrm D(\mathcal R^{1/4})=\V,
 \qquad
 \|\mathcal R^{1/4}u\|_2\sim E(u),
 \qquad u\in\V.
\]
This identifies the full fractional domain as well as its
homogeneous norm. The smoothing family
\[P_\lambda:=e^{-\lambda^4\mathcal R}\] then provides smooth
approximations at intrinsic scale $\lambda$ and the associated
smoothing square-function estimates
(Lemma~\ref{aux:rockland-smoothing}). 

Convolution on the Kolmogorov group provides smooth
approximations compatible with the intrinsic translations
and kinetic dilations. Our choice of the convolution family
generated by $\mathcal R$ implements an intrinsic version
of the Littlewood--Paley approach used in elliptic and
parabolic square-root arguments. The spectral theorem
supplies square-function estimates for functions of
$\mathcal R$, while the energy comparison and
Riesz-transform bounds convert these into the differential
smoothing estimates needed below. Moreover, $P_\lambda$
commutes with the fractional powers of $\mathcal R$,
a property used when removing the smoothing.
These additional properties do not follow merely from
smoothness and normalization of an intrinsic convolution
kernel. Horizontal differentiation still does not commute
with $P_\lambda$, and the resulting commutator requires
a separate estimate.

To explain how these estimates enter the quadratic estimate, we introduce
\[
 \Theta_\lambda F
 :=-\lambda R_\lambda\operatorname{div}_X(AF).
\]
Here the divergence is interpreted as an element of the continuous antidual $\V'$ of $\V$,
and the resolvent is extended variationally to this space. For $u\in\V$, the variational resolvent identity gives
the exact decomposition
\[
 Q_\lambda u
 =Q_\lambda(I-P_\lambda)u
  +\lambda R_\lambda TP_\lambda u
  +\Theta_\lambda\nabla_XP_\lambda u.
\]
The first term is the error introduced by smoothing, while the
second and third terms correspond to transport and diffusion.
The bounds
$\|Q_\lambda\|_{2\to2}\le2/\lambda$ and
$\|R_\lambda\|_{2\to2}\le1$, together with the Rockland
square-function estimates for
$\lambda^{-1}(I-P_\lambda)u$ and $\lambda TP_\lambda u$,
control the first two terms. The main remaining task is therefore
to estimate $\Theta_\lambda\nabla_XP_\lambda u$.

To obtain the required off-diagonal estimates, we localize
the resolvent equation using smooth cutoffs $\chi$ adapted
to the kinetic geometry. The resulting Caccioppoli estimate
includes contributions from both $\nabla_X\chi$ and the
full transport derivative
$T\chi=\partial_t\chi-X\cdot\nabla_Y\chi$.
Iterating this estimate and using the global resolvent bounds
yields off-diagonal estimates (Lemma~\ref{lem:off-diagonal}) of arbitrary polynomial order for $R_\lambda$,
$\Theta_\lambda$ and their scaled horizontal gradients.
These estimates allow us to define the principal coefficients
$\gamma_\lambda$ through the action of $\Theta_\lambda$ on
constant vector fields, interpreted by local $\mathrm L^2$
limits of cutoff approximations. Let $\mathsf E_\lambda$
denote intrinsic dyadic averaging at scale $\lambda$.
The off-diagonal bounds and a local kinetic Poincar\'e estimate
(Lemma~\ref{wt:poincare}) give
\[
 \|\Theta_\lambda F-\gamma_\lambda\mathsf E_\lambda F\|_2^2
 \lesssim
 \lambda^2\|\nabla_XF\|_2^2+\lambda^4\|TF\|_2^2,
\]
initially for Schwartz vector fields $F$ (Lemma~\ref{lem:local-principal-part}). The two powers of
$\lambda$ reflect the different intrinsic degrees of horizontal
differentiation and transport. Applying this estimate to
$F=\nabla_XP_\lambda u$ and using
the smoothing estimates for $\lambda\nabla_X^2P_\lambda u$
and $\lambda^2T\nabla_XP_\lambda u$, together with the bounds
for the first two terms in the decomposition of $Q_\lambda u$,
gives the smoothed principal-part approximation (Lemma~\ref{lem:smoothed-principal-part}):
\[
 \int_0^\infty
 \|Q_\lambda u
   -\gamma_\lambda\mathsf E_\lambda\nabla_XP_\lambda u\|_2^2
 \frac{\d\lambda}{\lambda}
 \lesssim E(u)^2,
 \qquad u\in\V.
\]

We next remove $P_\lambda$ from the principal part. This requires
additional estimates because
$[\partial_{X_j},T]=-\partial_{Y_j}$ and horizontal
differentiation does not commute with $P_\lambda$.
Since $\mathsf E_\lambda^2=\mathsf E_\lambda$, the error
decomposes as
\[
 \begin{split}
 \gamma_\lambda\mathsf E_\lambda
       (\nabla_XP_\lambda-\nabla_X)u
 &=\gamma_\lambda\mathsf E_\lambda
       [\nabla_X,P_\lambda]u+\gamma_\lambda\mathsf E_\lambda
       (P_\lambda-\mathsf E_\lambda)\nabla_Xu.
 \end{split}
\]
An annular cancellation argument controls the square function
of $[\nabla_X,P_\lambda]u$
(Lemma~\ref{lem:gradient-commutator}).
A separate dyadic approximation estimate
(Lemma~\ref{lem:dyadic-unsmoothing})
controls $(P_\lambda-\mathsf E_\lambda)F$ for arbitrary
$\mathrm L^2$ inputs and is applied componentwise to
$F=\nabla_Xu$. Together with the one-scale boundedness of
$\gamma_\lambda\mathsf E_\lambda$, these estimates give
(Proposition~\ref{prop:unsmoothed-principal-part})
\[
 \int_0^\infty
 \|Q_\lambda u-\gamma_\lambda\mathsf E_\lambda\nabla_Xu\|_2^2
       \frac{\d\lambda}{\lambda}
 \lesssim E(u)^2,\qquad u\in\V.
\]
This approximation is established independently of the Carleson
estimate.

It therefore remains to control the principal term
$\gamma_\lambda\mathsf E_\lambda\nabla_Xu$. To accomplish this we construct localized affine resolvent tests that satisfy average accretivity
and square-function bounds. A stopping-time argument then proves
that $|\gamma_\lambda(z)|^2\d z\,\d\lambda/\lambda$ is a
Carleson measure (Theorem~\ref{wt:carleson}).
The dyadic Carleson embedding theorem bounds
the square function of the principal term by
$\|\nabla_Xu\|_2^2$, completing the quadratic estimate for $H_A$.
Repeating the argument with $A^*$ and $-T$ gives the adjoint
estimate, using the same Rockland smoothing operator.

Finally, functional calculus gives the upper square-root bounds.
Form duality and hidden coercivity give the lower bounds, and
resolvent approximation proves equality of the square-root
domains with $\V$, as shown in the proof of
Proposition~\ref{prop:quadratic-reduction}.

\subsection{Organization of the paper}

Section~\ref{sec:geometry-energy} introduces the geometry and the
transport energy and proves the uncertainty estimate.
Section~\ref{sec:variational-realization} constructs the operators
and their energy resolvents. Section~\ref{sec:quadratic-reduction}
reduces the theorem to the two quadratic estimates.
Section~\ref{sec:flat-rockland} proves the Rockland comparison and
smoothing estimates. Sections~\ref{sec:localization} and
\ref{sec:unsmoothing} establish the principal-part approximation.
Section~\ref{sec:cm} proves the Carleson estimate and completes
the proof of Theorem~\ref{thm:main}. The shifted estimates and
Riesz transforms are recorded in Section~\ref{sec:consequences}.
Appendix~\ref{app:common-core} proves
Lemma~\ref{wr:common-core} from Section~\ref{sec:flat-rockland},
which establishes the invariance of $\mathscr S_0(G)$ under
complex powers of $\mathcal R$ and its common-core properties.
Appendix~\ref{app:type-zero-commutator} provides the proof
of the type-zero commutator estimate stated in
Lemma~\ref{lem:type-zero-commutator}
in Section~\ref{sec:unsmoothing}.

\section{Geometry and the transport energy}
\label{sec:geometry-energy}

We describe the homogeneous geometry underlying the Kolmogorov
operator and the spectral realization of the transport field.
We then establish the basic properties of the energy space $\V$,
including density of smooth functions and the uncertainty estimate
that controls fractional regularity in the position variables.

\subsection{The Kolmogorov group and its geometry}
\label{subsec:kolmogorov-geometry}

We work on $G=\R^m_X\times\R^m_Y\times\R_t$ with Lebesgue
measure $\d z=\d X\d Y\d t$. All inner products are linear in
the first argument. We write $\langle\cdot,\cdot\rangle$ and
$\|\cdot\|_2$ for the scalar or finite-dimensional vector-valued
$\mathrm L^2(G)$ inner product and norm. We use the comparison
notation and the convention on implicit constants from
Remark~\ref{convention}.

The group law and dilations are
\begin{align}
 (X,Y,t)\circ(X',Y',t')
  &=(X+X',Y+Y'-t'X,t+t'),\label{eq:group-law}\\
 \delta_\lambda(X,Y,t)
  &=(\lambda X,\lambda^3Y,\lambda^2t),\qquad\lambda>0.
       \label{eq:dilation}
\end{align}
Lebesgue measure is the Haar measure, with homogeneous dimension
$\mathfrak q=4m+2$. The left-invariant fields $\partial_{X_j}$
and $T=\partial_t-X\cdot\nabla_Y$ have degrees one and two, and
\begin{equation}\label{eq:bracket}
 [\partial_{X_j},T]=-\partial_{Y_j},\qquad
 T\partial_{X_j}=\partial_{X_j}T+\partial_{Y_j}.
\end{equation}
We fix a symmetric homogeneous quasi-norm $\rho$, smooth away from
the identity $0$. Thus $\rho\colon G\to[0,\infty)$ is continuous,
vanishes only at $0$, and satisfies
\[
 \rho(z^{-1})=\rho(z),\qquad
 \rho(\delta_\lambda z)=\lambda\rho(z),\qquad
 \rho(z\circ w)\lesssim\rho(z)+\rho(w),
 \qquad z,w\in G,\quad\lambda>0.
\]
For example, one may take
\[
 \rho(X,Y,t)
 :=\bigl(|X|^{12}
       +\bigl|Y+\tfrac12tX\bigr|^4+t^6\bigr)^{1/12}.
\]
The centred coordinate $Y+\tfrac12tX$ changes sign under
group inversion, which gives the required symmetry.
This quasi-norm satisfies
\[
 \rho(X,Y,t)\sim |X|+|Y|^{1/3}+|t|^{1/2},
\]
making the intrinsic scales explicit.   Hence
\begin{equation}\label{eq:kinetic-distance}
d(z,z'):=\rho(z^{-1}\circ z')
\end{equation}
is a symmetric, left-invariant
quasi-metric. Its balls
$B(z,\lambda):=\{z'\in G:d(z,z')<\lambda\}$ satisfy
\begin{equation}\label{eq:kinetic-volume}
 |B(z,\lambda)|=|B(0,1)|\lambda^{\mathfrak q}
\end{equation}
by left invariance of Lebesgue measure and homogeneity.
In particular, Lebesgue measure is doubling, and
$(G,d,\d z)$ is a space of homogeneous type.

\subsection{The closed derivatives and the energy space}

We take $\nabla_X$ with its maximal distributional domain in
$\mathrm L^2(G)$. This is a closed operator, since convergence
in $\mathrm L^2(G)$ implies convergence of distributions.
The maps
\begin{equation}\label{eq:transport-flow}
 (\Phi_su)(X,Y,t)=u(X,Y-sX,t+s),\qquad s\in\R,
\end{equation}
form a strongly continuous unitary group on $\mathrm L^2(G)$.
Indeed, their Jacobians are one, and strong continuity follows
first for compactly supported continuous functions and then by
density. Their infinitesimal generator is the skew-adjoint realization
of $T$. Consequently, $\Phi_s=e^{sT}$, where the exponential
is defined by the spectral calculus for the self-adjoint
operator $-iT$.

We use the unitary partial Euclidean Fourier transform in
$(Y,t)$, leaving $X$ unchanged. Define
\[
 \widehat u(X,\eta,\tau)
 :=(2\pi)^{-(m+1)/2}
   \int_{\R^m\times\R}
   e^{-i(Y\cdot\eta+t\tau)}u(X,Y,t)\d Y\d t,\qquad u\in\mathscr S(G).
\]
By Plancherel's theorem, this transform extends uniquely
to a unitary operator on $\mathrm L^2(G)$.
Throughout the paper, we use the notation
\[
 q(X,\eta,\tau):=\tau-X\cdot\eta,
 \qquad a(\eta):=|\eta|^{2/3}.
\]
The weight $a$ has the same intrinsic scaling as $q$ and
will enter both the uncertainty estimate and the construction
of the regularized transport multiplier below. Under the partial Fourier transform, $\Phi_s$ is represented
by multiplication by $e^{isq}$, and $T$ is unitarily equivalent
to multiplication by $iq$, with maximal domain
\[
 \mathrm D(T)
 =\{u\in\mathrm L^2(G):
 q\widehat u\in
 \mathrm L^2(\R^m_X\times\R^m_\eta\times\R_\tau)\}.
\]
Consequently, $T$ is skew-adjoint, and its nonnegative modulus
is defined by $|T|:=(T^*T)^{1/2}=(-T^2)^{1/2}$. The spectral theorem then identifies $|T|^{1/2}$, under
the same partial Fourier transform, with multiplication
by $|q|^{1/2}$.  We use the energy space
\begin{equation}\label{eq:energy-space-body}
 \V=\mathrm D(\nabla_X)\cap\mathrm D(|T|^{1/2}),\qquad
 E(u)=\|\nabla_Xu\|_2+\||T|^{1/2}u\|_2.
\end{equation}
Closedness of the operators $\nabla_X$ and $|T|^{1/2}$
shows that $\V$ is complete for the Hilbert norm
$\|\cdot\|_\V$ defined in the introduction. It embeds densely and continuously in $\mathrm L^2(G)$. The scalar identity
\[
 \int_\R\frac{|e^{isq}-1|^2}{s^2}\d s=2\pi|q|
\]
follows from one-dimensional Plancherel applied to the
indicator of an interval of length $|q|$. The spectral theorem and this identity give
\begin{equation}\label{eq:transport-differences}
 2\pi\||T|^{1/2}u\|_2^2=\int_\R\frac{\|\Phi_su-u\|_2^2}{|s|^2}\d s.
\end{equation}
This equality holds with extended values for every
$u\in\mathrm L^2(G)$. 

The transport form is
\begin{equation}\label{eq:transport-form}
 \mathfrak t(u,v)
 =\int_{\R^{2m+1}}
 i\sgn(q)|q|^{1/2}\widehat u\,
       \overline{|q|^{1/2}\widehat v}\,
 \d X\d\eta\d\tau .
\end{equation}
It is bounded and skew-Hermitian on $\V$. For
$u\in\V\cap\mathrm D(T)$ and $v\in\V$, it equals
$\langle Tu,v\rangle$.
For $u\in\V$ and $\varphi\in\mathrm C_0^\infty(G)$ it equals
$-\langle u,T\varphi\rangle$. $\sgn$ denotes the sign function, with the convention
$\sgn(0)=0$.

\begin{lem}\label{lem:form-density}
The classes $\mathrm C_0^\infty(G)$ and $\mathscr S(G)$ are
contained in $\V$ and dense in its Hilbert norm.
\end{lem}

\begin{proof}
For $v\in\mathscr S(G)$, both $\nabla_Xv$ and $Tv$ belong
to $\mathrm L^2(G)$. The spectral theorem gives
\[
 \||T|^{1/2}v\|_2^2\le\|Tv\|_2\|v\|_2.
\]
Thus $\mathrm C_0^\infty(G)\subset\mathscr S(G)\subset\V$.

We now prove density, beginning with approximation by
Schwartz functions. Fix $u\in\V$ and write $U=\widehat u$.
Choose cutoffs $\chi\in\mathrm C_0^\infty(\R^m)$ and
$\vartheta\in\mathrm C_0^\infty(\R^{m+1})$, taking values
in $[0,1]$ and equal to one near the origin. Set
\[
 \chi_n(X):=\chi(X/n),\qquad
 \vartheta_n(\eta,\tau):=\vartheta((\eta,\tau)/n),
 \qquad U_n:=\chi_n\vartheta_nU.
\]
Both cutoffs commute with multiplication by $|q|^{1/2}$,
and $\vartheta_n$ commutes with $\nabla_X$. Moreover,
\[
 \nabla_XU_n
 =\chi_n\vartheta_n\nabla_XU
   +(\nabla_X\chi_n)\vartheta_nU,
 \qquad
 \|(\nabla_X\chi_n)\vartheta_nU\|_2
 \lesssim n^{-1}\|U\|_2.
\]
Since $\chi_n\vartheta_n\to1$ pointwise, dominated
convergence and the preceding bound give
\[
 \|U_n-U\|_2
 +\|\nabla_X(U_n-U)\|_2
 +\||q|^{1/2}(U_n-U)\|_2
 \longrightarrow0.
\]
Each $U_n$ has compact support in $(X,\eta,\tau)$. Next, for each fixed $n$, mollify $U_n$ using a standard Euclidean
approximate identity in these variables, and denote the result
by $U_{n,\varepsilon}$. Then
$U_{n,\varepsilon}\in\mathrm C_0^\infty(\R^{2m+1})$.
Mollification converges in $\mathrm L^2$ and commutes with
distributional $X$ derivatives, so
\[
 \|U_{n,\varepsilon}-U_n\|_2
 +\|\nabla_X(U_{n,\varepsilon}-U_n)\|_2
 \longrightarrow0
 \qquad\text{as }\varepsilon\downarrow0.
\]
For $0<\varepsilon\le1$, all these functions are supported
in a fixed compact enlargement $K_n^*$ of the support of
$U_n$. Since $q=\tau-X\cdot\eta$ is bounded there,
\[
 \||q|^{1/2}(U_{n,\varepsilon}-U_n)\|_2
 \le
 \left(\sup_{K_n^*}|q|\right)^{1/2}
 \|U_{n,\varepsilon}-U_n\|_2
 \longrightarrow0.
\]
The bound may depend on $n$, which is harmless since $n$
is fixed during this approximation. Finally, choose $0<\varepsilon_n<1/n$ so that the sum of these
three mollification errors is less than $1/n$.
The inverse partial Fourier transforms of
$U_{n,\varepsilon_n}$ are Schwartz functions with compact
$X$ support, by smoothness and integration by parts in
the Fourier variables. Partial Plancherel and the preceding
estimates show that they converge to $u$ in $\V$.

It remains to approximate Schwartz functions by compactly
supported smooth functions in the $\V$ norm.
Fix $v\in\mathscr S(G)$ and choose
$\chi_\lambda\in\mathrm C_0^\infty(G)$ tending pointwise
to one as $\lambda\to\infty$, with
\[
 |\chi_\lambda|\le1,\qquad
 |\nabla_X\chi_\lambda|\lesssim\lambda^{-1},\qquad
 |T\chi_\lambda|\lesssim\lambda^{-2}.
\]
For example, take products of cutoffs in
$X/\lambda,Y/\lambda^3,t/\lambda^2$. Set $w_\lambda:=(\chi_\lambda-1)v$.
Dominated convergence gives $\|w_\lambda\|_2\to0$.
Moreover, the product rules
 \begin{align*}
 \nabla_Xw_\lambda
 &=(\chi_\lambda-1)\nabla_Xv
   +(\nabla_X\chi_\lambda)v,\quad Tw_\lambda
 =(\chi_\lambda-1)Tv+(T\chi_\lambda)v
 \end{align*}
and the cutoff bounds imply
\[
 \|\nabla_Xw_\lambda\|_2+\|Tw_\lambda\|_2
 \longrightarrow0,
\]
since $\nabla_Xv,Tv\in\mathrm L^2(G)$.
As $w_\lambda\in\mathscr S(G)\subset\mathrm D(T)$,
spectral Cauchy--Schwarz gives
\[
 \||T|^{1/2}w_\lambda\|_2^2
 \le\|Tw_\lambda\|_2\|w_\lambda\|_2
 \longrightarrow0.
\]
Thus $\chi_\lambda v\to v$ in $\V$. Finally, let $v_n\in\mathscr S(G)$ converge to $u$ in $\V$.
For each $n$, choose $\lambda_n$ sufficiently large that
\[
 \|\chi_{\lambda_n}v_n-v_n\|_\V<1/n.
\]
Then $\chi_{\lambda_n}v_n\in\mathrm C_0^\infty(G)$ and
\[
 \|\chi_{\lambda_n}v_n-u\|_\V
 \le 1/n+\|v_n-u\|_\V\longrightarrow0,
\]
proving the assertion.\end{proof}

\subsection{The uncertainty inequality}

\begin{lem}\label{lem:uncertainty}
For every $u\in\V$,
\begin{equation}\label{eq:uncertainty}
 \|D_Y^{1/3}u\|_2^2
 \lesssim\|\nabla_Xu\|_2^2+\||T|^{1/2}u\|_2^2,
 \qquad D_Y^{1/3}=(-\Delta_Y)^{1/6}.
\end{equation}
The implicit constant depends only on $m$.
\end{lem}

\begin{proof}
Fix $\eta\ne0$ and $\tau$, and put
$b=|\eta|$, $e=\eta/b$, $\lambda=b^{-1/3}$, and $a=b^{2/3}=\lambda^{-2}=|\eta|^{2/3}$. Recall the notation $q=\tau-X\cdot\eta$. Tile $\R^m_X$ by axis-parallel cubes of side $\lambda$.
For each such cube $Q$ let
\[
 E_Q=\{X\in Q:|\tau-X\cdot\eta|\ge c_m a\},
 \qquad c_m=(4\sqrt m)^{-1}.
\]
Choose $j$ with $|e_j|\ge m^{-1/2}$, which is possible
because $|e|=1$. For fixed values of the other coordinates,
\[
 q=\tau-b\sum_{k\ne j}e_kX_k-be_jX_j
\]
is an affine function of $X_j$ with slope $-be_j$.
Thus the excluded condition $|q|<c_m a$ restricts $X_j$
to an interval centered at $(\tau-b\sum_{k\ne j}e_kX_k)/(be_j)$  and of length
\[
 \frac{2c_m a}{b|e_j|}
 =\frac{2c_m\lambda}{|e_j|}
 \le 2c_m\sqrt m\,\lambda
 =\frac{\lambda}{2},
\]
where we used $a/b=\lambda$ and $c_m=(4\sqrt m)^{-1}$. Intersecting this interval with the corresponding side of
$Q$ can only decrease its length. Since that side has
length $\lambda$, the portion of each slice belonging
to $E_Q$ has length at least $\lambda/2$. Indeed, write $Q=I_j\times Q'_j$, with $Q'_j$ denoting the
product of the remaining coordinate intervals. Then Fubini's
theorem gives
\[
 |E_Q|
 =\int_{Q'_j}
   \bigl|\{x_j\in I_j:(x_j,X')\in E_Q\}\bigr|\,\d X'
 \ge \frac{\lambda}{2}|Q'_j|
 =\frac{|Q|}{2}.
\]
For $m=1$, this is simply the preceding interval estimate.
The bound is uniform in $\eta$, $\tau$ and $Q$.

For a smooth function $f$ on $\R^m$, set
\[
 c_Q:=\fint_Qf\d X,
 \qquad V_Q:=\int_Q|f-c_Q|^2\d X.
\]
Writing $c_Q=f(X)-(f(X)-c_Q)$ gives
\[
 |c_Q|^2\le2|f(X)|^2+2|f(X)-c_Q|^2.
\]
Integrating over $E_Q$ and using
$|E_Q|\ge |Q|/2$ and $E_Q\subset Q$, we obtain
\[
 \begin{split}
 |Q||c_Q|^2
 &\le2|E_Q||c_Q|^2\le4\biggl(\int_{E_Q}|f|^2\d X+V_Q\biggr).
 \end{split}
\]
The Euclidean Poincar\'e inequality gives
\[V_Q\lesssim\lambda^2\int_Q|\nabla_Xf|^2\d X.\]
Using this and the lower bound on $|q|$ on $E_Q$, we obtain
\[
 \int_Q|f|^2\d X
 \lesssim\lambda^2\int_Q|\nabla_Xf|^2\d X
       +a^{-1}\int_Q|q||f|^2\d X.
\]
Multiply by $a=\lambda^{-2}$ and sum over the cubes to obtain
\begin{equation}\label{hsuu}
 a\int_{\R^m}|f|^2\d X
 \lesssim\int_{\R^m}
       \bigl(|\nabla_Xf|^2+|q||f|^2\bigr)\d X.
\end{equation}
For $u\in\mathscr S(G)$, now apply this estimate with
$f(X)=\widehat u(X,\eta,\tau)$. Since $a=|\eta|^{2/3}$
and since the implicit constant is independent of $(\eta,\tau)$,
we may integrate in these variables. The case $\eta=0$
causes no difficulty, since the left-hand side then vanishes.
Partial Fourier transformation commutes with $\nabla_X$, and
\[
 \widehat{D_Y^{1/3}u}=|\eta|^{1/3}\widehat u,
 \qquad
 \widehat{|T|^{1/2}u}
       =|\tau-X\cdot\eta|^{1/2}\widehat u.
\]
Thus partial Plancherel and \eqref{hsuu} give
\[
 \|D_Y^{1/3}u\|_2^2
 \lesssim\|\nabla_Xu\|_2^2+\||T|^{1/2}u\|_2^2,
 \qquad u\in\mathscr S(G).
\]
To extend the estimate to $u\in\V$, use
Lemma~\ref{lem:form-density} to choose $u_n\in\mathscr S(G)$
with $u_n\to u$ in $\V$. Applying the estimate to $u_n-u_k$
gives
\[
 \|D_Y^{1/3}u_n-D_Y^{1/3}u_k\|_2^2
 \lesssim
 \|\nabla_X(u_n-u_k)\|_2^2
 +\||T|^{1/2}(u_n-u_k)\|_2^2
 \longrightarrow0
\]
as $n,k\to\infty$. Hence $D_Y^{1/3}u_n$ converges in
$\mathrm L^2(G)$ to some $g$. Since $u_n\to u$ in
$\mathrm L^2(G)$ and $D_Y^{1/3}$ is closed, we obtain
$u\in\mathrm D(D_Y^{1/3})$ and $D_Y^{1/3}u=g$.
Passing to the limit in the estimate for $u_n$ proves
\eqref{eq:uncertainty}.
\end{proof}

\begin{rem}
The term \emph{uncertainty estimate} expresses a competition
between concentration near the characteristic hyperplane
and horizontal differentiation. For fixed $\eta\ne0$ and
$\tau$, making $\int_{\R^m}|q||f|^2\d X$ small relative
to $\|f\|_2^2$ forces concentration near
$q=\tau-X\cdot\eta=0$. Rescaling a fixed function in the normal direction to
width $\lambda$, while keeping its $\mathrm L^2$ norm
unchanged, multiplies the $\mathrm L^2$ norm of its
normal derivative by $\lambda^{-1}$.  Its contribution to the horizontal
energy therefore scales as $\lambda^{-2}$, whereas the
transport energy scales as $|\eta|\lambda$.
Balancing these scales gives $\lambda=|\eta|^{-1/3}$
and the energy scale $|\eta|^{2/3}$ in \eqref{hsuu}.
Partial Plancherel converts this into the control of
a position derivative of order $1/3$ in
\eqref{eq:uncertainty}.
\end{rem}

\section{Hidden coercivity and variational realization}
\label{sec:variational-realization}

Assume that $A\colon G\to\C^{m\times m}$ is measurable and
satisfies \eqref{eq:ellipticity}. Set $J_0=T|T|^{-1}$, with
multiplier zero at $q=0$, and define
\begin{equation}\label{eq:form-body}
 \mathfrak a_A(u,v)
 =\langle A\nabla_Xu,\nabla_Xv\rangle+\mathfrak t(u,v),
 \qquad u,v\in\V,
\end{equation}
where $\mathfrak t$ is the transport form defined in
\eqref{eq:transport-form}. This form is bounded on $\V$, with
\[
 \operatorname{Re}\mathfrak a_A(u,u)
 \ge a_0\|\nabla_Xu\|_2^2.
\]
To control the half-transport derivative as well, we modify
the test function as follows.

\begin{lem}\label{lem:kinetic-hidden-coercivity}
There is a boundedly invertible operator $S_+$ on $\V$, also
bounded on $\mathrm L^2(G)$, such that, for every $u\in\V$,
\begin{align}
 E(S_+u)&\lesssim E(u),\label{eq:symmetrizer-energy}\\
 \operatorname{Re}\mathfrak a_A(u,S_+u)
   &\gtrsim E(u)^2,\label{eq:kinetic-hidden-coercivity}\\
 \operatorname{Re}\langle u,S_+u\rangle
   &=\|u\|_2^2.\label{eq:symmetrizer-mass}
\end{align}
The implicit constants depend only on $m,a_0,a_1$.
For the adjoint form
$\mathfrak a_A^*(u,v)=\overline{\mathfrak a_A(v,u)}$,
the same assertions hold with an operator $S_-$.
\end{lem}

\begin{proof}
For $\varepsilon>0$, define
\[
 \widehat{J_\varepsilon u}
 =m_\varepsilon\widehat u,
 \qquad
 m_\varepsilon:=\frac{iq}{(q^2+\varepsilon^2a^2)^{1/2}},
 \qquad a:=|\eta|^{2/3},\qquad q:=\tau-X\cdot\eta,
\]
and set $m_\varepsilon=0$ when $q=0$. Since $|m_\varepsilon|\le1$ and
$\overline{m_\varepsilon}=-m_\varepsilon$, partial Plancherel
shows that $J_\varepsilon$ is a skew-adjoint contraction on
$\mathrm L^2(G)$. Moreover, $|T|^{1/2}$ has multiplier
$|q|^{1/2}$. The bound $|m_\varepsilon|\le1$ therefore ensures
that $J_\varepsilon$ preserves $\mathrm D(|T|^{1/2})$, and
commutation of the multipliers gives
\[
 |T|^{1/2}J_\varepsilon u
 =J_\varepsilon|T|^{1/2}u,
 \qquad u\in\mathrm D(|T|^{1/2}).
\]
Since $q=\tau-X\cdot\eta$, differentiation gives,
for $\eta\ne0$,
\[
 \nabla_Xm_\varepsilon
 =-\frac{i\varepsilon^2a^2\eta}
         {(q^2+\varepsilon^2a^2)^{3/2}},
 \qquad
 |\nabla_Xm_\varepsilon|
 \le\varepsilon^{-1}|\eta|^{1/3}.
\]
The Fourier-side product rule therefore yields
\[
 \widehat{\nabla_X(J_\varepsilon u)}
 =m_\varepsilon\widehat{\nabla_Xu}
   +(\nabla_Xm_\varepsilon)\widehat u.
\]
Applying partial Plancherel and the preceding multiplier bounds,
we obtain
\begin{equation}\label{aapa}
 \|\nabla_XJ_\varepsilon u\|_2
 \le\|\nabla_Xu\|_2
       +\varepsilon^{-1}\|D_Y^{1/3}u\|_2.
\end{equation}
This calculation is first justified for Schwartz functions
whose partial Fourier transforms are supported away from
$\eta=0$. To approximate an arbitrary $u\in\V$ by such
functions, choose $\chi\in\mathrm C_0^\infty(\R^m)$ with
$0\le\chi\le1$ and $\chi=1$ near the origin, and define
the Fourier multipliers $M_n$ by
\[
 \widehat{M_nw}(X,\eta,\tau)
 :=\bigl(1-\chi(n\eta)\bigr)\widehat w(X,\eta,\tau).
\]
These multipliers commute with $\nabla_X$ and
$|T|^{1/2}$ and are contractions on $\V$.
Since $1-\chi(n\eta)\to1$ for $\eta\ne0$, dominated
convergence gives $M_nw\to w$ in $\V$ for every $w\in\V$. By Lemma~\ref{lem:form-density}, choose
$v_n\in\mathscr S(G)$ with $v_n\to u$ in $\V$, and set
$u_n:=M_nv_n$. For each fixed $n$, the multiplier
$1-\chi(n\eta)$ is smooth with bounded derivatives of
every order, so $u_n\in\mathscr S(G)$. Moreover, $\widehat u_n$ vanishes in a neighbourhood of
$\eta=0$, so $u_n$ belongs to the class on which the calculation
is justified. Contractivity of $M_n$ on $\V$ gives
\begin{align*}
 \|u_n-u\|_\V
 &\le\|M_n(v_n-u)\|_\V+\|M_nu-u\|_\V\le\|v_n-u\|_\V+\|M_nu-u\|_\V
 \longrightarrow0.
 \end{align*}
For fixed $\varepsilon>0$, \eqref{aapa} and the kinetic
uncertainty estimate of Lemma~\ref{lem:uncertainty} imply
\[
 \|\nabla_XJ_\varepsilon(u_n-u_k)\|_2
 \lesssim_\varepsilon E(u_n-u_k).
\]
Thus $(\nabla_XJ_\varepsilon u_n)$ is Cauchy in
$\mathrm L^2(G;\C^m)$. Boundedness of $J_\varepsilon$ on
$\mathrm L^2(G)$ also gives
$J_\varepsilon u_n\to J_\varepsilon u$ in $\mathrm L^2(G)$.
Closedness of $\nabla_X$ therefore implies
$J_\varepsilon u\in\mathrm D(\nabla_X)$ and
\[
 \nabla_XJ_\varepsilon u_n
 \longrightarrow\nabla_XJ_\varepsilon u
 \quad\text{in }\mathrm L^2(G;\C^m).
\]
Furthermore, the uncertainty estimate gives
\[
 \|D_Y^{1/3}(u_n-u)\|_2
 \lesssim E(u_n-u)\longrightarrow0.
\]
Since $\nabla_Xu_n\to\nabla_Xu$ as well, both terms on the
right-hand side of \eqref{aapa} converge. Passing to the limit
in that estimate, applied to $u_n$, proves \eqref{aapa}
for every $u\in\V$. Finally, commutation with $|T|^{1/2}$ and
$\mathrm L^2$ contractivity give
\[
 \||T|^{1/2}J_\varepsilon u\|_2
 =\|J_\varepsilon|T|^{1/2}u\|_2
 \le\||T|^{1/2}u\|_2.
\]
Combining this with \eqref{aapa} and \eqref{eq:uncertainty}
yields
\[
 E(J_\varepsilon u)\lesssim_\varepsilon E(u),
 \qquad u\in\V.
\]
Together with $\mathrm L^2$ contractivity, this proves
boundedness of $J_\varepsilon$ on $\V$.

To quantify the coercivity estimate, set
\[
 g:=\|\nabla_Xu\|_2,
 \qquad
 h:=\||T|^{1/2}u\|_2,
 \qquad
 d_u:=\|D_Y^{1/3}u\|_2.
\]
The elementary inequality
\begin{equation}\label{Elementary}
 \frac{q^2}{(q^2+\varepsilon^2a^2)^{1/2}}
 \ge |q|-\varepsilon a
\end{equation}
controls the loss introduced by regularizing the transport
multiplier. Here and below, the quotient is defined to be zero
when $q=0$.  For the diffusion contribution, ellipticity and the bound for $\nabla_XJ_\varepsilon u$ in \eqref{aapa} give
\[
 \begin{split}
 \operatorname{Re}\langle
 A\nabla_Xu,\nabla_X((I+\delta J_\varepsilon)u)\rangle
 &\ge a_0g^2-\delta a_1g\|\nabla_XJ_\varepsilon u\|_2\ge a_0g^2-\delta a_1g^2
             -\delta a_1\varepsilon^{-1}g d_u.
 \end{split}
\]
For the transport contribution, recall that
$\operatorname{Re}\mathfrak t(u,u)=0$. Partial Plancherel, the definition of $m_\varepsilon$, and \eqref{Elementary} therefore yield
\[
 \begin{split}
 \operatorname{Re}\mathfrak t(u,(I+\delta J_\varepsilon)u)
 &=\delta\int_{\R^{2m+1}}
       \frac{q^2}{(q^2+\varepsilon^2a^2)^{1/2}}
       |\widehat u|^2\d X\d\eta\d\tau\\
 &\ge\delta\int_{\R^{2m+1}}
       (|q|-\varepsilon a)|\widehat u|^2
       \d X\d\eta\d\tau=\delta h^2-\delta\varepsilon d_u^2.
 \end{split}
\]
Here the last equality uses the Fourier multipliers
$|q|^{1/2}$ and $a^{1/2}=|\eta|^{1/3}$ of
$|T|^{1/2}$ and $D_Y^{1/3}$, respectively.
Adding the two contributions in the last two displays gives
\begin{equation}\label{corbound}
 \operatorname{Re}\mathfrak a_A(u,(I+\delta J_\varepsilon)u)
 \ge a_0g^2-\delta a_1g^2
       -\delta a_1\varepsilon^{-1}g d_u
       +\delta h^2-\delta\varepsilon d_u^2.
\end{equation}
We now choose $\varepsilon$ and $\delta$ so that the negative
terms can be absorbed. By \eqref{eq:uncertainty}, we may fix $C_0=C_0(m)>0$ such that
\begin{equation}\label{corbound+uu}
 d_u^2\le C_0(g^2+h^2).
\end{equation}
We first fix $\varepsilon>0$ sufficiently small that
$C_0\varepsilon\le1/8$. By Young's inequality,
\begin{equation}\label{corbound+uu+}
 \delta a_1\varepsilon^{-1}g d_u
 \le\frac{a_0}{4}g^2
       +\frac{\delta^2a_1^2}{a_0\varepsilon^2}d_u^2.
\end{equation}
Writing $K:=C_0a_1^2/(a_0\varepsilon^2)$ and using \eqref{corbound+uu+} and \eqref{corbound+uu} to bound both remaining terms
containing $d_u^2$, we obtain from \eqref{corbound} that
\begin{equation}\label{corbound+}
 \operatorname{Re}\mathfrak a_A(u,(I+\delta J_\varepsilon)u)
 \ge
 \biggl(\frac{3a_0}{4}-\delta a_1
       -\frac{\delta}{8}-\delta^2K\biggr)g^2+
 \biggl(\frac{7\delta}{8}-\delta^2K\biggr)h^2.
 \end{equation}
Now choose $\delta>0$ sufficiently small so that
\begin{equation}\label{Invert}
 \delta K\le\frac18,
 \qquad
 \delta(a_1+1/4)\le\frac{a_0}{4},
 \qquad
 \delta\|J_\varepsilon\|_{\V\to\V}<\frac12.
\end{equation}
The coefficients of $g^2$ and $h^2$ in \eqref{corbound+} are then at least
$a_0/2$ and $3\delta/4$, respectively. Since
$E(u)^2=(g+h)^2\le2(g^2+h^2)$, setting
$S_+:=I+\delta J_\varepsilon$ we see that \eqref{corbound+} gives
\[
 \operatorname{Re}\mathfrak a_A(u,S_+u)
 \ge\frac{a_0}{2}g^2+\frac{\delta}{2}h^2
 \ge c E(u)^2,
 \qquad
 c:=\frac14\min\{a_0,\delta\}>0.
\]
All choices depend only on $m,a_0,a_1$.

The last condition on $\delta$ in \eqref{Invert} ensures that $S_+$ is boundedly
invertible on $\V$ by a Neumann series. Moreover, skew-adjointness
of $J_\varepsilon$ gives
\[
 \operatorname{Re}\langle u,S_+u\rangle
 =\|u\|_2^2+\delta\operatorname{Re}\langle u,J_\varepsilon u\rangle
 =\|u\|_2^2.
\]
For the adjoint form, $A^*$ satisfies the same coefficient bounds
and the transport term has the opposite sign. The same argument
with $S_-:=I-\delta J_\varepsilon$ gives the corresponding
coercivity estimate and mass identity. The map $S_-$ is also
boundedly invertible on $\V$.
\end{proof}

\begin{prop}\label{prop:kinetic-realization}
Define the $\mathrm L^2$ realization $H_A$ of $\mathfrak a_A$ by
\[
 \begin{aligned}
 \mathrm D(H_A):=\bigl\{u\in\V:\ &
 \text{there exists }f\in\mathrm L^2(G)\text{ such that }\mathfrak a_A(u,v)=\langle f,v\rangle
 \text{ for every }v\in\V\bigr\},
 \end{aligned}
\]
and set $H_Au:=f$. The representing function $f$ is unique.
The operator $H_A$ is densely defined, closed and maximal
accretive on $\mathrm L^2(G)$, and $\mathrm D(H_A)$ is dense
in $\V$ with respect to $\|\cdot\|_\V$.  Furthermore, its domain and action admit the following distributional
characterization: for $u\in\V$ and $f\in\mathrm L^2(G)$,
one has $u\in\mathrm D(H_A)$ and $H_Au=f$ if and only if
\begin{equation}\label{eq:distribution-realization}
 Tu-\operatorname{div}_X(A\nabla_Xu)=f
 \quad\text{in }\mathrm D'(G).
\end{equation}
The Hilbert-space adjoint $H_A^*$ is the $\mathrm L^2$
realization of the adjoint form $\mathfrak a_A^*$ on $\V$.
It realizes the differential expression
\[
 H_A^*=-T-\operatorname{div}_X(A^*\nabla_X),
\]
with the analogous distributional characterization of its
domain and action.
\end{prop}

\begin{proof}
We use the Gelfand triple
\[
 \V\hookrightarrow\mathrm L^2(G)\hookrightarrow\V',
\]
where $\V'$ denotes the continuous antidual of $\V$ and
$\mathrm L^2(G)$ is the pivot space. Thus a function
$f\in\mathrm L^2(G)$ is identified with the functional
$v\mapsto\langle f,v\rangle$ on $\V$.
The first embedding is continuous. Since $\mathrm C_0^\infty(G)\subset\V$ and
$\mathrm C_0^\infty(G)$ is dense in $\mathrm L^2(G)$,
the space $\V$ is dense in $\mathrm L^2(G)$. This ensures
uniqueness of the representing function in the definition
of $H_A$. The second embedding also has dense range, since,
under the Hilbert-space identifications, its adjoint is the
injective embedding $\V\hookrightarrow\mathrm L^2(G)$.

For $\nu>0$, define the bounded operator
\[
 B_\nu\colon\V\to\V',
 \qquad
 (B_\nu u)(v):=\nu\langle u,v\rangle+\mathfrak a_A(u,v).
\]
We first show that $B_\nu$ is an isomorphism.
Lemma~\ref{lem:kinetic-hidden-coercivity} provides a boundedly
invertible map $S_+$ on $\V$ such that
\[
 \operatorname{Re}(B_\nu u)(S_+u)
 \ge\nu\|u\|_2^2+cE(u)^2
 \ge\min\{\nu,c\}\|u\|_\V^2.
\]
Consequently, the bounded sesquilinear form
$(u,v)\mapsto(B_\nu u)(S_+v)$ is coercive on $\V$.
For every $F\in\V'$, the Lax--Milgram theorem therefore gives
a unique $u\in\V$ satisfying
\[
 (B_\nu u)(S_+v)=F(S_+v)
 \qquad\text{for every }v\in\V.
\]
Since $S_+$ maps $\V$ onto itself, this is equivalent to
$B_\nu u=F$. Thus $B_\nu$ is bijective, and the
Lax--Milgram estimate gives boundedness of its inverse.
This proves that $B_\nu\colon\V\to\V'$ is a bounded isomorphism.

Identify $\mathrm L^2(G)$ with its image in $\V'$.
For $u\in\V$, we have
\[
 (B_1u-u)(v)=\mathfrak a_A(u,v),
 \qquad v\in\V.
\]
Since $u\in\mathrm L^2(G)$, the functional
$\mathfrak a_A(u,\cdot)$ is represented by an
$\mathrm L^2(G)$ function if and only if
$B_1u\in\mathrm L^2(G)$. The definition of $H_A$
therefore gives
\[
 \mathrm D(H_A)
 =\{u\in\V:B_1u\in\mathrm L^2(G)\}
 =B_1^{-1}\bigl(\mathrm L^2(G)\bigr).
\]

Since $\mathrm L^2(G)$ is dense in $\V'$ and $B_1$ is an
isomorphism, $\mathrm D(H_A)$ is dense in $\V$, and hence
in $\mathrm L^2(G)$. To prove closedness, suppose that $u_n\in\mathrm D(H_A)$,
$u_n\to u$ and $H_Au_n\to f$ in $\mathrm L^2(G)$.
Then
\[
 B_1u_n=u_n+H_Au_n\longrightarrow u+f
 \quad\text{in }\V'.
\]
Boundedness of $B_1^{-1}$ gives
$u_n\to B_1^{-1}(u+f)$ in $\V$. Comparing the limits in
$\mathrm L^2(G)$ shows that $u=B_1^{-1}(u+f)$.
Thus $u\in\mathrm D(H_A)$ and $H_Au=f$, proving closedness.

For $u\in\mathrm D(H_A)$, ellipticity gives
\[
 \operatorname{Re}\langle H_Au,u\rangle
 =\operatorname{Re}\mathfrak a_A(u,u)
 \ge a_0\|\nabla_Xu\|_2^2\ge0.
\]
Thus $H_A$ is accretive. Moreover, for every $\nu>0$ and
$f\in\mathrm L^2(G)$, the function $u=B_\nu^{-1}f$ belongs
to $\mathrm D(H_A)$ and satisfies $(\nu+H_A)u=f$.
Hence $\operatorname{Ran}(\nu+H_A)=\mathrm L^2(G)$, which
proves maximal accretivity.

We next verify the distributional characterization.
For $u\in\V$ and $\varphi\in\mathrm C_0^\infty(G)$,
the spectral definition of the transport form gives
\[
 \mathfrak t(u,\varphi)=-\langle u,T\varphi\rangle.
\]
It follows that the identity
$\mathfrak a_A(u,\varphi)=\langle f,\varphi\rangle$
for all such $\varphi$ is precisely
\eqref{eq:distribution-realization}. Conversely, if
$u\in\V$ satisfies \eqref{eq:distribution-realization},
then the form identity holds initially for
$\varphi\in\mathrm C_0^\infty(G)$. Both sides are continuous
in the $\V$-norm, so Lemma~\ref{lem:form-density} extends
the identity to every $v\in\V$. This proves the asserted distributional characterization
of the domain and action of $H_A$.

Finally, let $H_A^\dagger$ denote the $\mathrm L^2$ realization
of the adjoint form $\mathfrak a_A^*$. The preceding argument,
using $S_-$ from Lemma~\ref{lem:kinetic-hidden-coercivity},
shows that $H_A^\dagger$ is also densely defined, closed and
maximal accretive. For $f,g\in\mathrm L^2(G)$, set
\[
 u:=(I+H_A)^{-1}f,
 \qquad
 v:=(I+H_A^\dagger)^{-1}g.
\]
The two variational identities yield
\[
 \begin{split}
 \langle f,v\rangle
 &=\langle u,v\rangle+\mathfrak a_A(u,v)=\langle u,v\rangle+\overline{\mathfrak a_A^*(v,u)}
 =\langle u,g\rangle.
 \end{split}
\]
Since $f$ and $g$ are arbitrary,
\[
 \bigl((I+H_A)^{-1}\bigr)^*
 =(I+H_A^\dagger)^{-1}.
\]
This identifies $H_A^\dagger$ with the Hilbert-space adjoint
$H_A^*$, including equality of their operator domains.
Applying the distributional characterization to the adjoint
form gives the stated differential expression. The proof of the proposition is thereby complete.
\end{proof}

For $F\in\mathrm L^2(G;\C^m)$, define divergence as an
element of $\V'$ by
\[
 (\operatorname{div}_X F)(v)
 =-\langle F,\nabla_Xv\rangle.
\]
This agrees with the distributional divergence on smooth test functions.

\begin{lem}\label{lem:scaled-resolvent}
For $\lambda>0$, $f\in\mathrm L^2(G)$ and
$F\in\mathrm L^2(G;\C^m)$, the equation
\[
 (I+\lambda^2H_A)u=f+\lambda\operatorname{div}_X F
\]
has a unique variational solution in $\V$, and
\begin{equation}\label{eq:scaled-resolvent-energy}
 \|u\|_2+\lambda E(u)
 \lesssim\|f\|_2+\|F\|_2.
\end{equation}
For $F=0$, the solution is $u=R_\lambda f$, where
$R_\lambda:=(I+\lambda^2H_A)^{-1}$ is a contraction on
$\mathrm L^2(G)$. Moreover, the variational extension of
$R_\lambda$ to $\V'$ satisfies
\begin{equation}\label{eq:divergence-resolvent}
 \lambda\|R_\lambda\operatorname{div}_X F\|_2
 +\lambda^2E(R_\lambda\operatorname{div}_X F)
 \lesssim\|F\|_2.
\end{equation}
The same statements hold for $H_A^*$.
\end{lem}

\begin{proof}
Fix $\lambda>0$. By Lemma~\ref{lem:kinetic-hidden-coercivity}
and the Lax--Milgram argument in the proof of
Proposition~\ref{prop:kinetic-realization}, there is a unique
$u\in\V$ satisfying
\[
 \langle u,v\rangle+\lambda^2\mathfrak a_A(u,v)
 =\langle f,v\rangle-\lambda\langle F,\nabla_Xv\rangle
 \qquad\text{for every }v\in\V.
\]
Indeed, the right-hand side is a bounded antilinear functional
on $\V$, and the form on the left becomes coercive after
replacing the test function $v$ by $S_+v$. To obtain an estimate uniform in $\lambda$, take $v=S_+u$
and real parts. The bounds
\eqref{eq:symmetrizer-energy}--\eqref{eq:symmetrizer-mass},
together with boundedness of $S_+$ on $\mathrm L^2(G)$, give
\[
 \begin{split}
 \|u\|_2^2+\lambda^2E(u)^2
 &\lesssim
 \operatorname{Re}\bigl(
   \langle f,S_+u\rangle
   -\lambda\langle F,\nabla_XS_+u\rangle
 \bigr)\\
 &\lesssim\|f\|_2\|u\|_2+\lambda\|F\|_2E(u)\\
 &\le(\|f\|_2+\|F\|_2)(\|u\|_2+\lambda E(u)).
 \end{split}
\]
This proves \eqref{eq:scaled-resolvent-energy}, with implicit
constants independent of $\lambda$.

If $F=0$, Proposition~\ref{prop:kinetic-realization} identifies
the variational solution with $u=R_\lambda f$.
Taking $v=u$ in the variational identity and using ellipticity,
we obtain
\[
 \|u\|_2^2+a_0\lambda^2\|\nabla_Xu\|_2^2
 \le\operatorname{Re}\langle f,u\rangle
 \le\|f\|_2\|u\|_2.
\]
Hence $\|R_\lambda f\|_2\le\|f\|_2$, proving the contraction
bound. If $f=0$, the solution is
$u=\lambda R_\lambda\operatorname{div}_X F$, where
$R_\lambda$ is understood through its variational extension
to $\V'$. Substitution into
\eqref{eq:scaled-resolvent-energy} gives
\[
 \lambda\|R_\lambda\operatorname{div}_X F\|_2
 +\lambda^2E(R_\lambda\operatorname{div}_X F)
 \lesssim\|F\|_2,
\]
which is \eqref{eq:divergence-resolvent}.

Finally, Proposition~\ref{prop:kinetic-realization} identifies
$H_A^*$ with the realization of the adjoint form.
Repeating the argument with $S_-$ from
Lemma~\ref{lem:kinetic-hidden-coercivity} proves all the
corresponding assertions for $H_A^*$.
\end{proof}

\begin{rem}\label{rem:injective}
Both $H_A$ and $H_A^*$ are injective and have dense
range. Indeed, a vector in either kernel has $\nabla_Xu=0$ by
accretivity, so it is independent of $X$. Since Lebesgue measure
of $\R^m_X$ is infinite, such a function in $\mathrm L^2(G)$
must vanish. The adjoint identity and
$\overline{\Ran(H_A)}=(\ker H_A^*)^\perp$, together with the
corresponding formula for $H_A^*$, give the dense range assertions.
No bounded inverse is asserted.
\end{rem}

\section{Reduction to two quadratic estimates}
\label{sec:quadratic-reduction}

Theorem~\ref{thm:main} reduces to two quadratic estimates,
one for $H_A$ and one for $H_A^*$. The reduction follows
\cite[Section~6]{AtaeiEgertNystrom}, with the regularized
transport symmetrizers $S_\pm$ replacing the time
Hilbert-transform modification. We retain the details needed
to obtain the lower bounds and identify the full square-root
domains.

By maximal accretivity, $H_A$ and $H_A^*$ are sectorial of
angle $\pi/2$. We use their holomorphic functional calculi
on larger sectors
\[
 \Sigma_\theta
 :=\{z\in\C\setminus\{0\}:|\arg z|<\theta\},
 \qquad \pi/2<\theta<\pi.
\]
No sectoriality with angle strictly smaller than $\pi/2$ is
required. For background on the functional calculus and its
relation to quadratic estimates, see
\cite[Sections~D--G]{AlbrechtDuongMcIntosh} and
\cite[Chapters~3 and~7]{HaaseFunctionalCalculus}.

\begin{lem}\label{lem:universal-accretive-quadratic}
Let $B$ be maximal accretive on a complex Hilbert space
$\mathcal H$. Then
\begin{equation}\label{eq:abstract-square-root-quadratic}
 \|B^{1/2}u\|^2
 \sim
 \int_0^\infty
 \|\lambda B(I+\lambda^2B)^{-1}u\|^2
 \frac{\d\lambda}{\lambda},
 \qquad u\in\mathrm D(B^{1/2}),
\end{equation}
with universal comparison constants.
\end{lem}

\begin{proof}
Set $\psi(z):=z^{1/2}(1+z)^{-1}$, using the principal
branch of the square root. The bounded holomorphic functional
calculus for maximal accretive operators and its Hilbert-space
quadratic-estimate characterization give
\[
 \int_0^\infty
 \|\psi(\lambda^2B)f\|^2\frac{\d\lambda}{\lambda}
 \sim\|f\|^2,
 \qquad f\in\overline{\Ran(B)};
\]
see \cite[Theorems~7.1.7 and~7.3.1]{HaaseFunctionalCalculus}.
Here a sector angle $\theta\in(\pi/2,\pi)$ is fixed once
and for all. The Cayley transform and von Neumann's inequality
give a contractive half-plane calculus, so the comparison
constants are universal once $\theta$ is fixed.
For a maximal accretive operator, the decomposition
$\mathcal H=\ker B\oplus\overline{\Ran(B)}$ is orthogonal,
since $\ker B=\ker B^*$, and reduces $B$. Its restriction to $\overline{\Ran(B)}$
is injective and maximal accretive with dense range, so the
quadratic-estimate characterization applies on this subspace. For $u\in\mathrm D(B^{1/2})$, one has
$B^{1/2}u\in\overline{\Ran(B)}$, and the functional calculus
gives
\[
 \psi(\lambda^2B)B^{1/2}u
 =\lambda B(I+\lambda^2B)^{-1}u.
\]
Applying the preceding estimate to $f=B^{1/2}u$ proves
\eqref{eq:abstract-square-root-quadratic}.
\end{proof}

The equivalence \eqref{eq:abstract-square-root-quadratic}
compares the square function with $\|B^{1/2}u\|$, but does
not identify this norm with the differential energy $E(u)$.
Establishing that identification for $B=H_A$ and $B=H_A^*$
is the remaining task.

\begin{prop}\label{prop:quadratic-reduction}
The following assertions are equivalent.
\begin{enumerate}[label=(\roman*)]
 \item There is a constant $C_0>0$ such that, for $B=H_A$
 and $B=H_A^*$,
 \begin{equation}\label{eq:quadratic-criterion}
  \int_0^\infty
  \|\lambda B(I+\lambda^2B)^{-1}u\|_2^2
  \frac{\d\lambda}{\lambda}
  \le C_0E(u)^2,
  \qquad u\in\V.
 \end{equation}
 \item One has
 \[
  \mathrm D(H_A^{1/2})
  =\mathrm D((H_A^*)^{1/2})=\V
 \]
 and
 \[
  \|H_A^{1/2}u\|_2\sim E(u)
  \sim\|(H_A^*)^{1/2}u\|_2,
  \qquad u\in\V.
 \]
\end{enumerate}
The comparison constants in (ii) depend only on
$m,a_0,a_1,C_0$. Conversely, the upper constants in (ii)
determine an admissible constant $C_0$ in (i).
\end{prop}

\begin{proof}
Assume (i) and write $H:=H_A$. Throughout this implication,
implicit constants may additionally depend on $C_0$.
For $u\in\mathrm D(H)$,
Lemma~\ref{lem:universal-accretive-quadratic} and
\eqref{eq:quadratic-criterion} give
\begin{equation}\label{Estsquare}
 \|H^{1/2}u\|_2^2
 \lesssim
 \int_0^\infty
 \|\lambda H(I+\lambda^2H)^{-1}u\|_2^2
 \frac{\d\lambda}{\lambda}
 \lesssim E(u)^2.
\end{equation}
By Proposition~\ref{prop:kinetic-realization},
$\mathrm D(H)$ is dense in $\V$. Closedness of $H^{1/2}$
therefore extends the estimate in \eqref{Estsquare} to every $u\in\V$ and gives
$\V\subset\mathrm D(H^{1/2})$. Applying the same argument
to $H^*$ yields
$\V\subset\mathrm D((H^*)^{1/2})$ and
\begin{equation}\label{eq:two-root-upper-bounds}
 \|H^{1/2}u\|_2+\|(H^*)^{1/2}u\|_2
 \lesssim E(u),
 \qquad u\in\V.
\end{equation}
We obtain the lower estimates from hidden coercivity.
Let $u\in\mathrm D(H)$. Since
$S_+u\in\V\subset\mathrm D((H^*)^{1/2})$, the
fractional-power identities
$H=(H^{1/2})^2$ and
$(H^{1/2})^*=(H^*)^{1/2}$ give
\[
 \mathfrak a_A(u,S_+u)
 =\langle Hu,S_+u\rangle
 =\langle H^{1/2}u,(H^*)^{1/2}S_+u\rangle.
\]
Using Lemma~\ref{lem:kinetic-hidden-coercivity},
\eqref{eq:two-root-upper-bounds} and
\eqref{eq:symmetrizer-energy}, we conclude that
\[
 \begin{split}
 E(u)^2
 \lesssim\operatorname{Re}\mathfrak a_A(u,S_+u)
 &\le\|H^{1/2}u\|_2\|(H^*)^{1/2}S_+u\|_2\\
 &\lesssim\|H^{1/2}u\|_2E(S_+u)
 \lesssim\|H^{1/2}u\|_2E(u).
 \end{split}
\]
Hence
\[
 E(u)\lesssim\|H^{1/2}u\|_2,
 \qquad u\in\mathrm D(H).
\]
The adjoint form and $S_-$ similarly give
\[
 E(u)\lesssim\|(H^*)^{1/2}u\|_2,
 \qquad u\in\mathrm D(H^*).
\]
It remains to extend these lower bounds and prove the reverse
domain inclusions. Let $u\in\mathrm D(H^{1/2})$ and set
\[
 u_n:=n(n+H)^{-1}u\in\mathrm D(H)\subset\V.
\]
Strong convergence of the resolvents and their commutation
with fractional powers give
\[
 u_n\longrightarrow u,
 \qquad
 H^{1/2}u_n
 =n(n+H)^{-1}H^{1/2}u
 \longrightarrow H^{1/2}u
 \quad\text{in }\mathrm L^2(G).
\]
The lower bound applied to $u_n-u_k$ shows that $(u_n)$ is
Cauchy in the full norm of $\V$. Completeness and the
continuous embedding $\V\hookrightarrow\mathrm L^2(G)$
therefore imply that $u\in\V$ and $u_n\to u$ in $\V$.
Passing to the limit also gives
$E(u)\lesssim\|H^{1/2}u\|_2$.
Thus $\mathrm D(H^{1/2})=\V$, with the required norm
equivalence. The same argument for $H^*$ proves the
corresponding assertion for the adjoint, establishing (ii).

Conversely, assume (ii). For $B=H_A$ and $B=H_A^*$,
Lemma~\ref{lem:universal-accretive-quadratic} gives
\[
 \int_0^\infty
 \|\lambda B(I+\lambda^2B)^{-1}u\|_2^2
 \frac{\d\lambda}{\lambda}
 \lesssim\|B^{1/2}u\|_2^2
 \lesssim E(u)^2,
 \qquad u\in\V.
\]
This proves (i), with $C_0$ controlled by the upper constants
in (ii).
\end{proof}

\section{The Rockland operator and intrinsic smoothing}
\label{sec:flat-rockland}

In this section we construct an intrinsic approximation to the identity whose square
functions are controlled by the energy $E$. The comparison operator is
independent of $A$, and its construction does not use a square-root
estimate for $H_A$.

\subsection{The positive Rockland operator}

Put $P=-\Delta_X$ and consider the left-invariant differential
operator
\[
 \mathcal R=P^2+T^\ast T=P^2-T^2.
\]
Both $P$ and $T$ have homogeneous degree two under the kinetic
dilations, so $\mathcal R$ has homogeneous degree four.
Moreover, integration by parts gives
\[
 \langle\mathcal Ru,u\rangle
 =\|Pu\|_2^2+\|Tu\|_2^2\ge0,
 \qquad u\in\mathrm C_0^\infty(G).
\]
Thus $\mathcal R$ combines horizontal differentiation and
transport in a nonnegative homogeneous operator. This choice is adapted to the energy space as its fourth root,
defined through the nonnegative self-adjoint realization
established below, has the same homogeneous degree as
$\nabla_X$ and $|T|^{1/2}$. The energy comparison below
makes this correspondence precise.
Furthermore, $\mathcal R$ is unchanged when $T$ is replaced
by $-T$, allowing the same smoothing operator to be used
for both $H_A$ and its adjoint.

In the following we first verify directly that $\mathcal R$ satisfies the
Rockland condition. The argument uses positivity and the
fact that $\partial_{X_1},\ldots,\partial_{X_m},T$
generate the Lie algebra of $G$. This will allow us to apply the general theory of positive
Rockland operators presented in \cite{FischerRuzhansky}
and \cite[Chapter~4]{FischerRuzhanskyBook}.

The Rockland condition requires $\mathrm d\pi(\mathcal R)$
to be injective on the space of smooth vectors of every
nontrivial irreducible unitary representation $\pi$ of
the group $(G,\circ)$,  see
\cite[Definition~2.7]{FischerRuzhansky}. Here a unitary representation assigns to each $z\in G$
a unitary operator $\pi(z)$ on a complex Hilbert space
$\mathcal H_\pi$, with
\[
 \pi(z\circ w)=\pi(z)\pi(w),
 \qquad \pi(e)=I.
\]
We also require strong continuity, meaning that
$z\mapsto\pi(z)v$ is continuous for every $v\in\mathcal H_\pi$.
The representation is irreducible if the only closed subspaces
invariant under all the operators $\pi(z)$ are $\{0\}$ and
$\mathcal H_\pi$. Nontriviality excludes the one-dimensional
representation $\pi(z)=1$ for every $z\in G$. A vector $v\in\mathcal H_\pi$ is smooth if its orbit map
$z\mapsto\pi(z)v$ is smooth as an $\mathcal H_\pi$-valued map.
On the space $\mathcal H_\pi^\infty$ of smooth vectors,
the infinitesimal representation is defined by
\[
 \mathrm d\pi(Z)v
 :=\left.\frac{\d}{\d s}\pi(\exp(sZ))v\right|_{s=0},
 \qquad v\in\mathcal H_\pi^\infty,
\]
where $Z$ is a Lie algebra element, identified with a
left-invariant vector field, and $\exp(sZ)$ is its
one-parameter subgroup. These infinitesimal operators preserve $\mathcal H_\pi^\infty$
and satisfy
\[
 [\mathrm d\pi(Z),\mathrm d\pi(W)]
 =\mathrm d\pi([Z,W]).
\]
They therefore extend by linearity and composition to
left-invariant differential operators, defining in particular
$\mathrm d\pi(\mathcal R)$ on $\mathcal H_\pi^\infty$, see \cite[Proposition~1.7.3]{FischerRuzhanskyBook}.

Fix a nontrivial irreducible unitary representation $\pi$
on $\mathcal H_\pi$, and suppose that
$v\in\mathcal H_\pi^\infty$ satisfies
$\mathrm d\pi(\mathcal R)v=0$.
The operator $\mathrm d\pi(\mathcal R)$ is obtained by replacing
each vector field in $\mathcal R$ by its infinitesimal generator
in the representation $\pi$, preserving the order of composition.
Recalling that $P=-\Delta_X$, we have
\[
 \mathrm d\pi(\mathcal R)
 =\mathrm d\pi(P)^2-\mathrm d\pi(T)^2
 =\biggl(-\sum_{j=1}^m
       \mathrm d\pi(\partial_{X_j})^2\biggr)^2
       -\mathrm d\pi(T)^2
 \quad\text{on }\mathcal H_\pi^\infty.
\]
All compositions are well defined on this space, since
infinitesimal generators preserve smooth vectors. Unitarity implies that
 the infinitesimal  generators
$\mathrm d\pi(\partial_{X_j})$ and $\mathrm d\pi(T)$
are skew-symmetric  on $\mathcal H_\pi^\infty$. Indeed,
differentiating
\[
 \langle\pi(\exp(sZ))v,\pi(\exp(sZ))w\rangle_{\mathcal H_\pi}
 =\langle v,w\rangle_{\mathcal H_\pi}
\]
at $s=0$ gives
\[
 \langle\mathrm d\pi(Z)v,w\rangle_{\mathcal H_\pi}
 =-\langle v,\mathrm d\pi(Z)w\rangle_{\mathcal H_\pi},
 \qquad v,w\in\mathcal H_\pi^\infty.
\]
Consequently, $\mathrm d\pi(P)$
is symmetric, and these symmetry identities give
\[
 0=\langle\mathrm d\pi(\mathcal R)v,v\rangle_{\mathcal H_\pi}
  =\|\mathrm d\pi(P)v\|_{\mathcal H_\pi}^2
   +\|\mathrm d\pi(T)v\|_{\mathcal H_\pi}^2.
\]
Hence both terms vanish, so
$\mathrm d\pi(P)v=\mathrm d\pi(T)v=0$. This gives
\[
 0=\langle\mathrm d\pi(P)v,v\rangle_{\mathcal H_\pi}
  =\sum_{j=1}^m
    \|\mathrm d\pi(\partial_{X_j})v\|_{\mathcal H_\pi}^2.
\]
Thus $\mathrm d\pi(\partial_{X_j})v=0$ for every $j$.
Since the infinitesimal representation preserves Lie brackets,
the relations $[\partial_{X_j},T]=-\partial_{Y_j}$ imply
\[
 \mathrm d\pi(\partial_{Y_j})v
 =-[\mathrm d\pi(\partial_{X_j}),\mathrm d\pi(T)]v=0.
\]
The fields $\partial_{X_j}$, $\partial_{Y_j}$ and $T$
span the Lie algebra of $G$. Consequently,
$\mathrm d\pi(Z)v=0$ for every Lie algebra element $Z$,
so every one-parameter subgroup fixes $v$.
Since $G$ is connected, these subgroups generate $G$,
and therefore $\pi(z)v=v$ for every $z\in G$.
If $v\ne0$, its span would be an invariant line on which
$\pi$ acts trivially. Irreducibility would then force
$\mathcal H_\pi$ to equal this line, making $\pi$ the
trivial representation. This contradicts the choice of $\pi$.
Hence $v=0$, proving the Rockland condition for $\mathcal R$.

By the realization theorem for positive Rockland operators
\cite[Theorem~3.3(4)--(5)]{FischerRuzhansky}, the operator
$\mathcal R$, initially defined on $\mathrm C_0^\infty(G)$,
is essentially self-adjoint on $\mathrm L^2(G)$.
We continue to denote its nonnegative self-adjoint closure
by $\mathcal R$. The spectral theorem defines its heat semigroup
$e^{-t\mathcal R}$, $t>0$.

The Rockland heat-kernel theorem
\cite[Theorem~2.8]{FischerRuzhansky} gives Schwartz kernels $h_s$ for
$e^{-s\mathcal R}$, with
\[
 h_s(z)=s^{-\mathfrak q/4}h_1(\delta_{s^{-1/4}}z),
 \qquad \int_Gh_s\d z=1.
\]
Unlike the classical Gaussian heat kernel, $h_s$ changes sign. Indeed,
integrating $h_s$ in $(Y,t)$  gives the sign-changing biharmonic
heat kernel. The arguments below use its normalization and Schwartz bounds,
together with the spectral calculus of $\mathcal R$, and do not
require pointwise positivity. We parametrize this semigroup by the intrinsic scale
$\lambda>0$ and set 
\begin{equation}\label{eq:rockland-approximation}
 P_\lambda:=e^{-\lambda^4\mathcal R},
 \qquad p_\lambda:=h_{\lambda^4}.
\end{equation}
We use the convolution convention
\[
 (f*g)(z)
 =\int_G f(\zeta)g(\zeta^{-1}\circ z)\d\zeta.
\]
Since $\mathcal R$ is left-invariant, its heat semigroup
commutes with left translations and is represented by right
convolution, $P_\lambda f=f*p_\lambda$. The heat-kernel scaling gives
\[
 p_\lambda(z)
 =\lambda^{-\mathfrak q}p_1(\delta_{1/\lambda}z).
\]
More generally, Hulanicki's spectral-multiplier theorem
\cite[Theorem~4.5.1 and Corollary~4.5.2]{FischerRuzhanskyBook}
shows that, for $\phi\in\mathscr S(\R)$, the operator
$\phi(\mathcal R)$ has a Schwartz convolution kernel.
The map from $\phi$ to this kernel is continuous in the
Schwartz topologies. The heat kernel is included by taking
a Schwartz function agreeing with $s\mapsto e^{-s}$ on
$[0,\infty)$.  Furthermore, $\mathcal R$ is injective. Indeed, if
$u\in\mathrm D(\mathcal R)$ and $\mathcal Ru=0$, the spectral
theorem gives $u=e^{-s\mathcal R}u=u*h_s$ for every $s>0$. As $\|h_s\|_2=s^{-\mathfrak q/8}\|h_1\|_2$, it follows that
$\|u\|_\infty\lesssim s^{-\mathfrak q/8}\|u\|_2$ for every $s>0$. Letting $s\to\infty$ yields $u=0$, proving injectivity. Since $\mathcal R$ is injective, the spectral theorem defines
its complex powers $\mathcal R^\alpha$, $\alpha\in\C$,
on their natural spectral domains. In particular,
$\mathcal R^{i\tau}$ is unitary on $\mathrm L^2(G)$
for every $\tau\in\R$.

\begin{rem}
The Rockland condition is an analogue of ellipticity adapted
to the group structure. On Euclidean space, a homogeneous
constant-coefficient operator is elliptic when its Fourier
symbol does not vanish at any nonzero frequency.
On $G$, irreducible unitary representations play the role
of Fourier frequencies, and the corresponding symbols
$\mathrm d\pi(\mathcal R)$ are operators. The Rockland
condition requires these operators to be injective on
smooth vectors for every nontrivial irreducible unitary representation. The analytic significance of this condition is given by
the Helffer--Nourrigat theorem
\cite[Theorem~4.1.12]{FischerRuzhanskyBook} which states, for homogeneous
left-invariant differential operators, that the Rockland condition
is equivalent to hypoellipticity. Thus, for a distribution
$u$ on an open set $O\subset G$,
\[
 \mathcal Ru\in\mathrm C^\infty(O)
 \quad\Longrightarrow\quad
 u\in\mathrm C^\infty(O).
\]
In our case, the bracket relations
$[\partial_{X_j},T]=-\partial_{Y_j}$ explain how horizontal
differentiation and transport together yield regularity
in all variables, despite the absence of an explicit
$Y$-diffusion term in $\mathcal R$. Although the fields $\partial_{X_1},\ldots,\partial_{X_m},T$
generate the Lie algebra, H\"ormander's classical
sum-of-squares theorem does not apply directly to
$\mathcal R$, since $P=-\Delta_X$ is a second-order
operator and $P^2$ is fourth order.
In the present argument, hypoellipticity follows instead
from the direct verification of the Rockland condition
and the Helffer--Nourrigat theorem.
\end{rem}

\subsection{Fractional powers and a common core}

Let $\mathscr S(G)$ be the Schwartz class and set
\[
 \mathscr S_0(G):=\left\{f\in\mathscr S(G):
       \int_G z^\beta f(z)\d z=0
       \text{ for every multi-index }\beta\right\}.
\]
Here $z^\beta$ denotes a coordinate monomial in
$z=(X,Y,t)$. Since $\mathcal R$ is a positive Rockland operator of
homogeneous degree four, the fractional-power theory in
\cite[Section~4.3]{FischerRuzhanskyBook} applies, see also
\cite[Section~3.2]{FischerRuzhansky}.
On $\mathrm L^2(G)$, these powers agree with those defined
by the spectral theorem. We use their general properties
and homogeneous kernel representations below. The role of $\mathscr S_0(G)$ is to provide a common smooth
class on which the fractional-power and Riesz-transform
identities can be justified. Its vanishing moments ensure
that applying negative powers of $\mathcal R$ preserves
Schwartz regularity. The following lemma makes this precise
by establishing invariance under all complex powers, together
with the density and core properties needed to extend the
resulting estimates.

\begin{lem}\label{wr:common-core}
The space $\mathscr S_0(G)$ is dense in $\V$ in its
inhomogeneous energy norm and is a core for
$\mathcal R^{1/4}$ on $\mathrm L^2(G)$, that is, it is
dense in $\mathrm D(\mathcal R^{1/4})$ with respect to
the graph norm
\[
 u\mapsto
 \bigl(\|u\|_2^2+\|\mathcal R^{1/4}u\|_2^2\bigr)^{1/2}.
\]
For every $\alpha\in\C$, one has
$\mathscr S_0(G)\subset\mathrm D(\mathcal R^\alpha)$,
and $\mathcal R^\alpha$ maps $\mathscr S_0(G)$ continuously
into itself for the topology inherited from $\mathscr S(G)$.
For every $f\in\mathscr S_0(G)$, the map
$\alpha\mapsto\mathcal R^\alpha f$ is holomorphic on
the whole complex plane, with values in $\mathscr S_0(G)$
endowed with the topology inherited from $\mathscr S(G)$.
On this space, the usual power and adjoint pairing
identities hold.
\end{lem}

The proof of Lemma~\ref{wr:common-core} is given in Appendix~\ref{app:common-core}.
It establishes density in $\V$ and the core property for
$\mathcal R^{1/4}$ independently of the energy comparison
proved below. Invariance under complex powers follows
from moment cancellation and a Mellin representation
of negative powers.

\subsection{Riesz transforms and the energy comparison}

Following \cite[Definition~2.3]{FischerRuzhansky}, we call
a distribution $K$ on $G$ a \emph{kernel of type $\nu$}
if it is homogeneous of degree $\nu-\mathfrak q$ under
the kinetic dilations and is smooth away from the identity.
In particular, a \emph{type-zero kernel} is homogeneous
of degree $-\mathfrak q$. The associated right convolution
operator acts by $f\mapsto f*K$.

\begin{lem}\label{wr:zero-order-and-bip}
If $D$ is a homogeneous left-invariant differential operator of degree
$d\in\{1,2,3,4\}$, then
\begin{equation}\label{wr:differential-riesz}
 \|D\mathcal R^{-d/4}f\|_2\lesssim_D\|f\|_2,
 \qquad f\in\mathscr S_0(G).
\end{equation}
These transforms have right convolution kernels of type zero.
\end{lem}

\begin{proof}
The homogeneous Sobolev estimate in
\cite[Theorem~4.14]{FischerRuzhansky} gives
\[
 \|Dg\|_2\lesssim_D\|\mathcal R^{d/4}g\|_2,
 \qquad g\in\mathscr S_0(G).
\]
For $f\in\mathscr S_0(G)$,
Lemma~\ref{wr:common-core} allows us to take
$g=\mathcal R^{-d/4}f\in\mathscr S_0(G)$.
The power identities then give
\[
 \|D\mathcal R^{-d/4}f\|_2
 \lesssim_D\|\mathcal R^{d/4}\mathcal R^{-d/4}f\|_2
 =\|f\|_2,
\]
which proves \eqref{wr:differential-riesz}. For the kernel assertion, note that
$0<d<\mathfrak q=4m+2$.
By \cite[Corollary~3.11]{FischerRuzhansky},
$\mathcal R^{-d/4}$ has a right convolution kernel $k_d$
of type $d$. Thus $k_d$ is smooth away from the identity
and homogeneous of degree $d-\mathfrak q$.
Left invariance of $D$ gives
\[
 D\mathcal R^{-d/4}f
 =D(f*k_d)=f*(Dk_d),
\]
where $Dk_d$ is understood distributionally.
By \cite[Proposition~2.4(3)]{FischerRuzhansky}, this
distribution is a kernel of type zero. In particular,
it is homogeneous of degree $-\mathfrak q$ and smooth
away from the identity.
\end{proof}

\begin{rem}\label{rem:rockland-riesz-transforms}
For a homogeneous left-invariant differential operator $D$
of degree $d$, we refer to  $D\mathcal R^{-d/4}$ as a \emph{Rockland Riesz transform}.
This construction generalizes the classical Riesz transform
$\partial_j(-\Delta)^{-1/2}$.
Since $\mathcal R$ has homogeneous degree four,
$\mathcal R^{-d/4}$ has degree $-d$, which compensates
for the degree of $D$. The composition therefore has
homogeneous degree zero. Examples include
\[
 \partial_{X_j}\mathcal R^{-1/4},
 \qquad
 T\mathcal R^{-1/2},
 \qquad
 \partial_{Y_j}\mathcal R^{-3/4},
\]
corresponding to differential operators of degrees one,
two and three, respectively.
These compositions are initially defined on $\mathscr S_0(G)$.
Lemma~\ref{wr:zero-order-and-bip} shows that they extend
to bounded operators on $\mathrm L^2(G)$.
\end{rem}

\begin{lem}\label{wr:fractional-transport-riesz}
For every $0\le\alpha\le2$,
\begin{equation}\label{wr:fractional-transport}
 \bigl\||T|^\alpha\mathcal R^{-\alpha/2}f\bigr\|_2
 \le C_\alpha\|f\|_2,
 \qquad f\in\mathscr S_0(G),
\end{equation}
where $C_\alpha$ depends only on $m$ and $\alpha$.
\end{lem}

\begin{proof}
By Lemma~\ref{wr:zero-order-and-bip}, the operator
$K_T:=|T|^2\mathcal R^{-1}=-T^2\mathcal R^{-1}$
is bounded on $\mathrm L^2(G)$. Fix $M\ge1$ with
$\|K_T\|_{\mathrm L^2\to\mathrm L^2}\le M$. For $f\in\mathscr S_0(G)$, consider
\[
 F(z)f:=|T|^{2z}\mathcal R^{-z}f,
 \qquad 0\le\operatorname{Re}z\le1.
\]
We first justify the analytic properties of this family.
Under the partial Fourier transform in $(Y,t)$, the operator
$|T|$ is multiplication by $|q|$, where $q=\tau-X\cdot\eta$.
Its zero set has Lebesgue measure zero, so $|T|$ is injective.
Thus the imaginary powers of both $|T|$ and $\mathcal R$ are unitary.
By Lemma~\ref{wr:common-core}, the map
$z\mapsto\mathcal R^{-z}f$ is holomorphic in the Schwartz topology,
and hence in the graph norm of $|T|^2=-T^2$.
On each compact subset of the open strip, the multipliers
$|q|^{2z}$ and their first $z$ derivatives are bounded by
$C(1+|q|^2)$, with $C$ depending on that subset.
This proves that $F(z)f$ is $\mathrm L^2$-valued and holomorphic
there. The bound $|q|^{2\operatorname{Re}z}\le1+|q|^2$
and dominated convergence give continuity on the closed strip.
For $z=x+ib$, the inequality $|q|^{2x}\le1+|q|^2$ and
the factorization
$|T|^2\mathcal R^{-z}=K_T\mathcal R^{1-z}$ give
\[
 \begin{split}
 \|F(z)f\|_2
 &\le\|\mathcal R^{-x}f\|_2
       +M\|\mathcal R^{1-x}f\|_2\le\|\mathcal R^{-1}f\|_2
       +(1+M)\|f\|_2+M\|\mathcal Rf\|_2.
 \end{split}
\]
Thus $F(z)f$ is bounded throughout the strip.
On its boundary lines,
\[
 F(ib)=|T|^{2ib}\mathcal R^{-ib},
 \qquad
 F(1+ib)=|T|^{2ib}K_T\mathcal R^{-ib},
\]
so the corresponding operator norms are at most $1$ and $M$. Fix $f\in\mathscr S_0(G)$ and $g\in\mathrm L^2(G)$, and set
\[
 h_{f,g}(z):=\langle F(z)f,g\rangle.
\]
The preceding argument shows that $h_{f,g}$ is holomorphic
in the open strip $0<\operatorname{Re}z<1$ and continuous
and bounded on its closure. Moreover, the boundary
estimates give
\[
 |h_{f,g}(ib)|\le\|f\|_2\|g\|_2,
 \qquad
 |h_{f,g}(1+ib)|\le M\|f\|_2\|g\|_2,
 \qquad b\in\R.
\]
Hadamard's three-lines theorem
\cite[Lemma~1.3.5]{GrafakosClassical} therefore yields
\[
 |h_{f,g}(\theta)|
 \le
 \bigl(\|f\|_2\|g\|_2\bigr)^{1-\theta}
 \bigl(M\|f\|_2\|g\|_2\bigr)^\theta
 =M^\theta\|f\|_2\|g\|_2,
 \qquad 0<\theta<1.
\]
Taking the supremum over $\|g\|_2=1$, we obtain
\[
 \|F(\theta)f\|_2\le M^\theta\|f\|_2.
\]
Setting $\theta=\alpha/2$ proves
\eqref{wr:fractional-transport} for $0<\alpha<2$,
with $C_\alpha=M^{\alpha/2}$.
The endpoint cases $\alpha=0$ and $\alpha=2$ follow
directly from the identity operator and the bound for
$K_T$, respectively.
\end{proof}

\begin{prop}\label{wr:energy-comparison}
One has
\begin{equation}\label{wr:energy-equivalence}
 \mathrm D(\mathcal R^{1/4})=\V,
 \qquad
 \|\mathcal R^{1/4}u\|_2\sim
       \|\nabla_Xu\|_2+\||T|^{1/2}u\|_2.
\end{equation}
Here the comparison constants depend only on $m$.
\end{prop}

\begin{proof}
For $u\in\mathscr S_0(G)$, put $f=\mathcal R^{1/4}u$.
Equations \eqref{wr:differential-riesz} and
\eqref{wr:fractional-transport}, with $\alpha=1/2$, give
\[
 \|\nabla_Xu\|_2+\||T|^{1/2}u\|_2
       \lesssim\|\mathcal R^{1/4}u\|_2.
\]
For the converse, take $g\in\mathscr S_0(G)$. Self-adjointness and
integration by parts give
\begin{align*}
 \langle\mathcal R^{1/4}u,g\rangle
 &=\langle\mathcal Ru,\mathcal R^{-3/4}g\rangle=\langle\nabla_Xu,\nabla_XP\mathcal R^{-3/4}g\rangle+\langle|T|^{1/2}u,|T|^{3/2}\mathcal R^{-3/4}g\rangle.
\end{align*}
The transport identity follows from the spectral calculus of $T$.
Apply \eqref{wr:differential-riesz} with $D=\nabla_XP$ and
\eqref{wr:fractional-transport} with $\alpha=3/2$. It follows that
\[
 |\langle\mathcal R^{1/4}u,g\rangle|
       \lesssim E(u)\|g\|_2.
\]
Density of $\mathscr S_0(G)$ in $\mathrm L^2(G)$ gives the reverse estimate.
The two estimates show that the inhomogeneous norms
\[
 \|u\|_2+\|\mathcal R^{1/4}u\|_2
 \quad\text{and}\quad
 \|u\|_2+E(u)
\]
are equivalent on $\mathscr S_0(G)$.
By Lemma~\ref{wr:common-core}, this space is dense in
$\mathrm D(\mathcal R^{1/4})$ and $\V$ for the respective norms.
Their completions agree as subspaces of $\mathrm L^2(G)$, proving
both the domain identity and the norm equivalence.
\end{proof}

\subsection{Square-function estimates for intrinsic smoothing}

For $s>0$, set
\[
 \psi_0(s)=s^{-1/4}(1-e^{-s}),\qquad
 \psi_1(s)=s^{1/4}e^{-s},\qquad
 \psi_2(s)=s^{1/2}e^{-s}.
\]
These functions belong to $\mathrm L^2((0,\infty),\d s/s)$.
The spectral theorem and injectivity of $\mathcal R$ give,
for every $f\in\mathrm L^2(G)$,
\begin{equation}\label{w:vertical-SF}
 \int_0^\infty\|\psi_i(\lambda^4\mathcal R)f\|_2^2
                         \frac{\d\lambda}{\lambda}
 =\frac14\left(\int_0^\infty|\psi_i(s)|^2\frac{\d s}{s}\right)
       \|f\|_2^2,
 \qquad i=0,1,2.
\end{equation}
Indeed, for each positive spectral value $r$, the substitution
$s=\lambda^4r$ gives
\[
 \int_0^\infty|\psi_i(\lambda^4r)|^2
                   \frac{\d\lambda}{\lambda}
 =\frac14\int_0^\infty|\psi_i(s)|^2\frac{\d s}{s}.
\]
Integrating against the spectral measure associated with $f$
and applying Tonelli's theorem proves the identity.
Injectivity of $\mathcal R$ ensures that its spectral
projection at zero vanishes, so this measure has total
mass $\|f\|_2^2$ on $(0,\infty)$.

\begin{lem}\label{aux:rockland-smoothing}
For every $u\in\V$,
\begin{equation}\label{w:smoothing-input}
 \begin{split}
 \int_0^\infty\bigl[
 &\lambda^{-2}\|(I-P_\lambda)u\|_2^2
  +\lambda^2\|TP_\lambda u\|_2^2+\lambda^2\|\nabla_X^2P_\lambda u\|_2^2
  +\lambda^4\|T\nabla_XP_\lambda u\|_2^2
 \bigr]\frac{\d\lambda}{\lambda}
       \lesssim E(u)^2.
 \end{split}
\end{equation}
The implicit constant depends only on $m$.
The kernels of $P_\lambda$ are scaled Schwartz functions of integral one.
\end{lem}

\begin{proof}
We first take $u\in\mathscr S_0(G)$ and put $v=\mathcal R^{1/4}u$.
For a homogeneous left-invariant differential operator $D$ of
degree $d_D=2$ or $3$, we have
\begin{align*}
 \lambda^{-1}(I-P_\lambda)u
       &=\psi_0(\lambda^4\mathcal R)v,\qquad \lambda^{d_D-1}DP_\lambda u
       =(D\mathcal R^{-d_D/4})
         (\lambda^4\mathcal R)^{(d_D-1)/4}
               e^{-\lambda^4\mathcal R}v.
\end{align*}
Using this, the differential Riesz-transform bound
\eqref{wr:differential-riesz} and \eqref{w:vertical-SF}, applied
with $D=T$, $D=\partial_{X_i}\partial_{X_j}$ and
$D=T\partial_{X_j}$, show that the left side of
\eqref{w:smoothing-input} is $\lesssim\|\mathcal R^{1/4}u\|_2^2$.
Proposition~\ref{wr:energy-comparison} bounds this by a constant
multiple of $E(u)^2$.
Approximation on the common core in Lemma~\ref{wr:common-core}
and closedness of the differential operators identify the
$\mathrm L^2$ limits at each fixed scale. Fatou's lemma then
extends the inequality to every $u\in\V$.
\end{proof}

\section{Localization and the smoothed principal part}
\label{sec:localization}

We now use the geometry, energy space, and realization established in
Sections~\ref{sec:geometry-energy}--\ref{sec:variational-realization}.
All norms and averages in this section use Lebesgue measure. Put
\[
 R_\lambda=(I+\lambda^2H_A)^{-1},\qquad
 \Theta_\lambda F=-\lambda R_\lambda\operatorname{div}_X(AF),\qquad
 Q_\lambda u=\lambda H_AR_\lambda u.
\]
The divergence resolvent is understood in the variational sense of
Lemma~\ref{lem:scaled-resolvent}. In particular, $\Theta_\lambda$
includes multiplication by $A$.
Lemma~\ref{lem:scaled-resolvent} and the coefficient bound give
\begin{equation}\label{wt:resolvent-input}
 \begin{split}
 \|R_\lambda f\|_2+\lambda E(R_\lambda f)
 &\lesssim\|f\|_2,\\
 \|\Theta_\lambda F\|_2+\lambda E(\Theta_\lambda F)
 &\lesssim\|F\|_2,
 \qquad \|R_\lambda\|_{2\to2}\le1.
 \end{split}
\end{equation}

\subsection{Dyadic geometry and averaging}
Fix a system of dyadic cubes
$\mathcal D=\bigcup_{k\in\mathbb Z}\mathcal D_k$
adapted to $(G,d,\d z)$, constructed as in
\cite{Christ,Korvenpaa}, with a fixed sufficiently small
parameter $\kappa\in(0,1)$.
Each generation $\mathcal D_k$ partitions $G$ up to a set
of measure zero, and any two cubes are either disjoint or
one contains the other.
For $Q\in\mathcal D_k$, write $\ell(Q):=\kappa^k$.
There is a point $z_Q\in Q$ such that
\[
 B(z_Q,c\ell(Q))\subset Q\subset B(z_Q,C\ell(Q)),
\]
where $c,C>0$ depend only on the fixed geometry.
Consequently, \eqref{eq:kinetic-volume} gives
$|Q|\sim\ell(Q)^{\mathfrak q}$. The cubes have a small-boundary
exponent $\eta>0$ satisfying
\begin{equation}\label{aux:small-boundary}
 |\{z\in Q:d(z,Q^c)<a\ell(Q)\}|
 \lesssim a^\eta|Q|\qquad(0<a<1).
\end{equation}
We use open representatives whose boundaries have measure zero.
For locally integrable $F$, write $F_Q=|Q|^{-1}\int_QF\d z$.
For each $\lambda>0$, let $k\in\mathbb Z$ be the unique integer
such that
\begin{equation}\label{aux:scale-convention}
 \kappa^{k+1}<\lambda\le\kappa^k.
\end{equation}
Given $\lambda$ and $k$, for a locally integrable function $F$, define its dyadic
average at scale $\lambda$ by
\[
 \mathsf E_\lambda F(z):=F_Q
 :=\frac1{|Q|}\int_Q F(w)\d w,
 \qquad z\in Q,\quad Q\in\mathcal D_k.
\]
For vector-valued functions, the average is taken componentwise.
Since $\mathcal D_k$ partitions $G$ up to a null set,
this defines $\mathsf E_\lambda F$ almost everywhere.
We set it equal to zero on the exceptional set.
The scale convention ensures that
\begin{equation}\label{averagingcubes}
 \lambda\le\ell(Q)<\kappa^{-1}\lambda,
\end{equation}
so the averaging cubes have intrinsic scale comparable to $\lambda$.
These averaging operators are orthogonal projections on
$\mathrm L^2(G)$, and
\begin{equation}\label{wt:averages-bounded}
 \|\mathsf E_\lambda F\|_2\le\|F\|_2.
\end{equation}
The associated dyadic maximal operator
$M_{\!d}F(z):=\sup_{Q\ni z}\fint_Q|F(w)|\d w$ satisfies
$\|M_{\!d}F\|_2\le2\|F\|_2$ by the martingale maximal inequality.

\subsection{Transport localization} The unitary transport flow and its spectral difference formula give
\begin{equation}\label{wt:transport-differences}
 \||T|^{1/2}u\|_2^2
   =\frac1{2\pi}\int_{\R}\|e^{sT}u-u\|_2^2
                                      \frac{\d s}{|s|^2}.
\end{equation}
Let $\chi$ be a bounded smooth function with bounded
$\nabla_X\chi$ and $T\chi$. The identity
\[
 e^{sT}(\chi u)-\chi u
 =(e^{sT}\chi)(e^{sT}u-u)
       +(e^{sT}\chi-\chi)u
\]
separates the transport difference of $u$ from the additional
difference of the cutoff. Along the transport flow,
\begin{equation}\label{Transportest}
 \|e^{sT}\chi-\chi\|_\infty
 \le\min\{2\|\chi\|_\infty,|s|\|T\chi\|_\infty\}.
\end{equation}
Indeed, the uniform bound follows from
$\|e^{sT}\chi\|_\infty=\|\chi\|_\infty$ and the triangle
inequality. For the linear bound, the fundamental theorem
of calculus along each transport trajectory gives
\[
 (e^{sT}\chi)(z)-\chi(z)
 =\int_0^s(e^{rT}(T\chi))(z)\d r,
\]
whose absolute value is at most $|s|\|T\chi\|_\infty$. Splitting the difference integral at the scale where the two
bounds in \eqref{Transportest} agree gives
\[
 \int_{\R}
 \frac{\|e^{sT}\chi-\chi\|_\infty^2}{|s|^2}\d s
 \lesssim\|\chi\|_\infty\|T\chi\|_\infty.
\]
Applying the triangle inequality in
$\mathrm L^2(\R,\d s/|s|^2;\mathrm L^2(G))$ and using
\eqref{wt:transport-differences}, we obtain
\begin{equation}\label{Multiplierbound1}
 \||T|^{1/2}(\chi u)\|_2
 \lesssim \|\chi\|_\infty\||T|^{1/2}u\|_2
   +\bigl(\|\chi\|_\infty\|T\chi\|_\infty\bigr)^{1/2}\|u\|_2.
\end{equation}
Together with the horizontal product rule
\begin{equation}\label{Multiplierbound2}
 \|\nabla_X(\chi u)\|_2
 \le\|\chi\|_\infty\|\nabla_Xu\|_2
       +\|\nabla_X\chi\|_\infty\|u\|_2,
\end{equation}
this proves that multiplication by $\chi$ is bounded on $\V$.

Suppose now that $\chi$ is real-valued. Since the transport field
has zero divergence, integration by parts gives
\begin{equation}\label{wt:transport-cutoff}
 \operatorname{Re}\langle Tu,\chi^2u\rangle
       =-\frac12\int T(\chi^2)|u|^2\d z,
       \qquad u\in\V,
\end{equation}
where $\langle Tu,\chi^2u\rangle$ is understood as the
half-transport form pairing $\mathfrak t(u,\chi^2u)$. To justify this identity, first take
$u\in\mathrm C_0^\infty(G)$. Since $\chi$ is real-valued
and
\[
 2\operatorname{Re}\bigl((Tu)\overline u\bigr)=T(|u|^2),
\]
we obtain
\[
 \begin{split}
 2\operatorname{Re}\mathfrak t(u,\chi^2u)
 &=\int_G\chi^2T(|u|^2)\d z=\int_G T(\chi^2|u|^2)\d z
       -\int_G T(\chi^2)|u|^2\d z=-\int_G T(\chi^2)|u|^2\d z.
 \end{split}
\]
The last equality uses that
$T=\partial_t-X\cdot\nabla_Y$ is divergence-free with
respect to Lebesgue measure. Indeed, the integral of
$Tg$ vanishes for every smooth compactly supported $g$,
by integration by parts in $t$ and $Y$.
Here $g=\chi^2|u|^2$ is compactly supported because $u$ is.
This proves \eqref{wt:transport-cutoff} for smooth
compactly supported functions. For general $u\in\V$, choose
$u_n\in\mathrm C_0^\infty(G)$ with $u_n\to u$ in $\V$,
as provided by Lemma~\ref{lem:form-density}.
The multiplier bounds in \eqref{Multiplierbound1} and \eqref{Multiplierbound2} give
\[
 \chi^2u_n\longrightarrow\chi^2u
 \quad\text{in }\V.
\]
Consequently, continuity of $\mathfrak t$ on
$\V\times\V$ yields
\[
 \begin{split}
 |\mathfrak t(u_n,\chi^2u_n)
             -\mathfrak t(u,\chi^2u)|
 &\le
 |\mathfrak t(u_n-u,\chi^2u_n)|
 +|\mathfrak t(u,\chi^2(u_n-u))|\\
 &\quad\lesssim_\chi
 \|u_n-u\|_\V
       \bigl(\|u_n\|_\V+\|u\|_\V\bigr)
 \longrightarrow0.
 \end{split}
\]
On the other side, boundedness of $T(\chi^2)$ and
Cauchy--Schwarz give
\[
 \left|\int_G T(\chi^2)
          \bigl(|u_n|^2-|u|^2\bigr)\d z\right|
 \le
 \|T(\chi^2)\|_\infty
 \bigl(\|u_n\|_2+\|u\|_2\bigr)\|u_n-u\|_2
 \longrightarrow0.
\]
Passing to the limit proves
\eqref{wt:transport-cutoff} for every $u\in\V$,
with the left-hand side understood through the
half-transport form pairing $\mathfrak t$.
Thus the transport term, although skew-adjoint globally,
contributes the explicit cutoff term after localization.

\subsection{Caccioppoli estimates} We now use the preceding localization results to derive
a local energy estimate, which we then iterate over
concentric kinetic balls.

Let $u\in\V$ satisfy
$(I+\lambda^2H_A)u=0$ in the variational sense in an open set
$O$, and let $\chi\in\mathrm C_0^\infty(O)$ be real-valued.
As proved above, the multiplier property and density make $\chi^2u$ an
admissible test. Expanding its horizontal gradient and using
\eqref{wt:transport-cutoff}, we obtain
\[
 \begin{split}
 \|\chi u\|_2^2
 +\lambda^2\operatorname{Re}
       \langle A\nabla_Xu,\chi^2\nabla_Xu\rangle
 &=-2\lambda^2\operatorname{Re}
       \langle A\nabla_Xu,\chi u\nabla_X\chi\rangle+\lambda^2\int\chi T\chi\,|u|^2\d z.
 \end{split}
\]
Uniform accretivity controls the horizontal term on the left.
The coefficient bound and Young's inequality give
\[
 2\bigl|\langle A\nabla_Xu,\chi u\nabla_X\chi\rangle\bigr|
 \le\frac{a_0}{2}\|\chi\nabla_Xu\|_2^2
       +C\int|\nabla_X\chi|^2|u|^2\d z.
\]
Absorbing the first term yields the local energy estimate
\begin{equation}\label{wt:caccioppoli}
 \|\chi u\|_2^2+\lambda^2\|\chi\nabla_Xu\|_2^2
 \lesssim\lambda^2\int
       (|\nabla_X\chi|^2+|\chi T\chi|)|u|^2\d z.
\end{equation}
In particular, the full transport derivative of the cutoff
is retained in this estimate.

We now apply \eqref{wt:caccioppoli} on concentric kinetic balls.
Suppose that $B(z_0,2r)\subset O$. We require a cutoff
$\chi\in\mathrm C_0^\infty(B(z_0,2r))$ satisfying
$0\le\chi\le1$ and $\chi=1$ on $B(z_0,r)$, with derivative
bounds determined by the intrinsic scale $r$ and independent
of the centre $z_0$. To construct such cutoffs, write
$z_0=(X_0,Y_0,t_0)$ and introduce the frozen coordinates
\[
 (\xi,\eta,\tau)
 =z_0^{-1}\circ z
 =\bigl(X-X_0,\,
        Y-Y_0+(t-t_0)X_0,\,
        t-t_0\bigr).
\]
In these coordinates, $\nabla_X=\nabla_\xi$,  $T=\partial_\tau-\xi\cdot\nabla_\eta$. Left translation and intrinsic dilation of a fixed smooth cutoff
therefore give cutoffs at radius $r$ satisfying $|\nabla_X\chi|\lesssim r^{-1}$, $|T\chi|\lesssim r^{-2}$. For cutoffs between two concentric balls, the constant also
depends on their fixed radius ratio.

We next record explicitly how \eqref{wt:caccioppoli} is iterated.
Suppose that $u$ solves the homogeneous resolvent equation in
$B(z_0,2r)$, and fix an integer $M\ge1$. Set
\[
 r_j=r\bigl(1+{j}/{M}\bigr),
 \qquad j=0,\ldots,M.
\]
For $j=0,\ldots,M-1$, choose cutoffs equal to one on
$B(z_0,r_j)$ and supported in $B(z_0,r_{j+1})$.
Their horizontal and transport derivatives satisfy
$|\nabla_X\chi|\lesssim_M r^{-1}$ and $|T\chi|\lesssim_M r^{-2}$,
respectively. Taking square roots in \eqref{wt:caccioppoli} gives
\begin{equation}\label{uubound}
 \|u\|_{\mathrm L^2(B(z_0,r_j))}
 +\lambda\|\nabla_Xu\|_{\mathrm L^2(B(z_0,r_j))}
 \lesssim_M\frac{\lambda}{r}
          \|u\|_{\mathrm L^2(B(z_0,r_{j+1}))}.
\end{equation}
For $j=1,\ldots,M-1$, discard the gradient term and
iterate the resulting $\mathrm L^2$ estimate. This gives
\begin{equation}\label{uubound+}
 \|u\|_{\mathrm L^2(B(z_0,r_1))}
 \lesssim_M
 \left(\frac{\lambda}{r}\right)^{M-1}
 \|u\|_{\mathrm L^2(B(z_0,r_M))}.
\end{equation}
When $M=1$, this step is simply the identity.
Combining the bound in \eqref{uubound+} with the estimate in \eqref{uubound} for $j=0$,
retaining both terms on its left-hand side, and using
$r_0=r$ and $r_M=2r$, we obtain
\begin{equation}\label{wt:iterated-caccioppoli}
 \|u\|_{\mathrm L^2(B(z_0,r))}
 +\lambda\|\nabla_Xu\|_{\mathrm L^2(B(z_0,r))}
 \lesssim_M\left(\frac{\lambda}{r}\right)^M
          \|u\|_{\mathrm L^2(B(z_0,2r))}.
\end{equation}
Thus each additional localization contributes a factor
$\lambda/r$, with constants allowed to depend on the
number of iterations.

\subsection{Off-diagonal estimates}

We now combine the iterated Caccioppoli estimate \eqref{wt:iterated-caccioppoli} with the
global resolvent bounds to obtain off-diagonal decay estimates
between separated sets. For measurable sets $E,F\subset G$, write
\[
 d(E,F):=\inf\{d(z,w):z\in E,\ w\in F\},
\]
with the convention $\inf\varnothing=\infty$.

\begin{lem}\label{lem:off-diagonal}
For every $N>0$, $\lambda>0$ and measurable sets $E,F\subset G$,
\begin{equation}\label{wt:off-diagonal}
 \begin{gathered}
 \|1_E\mathcal U_\lambda1_F\|_{2\to2}
 \lesssim_N \bigl(1+d(E,F)/\lambda\bigr)^{-N},\qquad \mathcal U_\lambda\in
 \{R_\lambda,\Theta_\lambda,\lambda\nabla_XR_\lambda,
                         \lambda\nabla_X\Theta_\lambda\}.
 \end{gathered}
\end{equation}
The operator norms are taken between the corresponding scalar
or vector-valued $\mathrm L^2$ spaces. The implicit constant
depends only on $N$, $m$, $a_0$, $a_1$ and the fixed geometric
choices.
\end{lem}

\begin{proof}
The assertion is immediate if either set is empty. Otherwise,
put $D:=d(E,F)$. If $D\le\lambda$, the estimate follows from
the global bounds \eqref{wt:resolvent-input}. We therefore
assume that $D>\lambda$. Set $r=D/4$ and choose a maximal family $\{z_\nu\}\subset E$
with $d(z_\nu,z_\mu)\ge r$ whenever $\nu\ne\mu$.
Maximality ensures that the balls $B(z_\nu,r)$ cover $E$.
The doubling property and the quasi-triangle inequality give
bounded overlap of their enlargements, so
\[
 E\subset\bigcup_\nu B(z_\nu,r),
 \qquad
 \sum_\nu 1_{B(z_\nu,2r)}\lesssim1.
\]
Since $z_\nu\in E$ and $2r<D$, every enlarged ball
$B(z_\nu,2r)$ is disjoint from $F$.  For $f\in\mathrm L^2(G)$ and
$g\in\mathrm L^2(G;\C^m)$, consider either
\[
 u=R_\lambda(1_Ff)
 \qquad\text{or}\qquad
 u=\Theta_\lambda(1_Fg).
\]
In both cases, the resolvent equation is homogeneous on every
$B(z_\nu,2r)$. Indeed, the scalar forcing $1_Ff$ vanishes
there. In the second case, the forcing is
$-\lambda\operatorname{div}_X(A1_Fg)$, and its variational
pairing with a test function $v$ supported in such a ball is
\[
 \langle-\lambda\operatorname{div}_X(A1_Fg),v\rangle
 =\lambda\langle A1_Fg,\nabla_Xv\rangle=0.
\]
This uses only the support of $1_Fg$ and requires no
differentiation of $A$. Fix an integer $M\ge1$. Apply
\eqref{wt:iterated-caccioppoli} on each ball, square the
estimate, and sum over $\nu$. The covering and bounded
overlap give
\[
 \begin{split}
 \|1_Eu\|_2^2+\lambda^2\|1_E\nabla_Xu\|_2^2
 &\le\sum_\nu\left(
       \|u\|_{\mathrm L^2(B(z_\nu,r))}^2
       +\lambda^2
        \|\nabla_Xu\|_{\mathrm L^2(B(z_\nu,r))}^2
       \right)\\
 &\lesssim_M
       \left(\frac{\lambda}{D}\right)^{2M}
       \sum_\nu\|u\|_{\mathrm L^2(B(z_\nu,2r))}^2\lesssim_M
       \left(\frac{\lambda}{D}\right)^{2M}\|u\|_2^2.
 \end{split}
\]
Here the fixed factor arising from $r=D/4$ is absorbed into
the implicit constant. Taking square roots and applying
\eqref{wt:resolvent-input} bounds the last norm by that of
the input. Choose an integer $M\ge\max\{1,N\}$. Since
$D>\lambda$,
\[
 \bigl({\lambda}/{D}\bigr)^M
 \le\bigl({\lambda}/{D}\bigr)^N
 \lesssim_N(1+D/\lambda)^{-N}.
\]
This proves \eqref{wt:off-diagonal} for all four operator
families.
\end{proof}

\begin{rem}
Lemma~\ref{lem:off-diagonal} asserts no off-diagonal estimate
for the nonlocal operator $|T|^{1/2}R_\lambda$, and no such
estimate is needed in the subsequent arguments.
This parallels the parabolic argument in
\cite{AtaeiEgertNystrom}, where no analogous off-diagonal
estimate involving the half-time derivative $D_t^{1/2}$
is required.
\end{rem}

\subsection{Local Poincar\'e estimates and principal parts} We first establish a local Poincar\'e estimate adapted to the
kinetic geometry. We then use the off-diagonal estimates to
define the principal coefficients associated with
$\Theta_\lambda$ and prove their local $\mathrm L^2$ bounds.
These results will be used to approximate $\Theta_\lambda$
by its principal part on intrinsically smoothed gradients.

\begin{rem}\label{Cara}
In the proof of the local Poincar\'e estimate we will use the auxiliary Carnot--Carath\'eodory distance $d_c$
generated by $\partial_{X_1},\ldots,\partial_{X_m},T$, assigning
degree one to each of these fields. An admissible absolutely
continuous curve
$\gamma(\theta)=(X(\theta),Y(\theta),t(\theta))$, $0\le\theta\le1$,
satisfies
\[
 \dot\gamma(\theta)
 =\sum_{j=1}^m a_j(\theta)
       \partial_{X_j}\big|_{\gamma(\theta)}
       +b(\theta)T\big|_{\gamma(\theta)}
\]
for measurable real controls $a(\theta)\in\R^m$ and
$b(\theta)\in\R$ with integrable control norm.
In coordinates, this means
\[
 \dot X=a,\qquad \dot t=b,\qquad \dot Y=-bX.
\]
Its control length is
\[
 L_c(\gamma)
 :=\int_0^1\bigl(|a(\theta)|^2+|b(\theta)|^2\bigr)^{1/2}
                                                      \d\theta.
\]
We define $d_c(z,w)$ as the infimum of these lengths over
admissible curves joining $z$ to $w$. Both signs of the
controls are allowed. The commutators
$[\partial_{X_j},T]=-\partial_{Y_j}$ supply the remaining
coordinate directions, and all further commutators vanish.
Thus these fields generate a step-two
stratification of the same underlying group $G$, with the
$Y$-directions having degree two. The associated auxiliary
dilations are therefore
\[
 \delta_s^c(X,Y,t)=(sX,s^2Y,st).
\]
Indeed, if $\gamma$ has controls $(a,b)$, then
$\delta_s^c\circ\gamma$ has controls $(sa,sb)$, as follows
directly from the coordinate equations above. Its control
length is consequently $sL_c(\gamma)$. Taking infima and
using the inverse dilation gives
\[
 d_c(\delta_s^c z,\delta_s^c w)=s\,d_c(z,w).
\]
Moreover, the generating fields are left-invariant, so left
translation preserves control lengths and hence $d_c$.
Writing $B_c(z,s):=\{w\in G:d_c(z,w)<s\}$, we obtain
\[
 B_c(z,s)=z\circ\delta_s^c\bigl(B_c(0,1)\bigr).
\]
The Jacobian of $\delta_s^c$ is $s^{m+2m+1}=s^{3m+1}$, so
\[
 |B_c(z,s)|=|B_c(0,1)|s^{3m+1},
 \qquad
 |B_c(z,2s)|=2^{3m+1}|B_c(z,s)|.
\]
In particular, Lebesgue measure is doubling for $d_c$.
This auxiliary geometry is used below at unit kinetic scale.
The passage to arbitrary kinetic scales uses the original
dilations $\delta_\lambda$.
\end{rem}

\begin{lem}\label{wt:poincare}
For every kinetic ball $B_\lambda=B(z_0,\lambda)$ and smooth $F$, with
$F_{B_\lambda}=|B_\lambda|^{-1}\int_{B_\lambda}F\d z$, one has
\begin{equation}\label{wt:poincare-estimate}
 \int_{B_\lambda}|F-F_{B_\lambda}|^2\d z
 \lesssim\lambda^2\int_{CB_\lambda}|\nabla_XF|^2\d z
       +\lambda^4\int_{CB_\lambda}|TF|^2\d z.
\end{equation}
Here $CB_\lambda:=B(z_0,C\lambda)$ for a fixed geometric
constant $C\ge1$. The implicit constant depends only on the
fixed geometry and is independent of $z_0$, $\lambda$ and $F$.
\end{lem}

\begin{proof}
We first prove the estimate on the kinetic unit ball
$B_1:=B(0,1)$ and then use left translation and kinetic dilation
to obtain the general case. For the unit-scale argument we use
the auxiliary distance $d_c$ and its balls $B_c$ from
Remark~\ref{Cara}. Choose a fixed radius $R_0>0$ such that
$B_1\subset D:=B_c(0,R_0)$. Such a choice is possible because
the kinetic unit ball has compact closure and $d_c$ induces
the Euclidean topology. Let $D^*:=B_c(0,aR_0)$, where $a>1$
is the fixed enlargement factor in the representation formula
\eqref{CarPoin} below. Since $D^*$ is bounded, there is also
a fixed constant $C_0$ such that
\[
 B_1\subset D\subset D^*\subset B(0,C_0).
\]
These choices depend only on the fixed geometries.

The last assertion of
\cite[Proposition~2.12]{FranchiLuWheeden}, applied to the
bracket-generating family
$\partial_{X_1},\ldots,\partial_{X_m},T$
and its control distance $d_c$, gives the pointwise estimate
\begin{equation}\label{CarPoin}
 |F(z)-F_D|
 \lesssim\int_{D^*}
       \frac{d_c(z,v)}{|B_c(z,d_c(z,v))|}
       \bigl(|\nabla_XF(v)|+|TF(v)|\bigr)\d v,
 \qquad z\in D.
\end{equation}
We first apply the cited result on a sufficiently small auxiliary
ball centred at the identity. An auxiliary dilation $\delta_s^c$
maps this ball onto $D$ and transfers the estimate to $D$. Under these dilations, all generating fields have degree one
and Haar measure scales with degree $3m+1$. The derivative,
kernel and measure therefore contribute the factors $s$,
$s^{3m}$ and $s^{-(3m+1)}$, respectively, which cancel.
The average and the fixed enlarged ball transform accordingly.
Thus the same representation estimate holds on $D$, with
constants independent of $F$ and of the centre and radius of $D$. For complex-valued functions,
apply the real-valued formula to the real and imaginary parts.
To estimate the integral on the right in \eqref{CarPoin}, put
\[
 h:=1_{D^*}\bigl(|\nabla_XF|+|TF|\bigr),
 \qquad
 M_ch(z):=\sup_{s>0}\frac1{|B_c(z,s)|}
                              \int_{B_c(z,s)}|h(v)|\d v.
\]
For $z\in D$ and $v\in D^*$, the triangle inequality gives
$d_c(z,v)\le R_*:=(1+a)R_0$. Set $r_j:=2^{-j}R_*$.
On the annulus
\[
 A_j(z):=\{v:r_{j+1}<d_c(z,v)\le r_j\},
\]
doubling gives
\[
 \frac{d_c(z,v)}{|B_c(z,d_c(z,v))|}
 \lesssim\frac{r_j}{|B_c(z,r_j)|}.
\]
Consequently,
\[
 \begin{split}
 \int_{D^*}
       \frac{d_c(z,v)}{|B_c(z,d_c(z,v))|}h(v)\d v
 &\lesssim\sum_{j\ge0}
       \frac{r_j}{|B_c(z,r_j)|}
                         \int_{B_c(z,r_j)}h(v)\d v\\
                         &\lesssim\sum_{j\ge0}r_jM_ch(z)
 =2R_*M_ch(z).
 \end{split}
\]
The $\mathrm L^2$ boundedness of the Hardy--Littlewood maximal
operator on the doubling space $(G,d_c,\d z)$ therefore yields
\[
 \int_D|F-F_D|^2\d z
 \lesssim R_*^2\|M_ch\|_2^2
 \lesssim\int_{D^*}
               \bigl(|\nabla_XF|^2+|TF|^2\bigr)\d z.
\]
Here the fixed factor $R_*^2$ has been absorbed into the constant.

We now transfer this estimate to the kinetic unit ball.
The average $F_{B_1}$ minimizes
$\int_{B_1}|F-c|^2\d z$ over constants $c$. Therefore, using the
inclusions above, we obtain
\[
 \begin{split}
 \int_{B_1}|F-F_{B_1}|^2\d z
 \le\int_{B_1}|F-F_D|^2\d z
 &\le\int_D|F-F_D|^2\d z\\
 &\lesssim\int_{B(0,C_0)}
               \bigl(|\nabla_XF|^2+|TF|^2\bigr)\d z.
 \end{split}
\]
This proves the required estimate at unit kinetic scale.
The same inequality holds for vector-valued $F$ by applying
the scalar estimate componentwise and summing.

Finally, let $B_\lambda=B(z_0,\lambda)$ and define
\[
 \widetilde F(w):=F(z_0\circ\delta_\lambda w),
 \qquad
 \delta_\lambda(X,Y,t)=(\lambda X,\lambda^3Y,\lambda^2t).
\]
The map $w\mapsto z_0\circ\delta_\lambda w$ sends $B_1$
onto $B_\lambda$ and $B(0,C_0)$ onto
$B(z_0,C_0\lambda)$. Left invariance of the fields and
their kinetic degrees give
\[
 \begin{split}
 \nabla_X\widetilde F(w)
   &=\lambda(\nabla_XF)(z_0\circ\delta_\lambda w),\qquad T\widetilde F(w)
   =\lambda^2(TF)(z_0\circ\delta_\lambda w).
 \end{split}
\]
Moreover, the change of variables has Jacobian
$\lambda^{\mathfrak q}$, where $\mathfrak q=4m+2$, and hence
$\widetilde F_{B_1}=F_{B_\lambda}$. Apply the unit-scale
inequality to $\widetilde F$ and multiply by $\lambda^{\mathfrak q}$.
Changing variables in the two derivative integrals produces the factors $\lambda^2$ and
$\lambda^4$, respectively, and proves
\eqref{wt:poincare-estimate}. The auxiliary and kinetic balls have been compared only at
unit kinetic scale. The passage to arbitrary centres and
radii uses the original kinetic translations and dilations,
so all constants are independent of $z_0$, $\lambda$ and $F$.
\end{proof}

We next construct the principal coefficient of $\Theta_\lambda$
using the off-diagonal estimates. We then combine these estimates
with Lemma~\ref{wt:poincare} to obtain the principal-part
approximation. Define
$\gamma_\lambda\colon G\longrightarrow\C^{1\times m}$ to be a measurable
row-vector-valued function defined almost everywhere, as follows. Since constant vector fields do not belong to
$\mathrm L^2(G;\C^m)$, their images under $\Theta_\lambda$
are defined by localization. Let $\mathbf e_1,\ldots,\mathbf e_m$
be the standard basis of $\C^m$, and choose
$\chi_n\in\mathrm C_0^\infty(G)$ with $0\le\chi_n\le1$ and
$\chi_n=1$ on $B(0,r_n)$, where $r_n\uparrow\infty$. For each $\lambda>0$ and $j=1,\ldots,m$, define the scalar
function $\gamma_{\lambda,j}\in\mathrm L^2_{\mathrm{loc}}(G)$ by
\begin{equation}\label{limit1}
 \gamma_{\lambda,j}
 :=\lim_{n\to\infty}\Theta_\lambda(\chi_n\mathbf e_j)
 \quad\text{in }\mathrm L^2_{\mathrm{loc}}(G).
\end{equation}
After choosing measurable representatives, set $\gamma_\lambda(z)
 :=\bigl(\gamma_{\lambda,1}(z),\ldots,\gamma_{\lambda,m}(z)\bigr)$. Thus, for a constant column $e=(e_1,\ldots,e_m)^T\in\C^m$,
\begin{equation}\label{limit2}
 \gamma_\lambda(z)e
 =\sum_{j=1}^m\gamma_{\lambda,j}(z)e_j,
 \qquad
 \gamma_\lambda e
 =\lim_{n\to\infty}\Theta_\lambda(\chi_ne)
       \quad\text{in }\mathrm L^2_{\mathrm{loc}}(G).
\end{equation}
The row acts complex-linearly, with no conjugation in this
product. We use the notation
\[
 |\gamma_\lambda(z)|^2
 :=\sum_{j=1}^m|\gamma_{\lambda,j}(z)|^2.
\]

To justify the limits in \eqref{limit1} and \eqref{limit2},
fix an intrinsic dyadic cube $Q$ from the dyadic system
introduced above, and a scale $\lambda$ satisfying \(\kappa\ell(Q)<\lambda\le\ell(Q)\). Choose a sufficiently
large fixed geometric constant $C_*$ and set
\begin{align}\label{Setsuu}
 \mathcal A_0&:=B(z_Q,C_*\ell(Q)),\notag\\
 \mathcal A_j&:=B(z_Q,C_*2^j\ell(Q))
       \setminus B(z_Q,C_*2^{j-1}\ell(Q)),\qquad j\ge1.
 \end{align}
The choice of $C_*$ ensures that $Q\subset\mathcal A_0$
and $d(Q,\mathcal A_j)\gtrsim2^j\ell(Q)$ for $j\ge1$.
The off-diagonal estimate \eqref{wt:off-diagonal} and
the volume growth give
\begin{equation}\label{wt:gamma-annular-bound}
 \begin{split}
 \|1_Q\Theta_\lambda(e1_{\mathcal A_j})\|_2
 &\lesssim_N
 2^{-jN}|e|\,|\mathcal A_j|^{1/2}\lesssim_N
 2^{-j(N-\mathfrak q/2)}|e|\,|Q|^{1/2},
 \qquad j\ge0,
 \end{split}
\end{equation}
where the case $j=0$ follows from the global resolvent
bound. The same estimate holds with $1_{\mathcal A_j}$
replaced by $b1_{\mathcal A_j}$ for any bounded measurable
scalar function $b$, with the right-hand side multiplied
by $\|b\|_\infty$. 

Fix $N>\mathfrak q$. Summing the geometric majorant in
\eqref{wt:gamma-annular-bound} gives, for every integer
$J\ge0$,
\[
 \sum_{j=J}^\infty
 \|1_Q\Theta_\lambda(eb1_{\mathcal A_j})\|_2
 \lesssim_N
 \|b\|_\infty |e|\,|Q|^{1/2}
 2^{-J(N-\mathfrak q/2)}.
\]
In particular, the series
\[
 \sum_{j=0}^\infty
 1_Q\Theta_\lambda(e1_{\mathcal A_j})
\]
converges absolutely in $\mathrm L^2(Q)$. To apply the tail bound to the cutoff approximants,
take $b$ to be the difference of two cutoffs.
Their uniform boundedness controls $\|b\|_\infty$,
and $b$ vanishes on any fixed finite union of the
annuli once both cutoffs equal one on a sufficiently
large ball. The preceding estimate therefore shows
that the approximants form a Cauchy sequence in
$\mathrm L^2(Q)$. The same argument compares any two
admissible cutoff families, proving cutoff independence.
Since $Q$ is arbitrary, the defining limits in
\eqref{limit1} and \eqref{limit2} exist in
$\mathrm L^2_{\mathrm{loc}}(G)$ and are independent
of the chosen cutoffs.

Finally, summing \eqref{wt:gamma-annular-bound} over $j\ge0$
and passing to the cutoff limit also gives
\[
 \|1_Q\gamma_\lambda e\|_2
 \lesssim |e|\,|Q|^{1/2}.
\]
Applying this bound to the standard basis
$\mathbf e_1,\ldots,\mathbf e_m$ of $\C^m$ yields
\begin{equation}\label{wt:gamma-one-scale}
 \int_Q|\gamma_\lambda|^2\d z
 =\sum_{k=1}^m\|1_Q\gamma_\lambda\mathbf e_k\|_2^2
 \lesssim |Q|,
 \qquad \kappa\ell(Q)<\lambda\le\ell(Q).
\end{equation}
Moreover, this local estimate implies a uniform $\mathrm L^2$ bound
for $\gamma_\lambda\mathsf E_\lambda$.
Indeed, fix $\lambda>0$ and let $Q$ range over the dyadic
generation defining $\mathsf E_\lambda$, so that
$\kappa\ell(Q)<\lambda\le\ell(Q)$ for every such cube.
For $F\in\mathrm L^2(G;\C^m)$, we have
$\mathsf E_\lambda F=F_Q$ almost everywhere on $Q$.
Pointwise Cauchy--Schwarz and \eqref{wt:gamma-one-scale}
therefore give
\begin{equation}\label{wt:gamma-average-bounded}
 \begin{split}
 \|\gamma_\lambda\mathsf E_\lambda F\|_2^2
 =\sum_Q\int_Q|\gamma_\lambda(z)F_Q|^2\d z&\le\sum_Q|F_Q|^2\int_Q|\gamma_\lambda(z)|^2\d z\lesssim\sum_Q|Q||F_Q|^2
 \le\|F\|_2^2.
 \end{split}
\end{equation}
The last inequality follows from Jensen's inequality on
each cube and the fact that the cubes in this generation
partition $G$ up to a null set. The bound is uniform
in $\lambda>0$ and holds for every $F\in\mathrm L^2(G;\C^m)$.

We also record the measurability of the principal coefficients.
For each fixed cutoff $\chi_n$, the resolvent identity gives
strong $\mathrm L^2$ continuity of
$\lambda\mapsto\Theta_\lambda(\chi_n\mathbf e_j)$.
On every fixed kinetic ball, the same off-diagonal argument
used in \eqref{wt:gamma-annular-bound} gives an annular tail
bound that is uniform when $\lambda$ ranges over a compact
subinterval of $(0,\infty)$.
The approximants therefore converge locally in
$\mathrm L^2$, uniformly on such scale intervals.
Consequently, $\lambda\mapsto\gamma_\lambda\mathbf e_j$ is
continuous with values in $\mathrm L^2_{\mathrm{loc}}(G)$.
In particular, $\gamma_\lambda(z)$ admits a jointly
measurable representative in $(z,\lambda)$, which we
fix throughout.

With the principal coefficients thus defined, we turn to
the estimates needed for the principal-part approximation.

\begin{lem}\label{lem:local-principal-part}
For every $\lambda>0$ and
$F\in\mathscr S(G;\C^m)$, one has
\begin{equation}\label{wt:local-principal-part}
 \|\Theta_\lambda F
       -\gamma_\lambda\mathsf E_\lambda F\|_2^2
 \lesssim
 \lambda^2\|\nabla_XF\|_2^2
 +\lambda^4\|TF\|_2^2.
\end{equation}
Here the derivatives of $F$ are taken componentwise.
\end{lem}
\begin{proof} We begin with $F\in\mathscr S(G;\C^m)$.  For an averaging cube $Q$ at scale $\lambda$, as in
\eqref{averagingcubes}, we have
$\mathsf E_\lambda F=F_Q$ almost everywhere on $Q$.
By definition, $\gamma_\lambda F_Q$ represents the action
of $\Theta_\lambda$ on the constant vector field $F_Q$.
Linearity therefore gives the exact cancellation
\begin{equation}\label{AAAEst-}
 1_Q(\Theta_\lambda F-\gamma_\lambda\mathsf E_\lambda F)
 =1_Q\Theta_\lambda(F-F_Q).
\end{equation}
Since $F-F_Q$ need not belong to $\mathrm L^2(G)$, the
right-hand side is understood as the local $\mathrm L^2$
limit of $1_Q\Theta_\lambda(\chi_n(F-F_Q))$.
This limit follows from $\chi_nF\to F$ in $\mathrm L^2(G)$
and the defining local convergence
$\Theta_\lambda(\chi_nF_Q)\to\gamma_\lambda F_Q$. For the local principal-part estimate, set
\[
 B_j:=B(z_Q,C_*2^j\ell(Q)),
 \qquad j\ge0,
\]
using the sufficiently large fixed constant $C_*$ chosen in the definitions in \eqref{Setsuu}. Thus $Q\subset B_0$, and the annuli
$\mathcal A_j=B_j\setminus B_{j-1}$, $j\ge1$, satisfy
$d(Q,\mathcal A_j)\gtrsim2^j\ell(Q)$. We first compare the oscillation around the cube average
$F_Q$ with that around the ball average $F_{B_j}$.
Since $Q\subset B_j$, Cauchy--Schwarz gives
\[
 \begin{split}
 |F_Q-F_{B_j}|
 =\left|\frac1{|Q|}
          \int_Q(F-F_{B_j})\d z\right|&\le |Q|^{-1/2}
          \|F-F_{B_j}\|_{\mathrm L^2(Q)}\\
          &\le |Q|^{-1/2}
          \|F-F_{B_j}\|_{\mathrm L^2(B_j)}.
 \end{split}
\]
Using
$F-F_Q=(F-F_{B_j})+(F_{B_j}-F_Q)$
and the triangle inequality, we therefore obtain
 \begin{align}\label{AAAEst}
 \|F-F_Q\|_{\mathrm L^2(B_j)}
 &\le \|F-F_{B_j}\|_{\mathrm L^2(B_j)}
       +|B_j|^{1/2}|F_{B_j}-F_Q|\notag\\
 &\le
 \bigl(1+\bigl({|B_j|}/{|Q|}\bigr)^{1/2}\bigr)
       \|F-F_{B_j}\|_{\mathrm L^2(B_j)}\notag\\
       &\lesssim
 \bigl({|B_j|}/{|Q|}\bigr)^{1/2}
       \|F-F_{B_j}\|_{\mathrm L^2(B_j)},
 \end{align}
where the last step uses $|B_j|\ge|Q|$. Using \eqref{AAAEst-} and the off-diagonal estimate
\eqref{wt:off-diagonal}, followed by \eqref{AAAEst}
and the local Poincar\'e estimate
\eqref{wt:poincare-estimate}, we obtain
\[
 \begin{split}
 \|1_Q(\Theta_\lambda F
          -\gamma_\lambda\mathsf E_\lambda F)\|_2
 &\quad\lesssim_N
 \sum_{j\ge0}2^{-jN}
       \|F-F_Q\|_{\mathrm L^2(B_j)}\\
 &\quad\lesssim_N
 \sum_{j\ge0}2^{-jN}2^{j\mathfrak q/2}
       \|F-F_{B_j}\|_{\mathrm L^2(B_j)}\\
 &\quad\lesssim_N
 \sum_{j\ge0}2^{-jN}2^{j\mathfrak q/2}
 \left(
 2^j\lambda\|\nabla_XF\|_{\mathrm L^2(CB_j)}
 +2^{2j}\lambda^2\|TF\|_{\mathrm L^2(CB_j)}
 \right).
 \end{split}
\]
Here the term $j=0$ is controlled by the global resolvent
bound. In the last step, we used that the radius of $B_j$
is comparable to $2^j\lambda$, since
$\lambda\sim\ell(Q)$. The factors $2^j\lambda$ and
$2^{2j}\lambda^2$ therefore come from the horizontal
and transport terms, respectively, in the local
Poincar\'e estimate.

We now sum the local estimates over the cubes $Q$ in
the averaging generation. These cubes partition $G$
up to a set of measure zero. For each fixed $j$, the
corresponding enlarged balls $CB_j$ have overlap
multiplicity $\lesssim2^{j\mathfrak q}$. Consequently, for any scalar- or vector-valued
$g\in\mathrm L^2(G)$,
\begin{equation}\label{overlap}
 \sum_Q\|g\|_{\mathrm L^2(CB_j)}^2
 \lesssim 2^{j\mathfrak q}\|g\|_2^2.
\end{equation}
Taking the $\ell^2$ norm over $Q$, applying Minkowski's
inequality to the sum over $j$, and using the overlap
bound in \eqref{overlap} which contributes another factor $2^{j\mathfrak q/2}$, we obtain
\[
 \|\Theta_\lambda F-\gamma_\lambda\mathsf E_\lambda F\|_2
 \lesssim_N
 \sum_{j\ge0}2^{-jN}
 \left(
 2^{j(\mathfrak q+1)}\lambda\|\nabla_XF\|_2
 +2^{j(\mathfrak q+2)}\lambda^2\|TF\|_2
 \right).
\]
Choosing $N>\mathfrak q+2$, summing the geometric series,
and squaring gives \eqref{wt:local-principal-part}. \end{proof}

\subsection{The smoothed principal-part approximation}

We now combine the exact resolvent decomposition with
the local principal-part estimate and the Rockland
smoothing estimates.

\begin{lem}\label{lem:smoothed-principal-part}
For every $u\in\V$, one has
\begin{equation}\label{aux:smoothed-PP}
 \int_0^\infty
 \|Q_\lambda u-\gamma_\lambda\mathsf E_\lambda
                     \nabla_XP_\lambda u\|_2^2
                     \frac{\d\lambda}{\lambda}
 \lesssim E(u)^2.
\end{equation}
\end{lem}

\begin{proof}
Recall that $P_\lambda=e^{-\lambda^4\mathcal R}$, where
$\mathcal R$ is the Rockland operator introduced in
\eqref{Rockland}. The spectral theorem shows that $P_\lambda$
preserves $\mathrm D(\mathcal R^{1/4})=\V$ and commutes
with $\mathcal R^{1/4}$. Thus
\eqref{wr:energy-equivalence} gives
\[
 E(P_\lambda u)
 \lesssim\|\mathcal R^{1/4}P_\lambda u\|_2
 =\|P_\lambda\mathcal R^{1/4}u\|_2
 \le\|\mathcal R^{1/4}u\|_2
 \lesssim E(u),
 \qquad u\in\V.
\]
Moreover, since $T$ has homogeneous degree two,
\eqref{wr:differential-riesz} gives boundedness of
$T\mathcal R^{-1/2}$ on $\mathrm L^2(G)$. Combining this
with the spectral theorem yields
\[
 \begin{split}
 \|TP_\lambda u\|_2=\|(T\mathcal R^{-1/2})
           \mathcal R^{1/2}P_\lambda u\|_2\lesssim\|\mathcal R^{1/2}P_\lambda u\|_2\le\sup_{r>0}\bigl(r^{1/2}e^{-\lambda^4r}\bigr)\|u\|_2
 \lesssim\lambda^{-2}\|u\|_2.
 \end{split}
\]
The factorization is first justified on the common core
of Lemma~\ref{wr:common-core} and then extended by density
and closedness of $T$. In particular,
$P_\lambda u\in\V$ and $TP_\lambda u\in\mathrm L^2(G)$
for every $u\in\V$ and $\lambda>0$.

We now establish an exact decomposition of $Q_\lambda u$ valid for every
$u\in\V$. For $v\in\V$, we interpret $H_Av$ as the element of $\V'$
defined by
\[
 (H_Av)(w)=\mathfrak a_A(v,w),
 \qquad w\in\V.
\]
The variational extension of the resolvent satisfies
\begin{equation}\label{Identityuu}
 Q_\lambda v
 =\lambda^{-1}(v-R_\lambda v)
 =\lambda R_\lambda H_Av,
 \qquad v\in\V.
\end{equation}
On $\mathrm D(H_A)$, this is the usual commutation identity
$H_AR_\lambda v=R_\lambda H_Av$.
For general $v\in\V$, the expression $R_\lambda H_Av$
is understood through the variational extensions.
Applying the identity in \eqref{Identityuu} to $v=P_\lambda u$ and splitting
\[
 H_AP_\lambda u
 =TP_\lambda u
  -\operatorname{div}_X(A\nabla_XP_\lambda u)
 \quad\text{in }\V'
\]
we obtain
\begin{equation}\label{aux:Q-decomposition}
 Q_\lambda u
 =Q_\lambda(I-P_\lambda)u
       +\lambda R_\lambda TP_\lambda u
       +\Theta_\lambda\nabla_XP_\lambda u.
\end{equation}
This argument uses only $P_\lambda u\in\V$, interpreting
$H_AP_\lambda u$ as an element of $\V'$.
It does not require $P_\lambda u\in\mathrm D(H_A)$,
or equivalently $H_AP_\lambda u\in\mathrm L^2(G)$,
which need not follow from smoothing when $A$ is merely measurable.

To prove the square-function estimate, we now take
$u\in\mathscr S_0(G)$.
Since the Rockland heat operator
$P_\lambda=e^{-\lambda^4\mathcal R}$ preserves
$\mathscr S_0(G)$, we have
$\nabla_XP_\lambda u\in\mathscr S(G;\C^m)$
for $u\in\mathscr S_0(G)$. Hence Lemma~\ref{lem:local-principal-part} applies.
We will remove this restriction by density at the end
of the proof. The bounds
\[
 \|Q_\lambda\|_{2\to2}\le2/\lambda,
 \qquad
 \|R_\lambda\|_{2\to2}\le1
\]
and \eqref{aux:Q-decomposition} yield
\begin{align*}
 \|Q_\lambda u-\gamma_\lambda\mathsf E_\lambda
                              \nabla_XP_\lambda u\|_2
 &\le
 2\lambda^{-1}\|(I-P_\lambda)u\|_2
 +\lambda\|TP_\lambda u\|_2+
 \|\Theta_\lambda\nabla_XP_\lambda u
       -\gamma_\lambda\mathsf E_\lambda
                         \nabla_XP_\lambda u\|_2.
\end{align*}
We apply \eqref{wt:local-principal-part} with
$F=\nabla_XP_\lambda u$ to control the last term.
Squaring the preceding inequality and integrating
then gives
\begin{align}\label{aux:smoothed-PP-}
 \int_0^\infty
 \|Q_\lambda u-\gamma_\lambda\mathsf E_\lambda
                     \nabla_XP_\lambda u\|_2^2
                     \frac{\d\lambda}{\lambda}
 \lesssim&
 \int_0^\infty
 \Bigl(
 \|\lambda^{-1}(I-P_\lambda)u\|_2^2
 +\|\lambda TP_\lambda u\|_2^2
 \Bigr)\frac{\d\lambda}{\lambda}\notag\\
 &+
 \int_0^\infty
 \Bigl(
 \|\lambda\nabla_X^2P_\lambda u\|_2^2
 +\|\lambda^2T\nabla_XP_\lambda u\|_2^2
 \Bigr)\frac{\d\lambda}{\lambda}.
\end{align}
The four smoothing
estimates in \eqref{w:smoothing-input} now prove \eqref{aux:smoothed-PP} on
$\mathscr S_0(G)$.

To extend the estimate to an arbitrary $u\in\V$,
choose $u_n\in\mathscr S_0(G)$ converging to $u$
in $\V$, using Lemma~\ref{wr:common-core}.
For each fixed $\lambda>0$, boundedness of $Q_\lambda$
on $\mathrm L^2(G)$ gives
\[
 \|Q_\lambda(u_n-u)\|_2
 \le 2\lambda^{-1}\|u_n-u\|_2
 \longrightarrow0.
\]
Moreover, \eqref{wt:gamma-average-bounded} and boundedness
of $P_\lambda$ in the energy norm imply
\[
 \begin{split}
 \|\gamma_\lambda\mathsf E_\lambda
       \nabla_XP_\lambda(u_n-u)\|_2
 &\lesssim
 \|\nabla_XP_\lambda(u_n-u)\|_2\lesssim E(u_n-u)\longrightarrow0.
 \end{split}
\]
Thus the expression inside the square-function norm
converges in $\mathrm L^2(G)$ at every fixed scale.
Fatou's lemma, together with $E(u_n)\to E(u)$,
now proves \eqref{aux:smoothed-PP} for every $u\in\V$.
\end{proof}

\section{Removing the intrinsic smoothing}\label{sec:unsmoothing}

The purpose of this section is to remove $P_\lambda$ from
the smoothed principal-part estimate in Lemma~\ref{lem:smoothed-principal-part} by replacing
$\gamma_\lambda\mathsf E_\lambda\nabla_XP_\lambda u$
by $\gamma_\lambda\mathsf E_\lambda\nabla_Xu$.
Since $\mathsf E_\lambda^2=\mathsf E_\lambda$, the resulting
error decomposes as
\[
 \begin{split}
 \gamma_\lambda\mathsf E_\lambda
       (\nabla_XP_\lambda-\nabla_X)u
 &=\gamma_\lambda\mathsf E_\lambda
       [\nabla_X,P_\lambda]u+\gamma_\lambda\mathsf E_\lambda
       (P_\lambda-\mathsf E_\lambda)\nabla_Xu.
 \end{split}
\]
The first term arises because horizontal differentiation
does not commute with the intrinsic smoothing. The second term compares smoothing with dyadic averaging at the same scale.
By the one-scale boundedness of
$\gamma_\lambda\mathsf E_\lambda$, it suffices to establish
square-function estimates for $[\nabla_X,P_\lambda]u$
and $(P_\lambda-\mathsf E_\lambda)\nabla_Xu$. Consequently, the resulting approximation involves the horizontal gradient of
$u$ directly and provides the input for the localized testing
argument in the next section.

All norms, convolutions and dyadic averages in this section use
Lebesgue measure, which is the Haar measure on $G$. We continue to use the convention
\[
 (f*k)(z)=\int_G f(\zeta)k(\zeta^{-1}\circ z)\d\zeta.
\]
The quasi-norm $\rho$ from Subsection~\ref{subsec:kolmogorov-geometry}
is symmetric, and inversion preserves Haar measure.
Recall the scale kernels
$p_\lambda=h_{\lambda^4}$, where $h_t$ is the heat-time kernel
of $e^{-t\mathcal R}$.

The argument has two parts. Cancellation and almost orthogonality
control the commutator with horizontal differentiation. A separate
dyadic estimate compares the smooth approximation with conditional
expectation. Neither part uses a Carleson estimate.

\subsection{The convolution commutator kernel}

Recall that a distribution
kernel $K$ which is homogeneous of degree $-\mathfrak q$ and smooth away
from the identity is called a
\emph{type-zero kernel}.

\begin{lem}\label{lem:type-zero-commutator}
Let $\mathcal C$ be a bounded right-convolution operator on $\mathrm L^2(G)$ whose distribution kernel $K$ is of type zero.
Then, for every $f\in\mathrm L^2(G)$,
\begin{equation}\label{w:commutator-SF}
 \int_0^\infty
 \|[\mathcal C,P_\lambda]f\|_2^2\frac{\d\lambda}{\lambda}
 \lesssim\|f\|_2^2,
\end{equation}
where
$[\mathcal C,P_\lambda]=\mathcal CP_\lambda-P_\lambda\mathcal C$.
The implicit constant depends on $m$ and finitely many homogeneous
kernel seminorms of $K$.
\end{lem}

The proof of Lemma~\ref{lem:type-zero-commutator} is given in
Appendix~\ref{app:type-zero-commutator}.
Its main step is to establish cancellation, a finite
positive fractional moment, and integrable left and
right derivatives for the commutator kernel
\[
 \mathfrak h=p_1*K-K*p_1.
\]
These properties yield almost orthogonality between
different scales, from which the square-function
estimate follows by Schur's test.
Both left and right translation estimates enter because
convolution on $G$ is noncommutative.

We will use the following consequence of Lemma~\ref{lem:type-zero-commutator}.

\begin{lem}\label{lem:gradient-commutator}
For every $u\in\V$ we have
\begin{equation}\label{eq:gradient-unsmoothing}
 \int_0^\infty
 \|[\nabla_X,P_\lambda]u\|_2^2\frac{\d\lambda}{\lambda}
 \lesssim E(u)^2.
\end{equation}
\end{lem}
\begin{proof} Set
\[
 \mathcal C_j:=\partial_{X_j}\mathcal R^{-1/4},
 \qquad j=1,\ldots,m.
\]
Lemma~\ref{wr:zero-order-and-bip} shows that these are bounded
convolution operators with type-zero kernels. For these fixed
transforms, the required kernel seminorms depend only on $m$
and the fixed geometric choices.  For every $u$ in the common core $\mathscr S_0(G)$
of Lemma~\ref{wr:common-core}, we have $\partial_{X_j}u=\mathcal C_j\mathcal R^{1/4}u$ for $j=1,\ldots,m$. Since $P_\lambda$ commutes with powers of $\mathcal R$,
we consequently have
\[
 [\partial_{X_j},P_\lambda]u
 =[\mathcal C_j,P_\lambda]\mathcal R^{1/4}u.
\]
Applying Lemma~\ref{lem:type-zero-commutator}, summing over $j$, and
using \eqref{wr:energy-equivalence} gives
\[
 \int_0^\infty
 \|[\nabla_X,P_\lambda]u\|_2^2\frac{\d\lambda}{\lambda}
 \lesssim\|\mathcal R^{1/4}u\|_2^2
 \lesssim E(u)^2,
 \qquad u\in\mathscr S_0(G).
\]
To extend this estimate to $u\in\V$, choose
$u_n\in\mathscr S_0(G)$ converging to $u$ in $\V$.
The energy comparison and the spectral theorem give
$E(P_\lambda v)\lesssim E(v)$, uniformly in $\lambda$.
Together with $\mathrm L^2$ contractivity of $P_\lambda$,
this implies, for every fixed $\lambda>0$,
\[
 [\nabla_X,P_\lambda]u_n
 \longrightarrow[\nabla_X,P_\lambda]u
 \quad\text{in }\mathrm L^2(G;\C^m).
\]
Fatou's lemma now proves \eqref{eq:gradient-unsmoothing}
for every $u\in\V$.
\end{proof}

\subsection{Dyadic averages and smooth approximation}

The intrinsic dyadic cubes have the small-boundary property
\eqref{aux:small-boundary}, and $|Q|\sim\ell(Q)^{\mathfrak q}$.
Let $\mathsf E_j$ be conditional expectation at generation $j$.
These are orthogonal projections on $\mathrm L^2(G)$. Also
$\|P_\lambda\|_{2\to2}\le1$ by the spectral theorem.

\begin{lem}\label{lem:dyadic-unsmoothing}
For each $\lambda>0$, let $\mathsf E_\lambda=\mathsf E_j$,
where $j\in\mathbb Z$ is the unique integer satisfying
$\kappa^{j+1}<\lambda\le\kappa^j$.
Then, for every $f\in\mathrm L^2(G)$,
\begin{equation}\label{eq:dyadic-unsmoothing}
 \int_0^\infty
 \|(P_\lambda-\mathsf E_\lambda)f\|_2^2
 \frac{\d\lambda}{\lambda}
 \lesssim\|f\|_2^2.
\end{equation}
\end{lem}

\begin{proof}
Fix $j\in\mathbb Z$ and
$\kappa^{j+1}<\lambda\le\kappa^j$, so that
$\mathsf E_\lambda=\mathsf E_j$.
For $k\in\mathbb Z$, write
$\Delta_k:=\mathsf E_{k+1}-\mathsf E_k$.
We first establish the required off-scale estimates,
with constants independent of $j$, $k$ and $\lambda$.

If $k\ge j$, then $\mathsf E_j\Delta_k=0$ and $\Delta_kf$
has Haar mean zero on each cube $Q$ of generation $k$. For each fixed output point $z$, subtract the kernel value
$p_\lambda(z_Q^{-1}\circ z)$ corresponding to the centre
of the input cube $Q$. Since $\int_Q\Delta_kf\,\d \zeta=0$, we have
\begin{align}\label{Kernelest1}
 &\int_Q p_\lambda(\zeta^{-1}\circ z)\Delta_kf(\zeta)\d \zeta=
 \int_Q
 \bigl[p_\lambda(\zeta^{-1}\circ z)
       -p_\lambda(z_Q^{-1}\circ z)\bigr]
 \Delta_kf(\zeta)\d \zeta.
 \end{align}
When
$\kappa^k\le c_*\lambda$ for a sufficiently small fixed
$c_*>0$, Schwartz estimates using right-invariant derivatives give
\begin{align}\label{Kernelest2}
 |p_\lambda(\zeta^{-1}\circ z)-p_\lambda(z_Q^{-1}\circ z)|
 &\lesssim_N(\kappa^k/\lambda)^a \lambda^{-\mathfrak q}
                (1+d(\zeta,z)/\lambda)^{-N},\qquad \zeta\in Q,
\end{align}
with a fixed $a>0$ and arbitrarily large $N$.  The subscript records the possible dependence on $N$. Indeed, put $w=\zeta^{-1}\circ z$ and
$b=z_Q^{-1}\circ\zeta$. Then
$z_Q^{-1}\circ z=b\circ w$, while
$\rho(b)\lesssim\kappa^k$ by the diameter bound for $Q$.
The Schwartz bounds for $p_1=h_1$ and the scaling of
$p_\lambda$ give
\[
 |Z^Rp_\lambda(y)|
 \lesssim_N\lambda^{-\mathfrak q-v}
             (1+\rho(y)/\lambda)^{-N}
\]
for each homogeneous right-invariant basis field $Z^R$
of degree $v\in\{1,2,3\}$.
The corresponding mean-value inequality therefore yields
\begin{align}
 |p_\lambda(b\circ w)-p_\lambda(w)|
 &\lesssim_N
 \lambda^{-\mathfrak q}
 \sum_{v=1}^3
 \left(\frac{\rho(b)}{\lambda}\right)^v
 (1+\rho(w)/\lambda)^{-N}\lesssim_N
 \frac{\kappa^k}{\lambda}\lambda^{-\mathfrak q}
 (1+d(\zeta,z)/\lambda)^{-N}.
 \end{align}
Here the quasi-triangle inequality makes the decay weights
along the intervening translations comparable to that at $w$,
since those translations have size $O(\rho(b))$.
Choosing $c_*$ sufficiently small ensures
$\rho(b)/\lambda\le1$, so each higher-degree contribution
is bounded by $\rho(b)/\lambda$. For the right-invariant mean-value inequality used here, see
\cite[Remark~3.1.48(1), following Proposition~3.1.46]
{FischerRuzhanskyBook}. This proves \eqref{Kernelest2} with $a=1$. Writing $\mathcal D_k$ for the dyadic cubes of generation $k$,
define
\[
 K_{\lambda,k}(z,\zeta)
 :=\sum_{Q\in\mathcal D_k}1_Q(\zeta)
 \bigl[
 p_\lambda(\zeta^{-1}\circ z)
 -p_\lambda(z_Q^{-1}\circ z)
 \bigr].
\]
Only the cube containing $\zeta$ contributes to the sum for almost
every $\zeta\in G$.  When $\kappa^k\le c_*\lambda$, choosing $N>\mathfrak q$ and integrating
\eqref{Kernelest2}, using the volume bound
$|B(z,r)|\lesssim r^{\mathfrak q}$, gives the two bounds
required by Schur's test,
\begin{equation}\label{Schur}
 \sup_{z\in G}\int_G|K_{\lambda,k}(z,\zeta)|\d\zeta
 +\sup_{\zeta\in G}\int_G|K_{\lambda,k}(z,\zeta)|\d z
 \lesssim(\kappa^k/\lambda)^a.
\end{equation}
Since $\Delta_kf$ has mean zero on each cube of generation
$k$ and $\mathsf E_j\Delta_kf=0$ for $k\ge j$, Schur's test
therefore yields
\[
 \|(P_\lambda-\mathsf E_j)\Delta_kf\|_2
 \lesssim(\kappa^k/\lambda)^a\|\Delta_kf\|_2
 \lesssim\kappa^{a(k-j)}\|\Delta_kf\|_2,
\]
where the last step uses $\lambda>\kappa^{j+1}$. In the remaining range $\kappa^k>c_*\lambda$, the uniform
bound $\|p_\lambda\|_1=\|p_1\|_1$ gives
\[
 \|(P_\lambda-\mathsf E_j)\Delta_kf\|_2
 =\|P_\lambda\Delta_kf\|_2
 \lesssim\|\Delta_kf\|_2.
\]
Here $\kappa^{k-j}>c_*\kappa$, so this also implies the
desired estimate, with a fixed larger constant.
Combining the two ranges, we obtain
\begin{equation}\label{w:Haar-dyadic-fine}
 \|(P_\lambda-\mathsf E_j)\Delta_k f\|_2
 \lesssim\kappa^{a(k-j)}\|\Delta_kf\|_2,
 \qquad k\ge j.
\end{equation}

If $k<j$, then $\mathsf E_j\Delta_k=\Delta_k$, and
$\Delta_kf$ is constant on generation-$k+1$ cubes. For any
step function $v=\sum_Qc_Q1_Q$ in $\mathrm L^2(G)$ on cubes
of side $\ell$, the identity $\int_Gp_\lambda\d z=1$ and
Cauchy--Schwarz with $|p_\lambda|$ give
\[
 \begin{split}
 \|(P_\lambda-I)v\|_2^2
 &\lesssim\sum_Q|c_Q|^2\bigg[
   \int_Q\int_{Q^c}|p_\lambda(\zeta^{-1}\circ z)|\d\zeta\d z+\int_Q\int_{Q^c}|p_\lambda(z^{-1}\circ\zeta)|\d\zeta\d z\bigg].
 \end{split}
\]
These estimates use only the absolute value of the heat kernel.
For $z\in Q$, its Schwartz decay gives
\[
 \int_{Q^c}|p_\lambda(\zeta^{-1}\circ z)|\d\zeta
 \lesssim_N(1+d(z,Q^c)/\lambda)^{-N},
\]
and the same bound holds for the reversed kernel.
By the Haar small-boundary estimate,
\[
 |\{z\in Q:d(z,Q^c)\le a\ell\}|\lesssim a^\eta|Q|,
 \qquad 0<a<1.
\]
Sum the layers $d(z,Q^c)\le\lambda$ and
$2^l\lambda<d(z,Q^c)\le2^{l+1}\lambda$, with a fixed $N>\eta$.
Use the small-boundary estimate while $2^{l+1}\lambda<\ell$
and the trivial measure bound $|Q|$ for the remaining layers.
The resulting geometric sums bound both double integrals by
a fixed multiple of $(\lambda/\ell)^\eta|Q|$.
Taking $\ell\sim\kappa^{k+1}$ gives
\begin{equation}\label{w:Haar-dyadic-coarse}
 \|(P_\lambda-\mathsf E_j)\Delta_kf\|_2
  \lesssim\kappa^{\eta(j-k-1)/2}\|\Delta_kf\|_2,\qquad k<j.
\end{equation}
Comparable-scale cases are covered by the uniform operator bounds.
In particular, using \eqref{w:Haar-dyadic-fine} and \eqref{w:Haar-dyadic-coarse}  we conclude, for some $c>0$,
\begin{equation}\label{w:Haar-offscale-operator}
 \|(P_\lambda-\mathsf E_j)\Delta_k\|_{2\to2}
                         \lesssim\kappa^{c|j-k|}.
\end{equation}

The martingale differences are mutually orthogonal. Moreover,
$\mathsf E_j f\to f$ as $j\to\infty$, and
$\mathsf E_j f\to0$ as $j\to-\infty$, in $\mathrm L^2$.
For the first assertion, the cube diameters tend to zero. Approximate
by continuous compactly supported functions and use contractivity
of $\mathsf E_j$. For the second, take first $f$ bounded and compactly
supported. Since all generation-$j$ cubes have volume comparable to
$\kappa^{j\mathfrak q}$,
\[
 \|\mathsf E_j f\|_2^2
 =\sum_{Q:\,\ell(Q)=\kappa^j}\frac1{|Q|}
             \biggl|\int_Qf\d z\biggr|^2
 \lesssim\kappa^{-j\mathfrak q}\|f\|_1^2\longrightarrow0
 \qquad(j\to-\infty).
\]
Density and contractivity give the assertion for general $f$.
Consequently the telescope and orthogonality yield
\begin{equation}\label{w:martingale-square}
 f=\sum_k\Delta_kf\quad\text{in }\mathrm L^2,
 \qquad \sum_k\|\Delta_kf\|_2^2=\|f\|_2^2.
\end{equation}

Next suppose that $f$ is a finite sum of martingale
differences. Minkowski's inequality and
\eqref{w:Haar-offscale-operator}, together with
$\Delta_k^2=\Delta_k$, give
\[
 \|(P_\lambda-\mathsf E_j)f\|_2
 \lesssim
 \sum_k\kappa^{c|j-k|}\|\Delta_kf\|_2,
 \qquad \kappa^{j+1}<\lambda\le\kappa^j.
\]
Since $0<\kappa<1$ and $c>0$, the sequence
$(\kappa^{c|n|})_{n\in\mathbb Z}$ belongs to
$\ell^1(\mathbb Z)$. Integrating over the scale intervals
and applying the discrete Young inequality therefore yields
\[
 \begin{split}
 \sum_j\int_{\kappa^{j+1}}^{\kappa^j}
     \|(P_\lambda-\mathsf E_j)f\|_2^2
     \frac{\d\lambda}{\lambda}
 &\lesssim
 \sum_j\biggl(
     \sum_k\kappa^{c|j-k|}\|\Delta_kf\|_2
 \biggr)^2\lesssim\sum_k\|\Delta_kf\|_2^2
 =\|f\|_2^2.
 \end{split}
\]
Here the fixed factor $\log(1/\kappa)$ from integration
over each scale interval is absorbed into the implicit
constant, and the last equality follows from orthogonality
of the martingale differences. For general $f\in\mathrm L^2(G)$, let
$f_n:=\sum_{|k|\le n}\Delta_kf$. The martingale endpoint
convergence gives $f_n\to f$ in $\mathrm L^2(G)$.
On every finite scale interval $[a,b]\subset(0,\infty)$,
the bound $\|P_\lambda-\mathsf E_\lambda\|_{2\to2}\le2$ gives
\[
 \int_a^b
 \|(P_\lambda-\mathsf E_\lambda)(f_n-f)\|_2^2
                      \frac{\d\lambda}{\lambda}
 \le4\log(b/a)\|f_n-f\|_2^2\longrightarrow0.
\]
Passing to the limit in the estimate for $f_n$ proves the
same bound on $[a,b]$, uniformly in its endpoints.
Finally let $a\downarrow0$ and $b\uparrow\infty$ and use
monotone convergence to obtain \eqref{eq:dyadic-unsmoothing}.
\end{proof}

\subsection{The unsmoothed principal part}

In the following, the two preceding lemmas treat different errors.
Lemma~\ref{lem:gradient-commutator} controls the
commutator with the horizontal gradient, and
Lemma~\ref{lem:dyadic-unsmoothing} compares smooth averaging
with conditional expectation.
Together they give the following estimate.

\begin{prop}\label{prop:unsmoothed-principal-part}
For every $u\in\V$,
\begin{equation}\label{w:unsmoothed-principal-part}
 \int_0^\infty
 \|Q_\lambda u
   -\gamma_\lambda\mathsf E_\lambda\nabla_Xu\|_2^2
 \frac{\d\lambda}{\lambda}
 \lesssim E(u)^2.
\end{equation}
\end{prop}
\begin{proof}
The one-scale bound \eqref{wt:gamma-average-bounded}
gives $\|\gamma_\lambda\mathsf E_\lambda F\|_2\lesssim\|F\|_2$.
Since $\mathsf E_\lambda^2=\mathsf E_\lambda$, we have the identity
\[
 \begin{split}
 \gamma_\lambda\mathsf E_\lambda(\nabla_XP_\lambda-\nabla_X)u
 &=\gamma_\lambda\mathsf E_\lambda[\nabla_X,P_\lambda]u+\gamma_\lambda\mathsf E_\lambda
             (P_\lambda-\mathsf E_\lambda)\nabla_Xu.
 \end{split}
\]
Thus, applying the one-scale bound to each term and integrating
over $\lambda$ yields
\[
 \begin{split}
 \int_0^\infty
 \|\gamma_\lambda\mathsf E_\lambda
          (\nabla_XP_\lambda-\nabla_X)u\|_2^2
                                      \frac{\d\lambda}{\lambda}
 &\lesssim
 \int_0^\infty\|[\nabla_X,P_\lambda]u\|_2^2
                                      \frac{\d\lambda}{\lambda}+\int_0^\infty
       \|(P_\lambda-\mathsf E_\lambda)\nabla_Xu\|_2^2
                                      \frac{\d\lambda}{\lambda}\\
 &\lesssim E(u)^2+\|\nabla_Xu\|_2^2
 \lesssim E(u)^2.
 \end{split}
\]
Here the first integral is controlled by
Lemma~\ref{lem:gradient-commutator}, while the second
is controlled by Lemma~\ref{lem:dyadic-unsmoothing}, applied
componentwise to $\nabla_Xu$.
Finally, decompose the unsmoothed error as
\[
 \begin{split}
 Q_\lambda u-\gamma_\lambda\mathsf E_\lambda\nabla_Xu
 &=Q_\lambda u
       -\gamma_\lambda\mathsf E_\lambda\nabla_XP_\lambda u+\gamma_\lambda\mathsf E_\lambda
                      (\nabla_XP_\lambda-\nabla_X)u.
 \end{split}
\]
The smoothed principal-part estimate
\eqref{aux:smoothed-PP} of Lemma~\ref{lem:smoothed-principal-part} controls the square function of the
first term, and the estimate just proved controls that of
the second. Hence \eqref{w:unsmoothed-principal-part} follows.\end{proof}

\section{Localized testing and the Carleson estimate}
\label{sec:cm}
In Theorem~\ref{wt:carleson} below we establish the key Carleson estimate and use it to complete the proof of
Theorem~\ref{thm:main}. We first state the Carleson estimate and
deduce its consequences. The proof of Theorem~\ref{wt:carleson} is then given in
Subsections~\ref{sub1}--\ref{sub4}.

For a dyadic cube $Q$ put
\[
 W_Q=Q\times(\kappa\ell(Q),\ell(Q)],\qquad
 \mathsf T_Q=Q\times(0,\ell(Q)].
\]
The Whitney boxes partition the upper half-space modulo null sets,
and, modulo null sets, $\mathsf T_Q$ is the union of $W_R$ over descendants $R\subseteq Q$.
The row vector $\gamma_\lambda$ acts complex-linearly on columns
through the bilinear row--column product.
We will use the unsmoothed principal-part estimate
\eqref{w:unsmoothed-principal-part} and the resolvent bounds
\eqref{wt:resolvent-input}, both proved independently of the
Carleson estimate. The testing organization follows
\cite[Sections~8.1--8.2]{AtaeiEgertNystrom}. Here the transport
derivative of the cutoff is retained explicitly.

\begin{thm}\label{wt:carleson}
The principal coefficients satisfy
\begin{equation}\label{wt:carleson-conclusion}
 \sup_Q\frac1{|Q|}
       \int_0^{\ell(Q)}\int_Q|\gamma_\lambda(z)|^2\d z
                                 \frac{\d\lambda}{\lambda}
       \lesssim 1.
\end{equation}
The corresponding coefficients for $H_A^*$ satisfy the same estimate.
\end{thm}

Let $\nu$ be a positive Borel measure on $G\times(0,\infty)$ with
finite dyadic Carleson norm
\[
 A_\nu:=\sup_Q\frac{\nu(\mathsf T_Q)}{|Q|}<\infty.
\]
Then, for every $F\in\mathrm L^2(G)$,
\begin{equation}\label{wt:carleson-embedding}
 \iint_{G\times(0,\infty)}
       |\mathsf E_\lambda F(z)|^2\d\nu(z,\lambda)
       \lesssim A_\nu\|F\|_2^2.
\end{equation}
To prove \eqref{wt:carleson-embedding}, write
its left-hand side as $\sum_Q\nu(W_Q)|F_Q|^2$ and let $M_{\!d}F$ denote the dyadic maximal function
of $|F|$. The layer-cake formula and the Carleson condition give
\[
 \begin{split}
 \sum_Q\nu(W_Q)|F_Q|^2
 =\int_0^\infty
       2s\sum_{\{Q:\,|F_Q|>s\}}\nu(W_Q)\,\d s&\lesssim A_\nu
       \int_0^\infty 2s\,|\{M_{\!d}F>s\}|\,\d s\\
       &=A_\nu\|M_{\!d}F\|_2^2
 \lesssim A_\nu\|F\|_2^2.
 \end{split}
\]
For the first inequality, use $|F_Q|\le\fint_Q|F|\d z$
and group the cubes in the inner sum under the maximal
dyadic cubes on which $\fint_Q|F|\d z>s$.
The Carleson condition bounds the contribution beneath
each maximal cube by a structural constant times
$A_\nu|Q|$. These maximal cubes are disjoint and lie
in $\{M_{\!d}F>s\}$.
The final inequality is the dyadic maximal inequality,
which proves \eqref{wt:carleson-embedding}. The same estimate and deduction hold for finite-dimensional vector-valued
$F$.

Using Theorem~\ref{wt:carleson} we can now prove the central square function estimate.
\begin{prop}\label{prop:rough-quadratic}
For $B=H_A$ and $B=H_A^*$, and every $u\in\V$,
\begin{equation}\label{eq:rough-quadratic}
 \int_0^\infty
 \|\lambda B(I+\lambda^2B)^{-1}u\|_2^2\frac{\d\lambda}{\lambda}
 \lesssim E(u)^2.
\end{equation}
\end{prop}

\begin{proof}
For $H_A$, Theorem~\ref{wt:carleson} shows that the measure
\[
 \d\nu(z,\lambda)
 :=|\gamma_\lambda(z)|^2\d z\frac{\d\lambda}{\lambda}
\]
is a Carleson measure with constant $A_\nu\lesssim1$. Consequently, applying \eqref{wt:carleson-embedding}
to $F=\nabla_Xu$ gives
\[
 \begin{split}
 \int_0^\infty
   \|\gamma_\lambda\mathsf E_\lambda\nabla_Xu\|_2^2
                                      \frac{\d\lambda}{\lambda}
 &\le \iint_{G\times(0,\infty)}
   |\mathsf E_\lambda\nabla_Xu(z)|^2\d\nu(z,\lambda)\lesssim\|\nabla_Xu\|_2^2\le E(u)^2.
 \end{split}
\]
Combining this bound with \eqref{w:unsmoothed-principal-part}
proves \eqref{eq:rough-quadratic} for $H_A$.
For the adjoint, repeat the argument with $A^*$ and $-T$.
The form space, uncertainty inequality, and fixed Rockland smoothing
are unchanged. In particular, $(-T)^2=T^2$, the cutoff estimates
involve absolute values, and the bracket-generating geometry is the
same. All constants have the same structural dependence.
\end{proof}

\begin{proof}[Proof of Theorem~\ref{thm:main}]
Proposition~\ref{prop:kinetic-realization} constructs the maximal
accretive realizations and identifies the adjoint.
Proposition~\ref{prop:rough-quadratic} supplies the two estimates
required in Proposition~\ref{prop:quadratic-reduction}, with constants
depending only on $m,a_0,a_1$. That proposition gives the domain
identities and the two norm equivalences.
\end{proof}

\subsection{Resolvent tests and the transport cutoff error}\label{sub1}

Starting here, the proof of Theorem~\ref{wt:carleson} runs through Subsections~\ref{sub1}--\ref{sub4}. In what follows, $\varepsilon$ and $\delta$ denote new small
parameters chosen for testing and boundary localization.
Their values are chosen independently of the parameters
fixed in the symmetrizer construction in
Section~\ref{sec:variational-realization}.

For a cube $Q$ and a unit column $\zeta\in\C^m$, choose
a smooth outer cutoff $\chi_Q$ in the left-translated
coordinates $(\xi,\eta,\tau)=z_Q^{-1}\circ z$,
centred at $z_Q$ and scaled according to the kinetic dilations. We choose $\chi_Q$ so that it equals one on a neighbourhood
of the ball $B(z_Q,C\ell(Q))$ containing $Q$, is supported
in a fixed concentric enlargement of this ball, and satisfies
\[
 |\nabla_X\chi_Q|\lesssim\ell(Q)^{-1},
 \qquad
 |T\chi_Q|\lesssim\ell(Q)^{-2}.
\] Let $X_Q$ be the horizontal
coordinate of $z_Q$ and put
\[
 L_Q^\zeta(z):=\chi_Q(z)\bigl((X-X_Q)\cdot\zeta\bigr),
 \qquad \lambda_Q=\varepsilon\ell(Q),\qquad
 b_Q^\zeta=R_{\lambda_Q}L_Q^\zeta.
\]
Since $\chi_Q$ is supported in a fixed enlargement of a ball
at scale $\ell(Q)$, volume growth and the dyadic cube estimates
give $|\operatorname{supp}\chi_Q|\lesssim|Q|$.
Moreover, $|X-X_Q|\lesssim\ell(Q)$ on this support.
The cutoff bounds and the identities $TX_j=0$ therefore yield
\begin{align}\label{wt:test-cutoff-norms}
 \|L_Q^\zeta\|_2&\lesssim\ell(Q)|Q|^{1/2},\qquad
 \|\nabla_XL_Q^\zeta\|_2\lesssim|Q|^{1/2},\notag\\
 TL_Q^\zeta&=(T\chi_Q)(X-X_Q)\cdot\zeta,\qquad
 \|TL_Q^\zeta\|_2\lesssim\ell(Q)^{-1}|Q|^{1/2}.
\end{align}
Applying Cauchy--Schwarz in the spectral representation of
$|T|$ gives
\[\||T|^{1/2}L_Q^\zeta\|_2^2
 \le\|TL_Q^\zeta\|_2\|L_Q^\zeta\|_2,\qquad
E(L_Q^\zeta)\lesssim|Q|^{1/2}.\] We have the exact weak identity
\begin{align}\label{wt:exactid}
 b_Q^\zeta-L_Q^\zeta
 &=\lambda_Q^2R_{\lambda_Q}\operatorname{div}_X
                   (A\nabla_XL_Q^\zeta)
       -\lambda_Q^2R_{\lambda_Q}TL_Q^\zeta.
\end{align}
The second term is the transport cutoff error. Since
$\lambda_Q=\varepsilon\ell(Q)$, the bounds
\eqref{wt:resolvent-input} and \eqref{wt:test-cutoff-norms} give
\[
 \begin{split}
 \|\lambda_Q^2R_{\lambda_Q}TL_Q^\zeta\|_2
 &\lesssim\lambda_Q^2\ell(Q)^{-1}|Q|^{1/2}
 =\varepsilon\lambda_Q|Q|^{1/2},\\
 E(\lambda_Q^2R_{\lambda_Q}TL_Q^\zeta)
 &\lesssim\lambda_Q\ell(Q)^{-1}|Q|^{1/2}
 =\varepsilon|Q|^{1/2}.
 \end{split}
\]
Together with the divergence-forcing bound in
\eqref{wt:resolvent-input}, the weak identity in \eqref{wt:exactid} therefore implies,
uniformly for $0<\varepsilon\le1$,
\begin{equation}\label{wt:test-difference}
 \|b_Q^\zeta-L_Q^\zeta\|_2\lesssim\lambda_Q|Q|^{1/2},
 \qquad E(b_Q^\zeta-L_Q^\zeta)+E(b_Q^\zeta)\lesssim|Q|^{1/2}.
\end{equation}

Since $b_Q^\zeta\in\mathrm D(H_A)$ and
$H_Ab_Q^\zeta=\lambda_Q^{-2}(L_Q^\zeta-b_Q^\zeta)$,
contractivity of $R_\lambda$ and resolvent commutation on
$\mathrm D(H_A)$ give
\[
 \|\lambda H_AR_\lambda b_Q^\zeta\|_2
 \le\frac{\lambda}{\lambda_Q^2}
       \|L_Q^\zeta-b_Q^\zeta\|_2.
\]
Squaring and integrating, we obtain
\[
 \int_0^{\ell(Q)}\|Q_\lambda b_Q^\zeta\|_2^2
                                      \frac{\d\lambda}{\lambda}=\int_0^{\ell(Q)}\|\lambda H_AR_\lambda b_Q^\zeta\|_2^2
                                      \frac{\d\lambda}{\lambda}
 \le\frac{\ell(Q)^2}{2\lambda_Q^4}
        \|L_Q^\zeta-b_Q^\zeta\|_2^2
 \lesssim\varepsilon^{-2}|Q|.
\]
The last bound uses \eqref{wt:test-difference} and
$\lambda_Q=\varepsilon\ell(Q)$. Combining this estimate with
\eqref{w:unsmoothed-principal-part} and the energy bound
in \eqref{wt:test-difference} yields
\begin{equation}\label{wt:testing}
 J_Q^\zeta:=\int_0^{\ell(Q)}\int_Q
       |\gamma_\lambda\mathsf E_\lambda\nabla_Xb_Q^\zeta|^2\d z
                                       \frac{\d\lambda}{\lambda}
       \lesssim(1+\varepsilon^{-2})|Q|.
\end{equation}
Neither the Carleson estimate nor its embedding theorem has been
used in proving this testing bound.

\subsection{Average accretivity}\label{sub2}

Let $\eta>0$ be the small-boundary exponent in
\eqref{aux:small-boundary}, and let $\delta_0>0$ be a sufficiently
small fixed geometric constant. For $0<\delta<\delta_0$ we construct
$\phi=\phi_{Q,\delta}\in\mathrm C_0^\infty(Q)$ with
\[
 0\le\phi\le1,\qquad
 |\{z\in Q:\phi(z)\ne1\}|\lesssim\delta^\eta|Q|,
 \qquad |\nabla_X\phi|\lesssim(\delta\ell(Q))^{-1}.
\]
Indeed, fix a nonnegative $\varrho\in\mathrm C_0^\infty(G)$ of integral
one, supported in a sufficiently small ball about the identity.
Choose a sufficiently small fixed geometric constant $c>0$,
and for $\lambda_\delta=c\delta\ell(Q)$ set
$\varrho_{\lambda_\delta}=\lambda_\delta^{-\mathfrak q}
\varrho\circ\delta_{1/\lambda_\delta}$.
For a sufficiently large fixed geometric constant $C_*$,
with $c$ chosen small enough in terms of $C_*$, put
\[
 Q_{\lambda_\delta}=\{z\in Q:d(z,Q^c)>C_*\lambda_\delta\},\qquad
 \phi_{Q,\delta}=1_{Q_{\lambda_\delta}}*\varrho_{\lambda_\delta}.
\]
The quasi-triangle inequality ensures that this cutoff has compact
support in $Q$ and equals one outside a boundary strip whose width
is a fixed multiple of $\delta\ell(Q)$. A left-invariant horizontal
derivative acts on the second convolution factor, whose derivative
has $\mathrm L^1$ norm $\lesssim\lambda_\delta^{-1}$.
The strip estimate follows from
\eqref{aux:small-boundary}.

Set $h=L_Q^\zeta-b_Q^\zeta$ and
$E_\delta=\{\phi\ne1\}\cap Q$. Since
$\phi\in\mathrm C_0^\infty(Q)$, horizontal integration by parts gives
\[
 \int_Q\nabla_Xh\d z
 =\int_Q(1-\phi)\nabla_Xh\d z
       -\int_Qh\nabla_X\phi\d z.
\]
The first integrand on the right vanishes outside $E_\delta$. The cutoff bounds
and \eqref{wt:test-difference} therefore imply
\begin{align}\label{osc}
 \left|\fint_Q\nabla_Xh\d z\right|
 &\le\frac1{|Q|}
    \left(|E_\delta|^{1/2}\|\nabla_Xh\|_2
        +\|h\|_2\|\nabla_X\phi\|_2\right)\lesssim\delta^{\eta/2}
                   +\frac{\lambda_Q}{\delta\ell(Q)}
 =\delta^{\eta/2}+\varepsilon/\delta.
\end{align}
Let $C$ be the structural constant in the estimate in \eqref{osc}. Choose $\delta>0$ first so that $C\delta^{\eta/2}\le1/8$,
and then choose $0<\varepsilon\le1$ so that
$C\varepsilon/\delta\le1/8$. These choices depend only on
the structural parameters and make the preceding bound
at most $1/4$. Fix this value of $\varepsilon$. Set $g:=\nabla_Xb_Q^\zeta$ and
$g_Q:=\fint_Qg\d z$. Since $\nabla_XL_Q^\zeta=\zeta$
on $Q$ and $h=L_Q^\zeta-b_Q^\zeta$, we have $|\zeta-g_Q|\leq 1/4$ by \eqref{osc}. Using $|\zeta|=1$, Cauchy--Schwarz and
\eqref{wt:test-difference}, we therefore obtain
\begin{equation}\label{wt:accretivity}
\operatorname{Re}\langle g_Q,\zeta\rangle
 \ge 1-|g_Q-\zeta|\,|\zeta|
 \ge\frac34,
 \qquad \int_Q|g|^2\d z\le M|Q|,
\end{equation}
where $M$ depends only on $m,a_0,a_1$ and the fixed
geometric choices.

\subsection{Stopping cubes and packing}\label{sub3}

Cauchy--Schwarz and \eqref{wt:accretivity} imply, for $E\subset Q$,
\begin{equation}\label{wt:g-small-sets}
 \frac1{|Q|}\int_E|g|\d z
 \le C_0\left(\frac{|E|}{|Q|}\right)^{1/2},
 \qquad C_0=\max\{1,\sqrt M\}.
\end{equation}
Choose $K_0>2C_0$ so large that
$C_0(C_0/K_0)^{1/2}\le1/8$. We now select maximal strict dyadic
subcubes $Q_j\subset Q$ satisfying at least one of
\begin{equation}\label{wt:stopping-rules}
 \operatorname{Re}\langle g_{Q_j},\zeta\rangle\le\frac12,
 \qquad \fint_{Q_j}|g|\d z\ge K_0.
\end{equation}
The top cube satisfies none of these conditions.
We divide the selected cubes into two disjoint classes.
The first class consists of cubes satisfying the first
condition in \eqref{wt:stopping-rules} but not the second. The second class consists
of cubes satisfying the second condition in \eqref{wt:stopping-rules}, whether or not
they also satisfy the first.  Let $E_1$ and $E_2$ be the unions of the cubes in the first
and second classes, respectively, and set
$G_Q:=Q\setminus(E_1\cup E_2)$.
Summing the size condition over the disjoint cubes in the
second class and applying \eqref{wt:g-small-sets} with $E=Q$
gives
\begin{equation}\label{wt:stopping-rulesCons}
 |E_2|
 \le K_0^{-1}\int_{E_2}|g|\d z
 \le K_0^{-1}\int_Q|g|\d z
 \le\frac{C_0}{K_0}|Q|.
\end{equation}
Applying \eqref{wt:g-small-sets} again, now with $E=E_2$ and using \eqref{wt:stopping-rulesCons},
we obtain
 \begin{align}\label{SatEve}
 \int_{E_2}|g|\d z\le C_0\left(\frac{|E_2|}{|Q|}\right)^{1/2}|Q|\le C_0\left(\frac{C_0}{K_0}\right)^{1/2}|Q|
 \le\frac18|Q|.
 \end{align}
Using the inequality for $\operatorname{Re}\langle g_Q,\zeta\rangle$ in \eqref{wt:accretivity}, decomposing
$Q=E_1\cup E_2\cup G_Q$, summing the low-accretivity condition over the disjoint cubes
forming $E_1$, and recalling that $|\zeta|=1$, we obtain
\[
 \begin{aligned}
 \frac34|Q|
 &\le\frac12|E_1|
      +\int_{E_2}|g|\d z+\int_{G_Q}|g|\d z\le\frac58|Q|
      +C_0|Q|^{1/2}|G_Q|^{1/2}.
 \end{aligned}
\]
The last inequality uses $|E_1|\le|Q|$, the  bound in \eqref{SatEve}, and
\eqref{wt:g-small-sets} applied to $G_Q$. Consequently, with $\theta=(8C_0)^{-2}>0$,
\begin{equation}\label{eq:stopping-packing}
 |G_Q|\ge\theta|Q|,\qquad
 \sum_j|Q_j|=|Q|-|G_Q|\le(1-\theta)|Q|.
\end{equation}
Every descendant $R$ not contained in a stopping cube satisfies
\begin{equation}\label{wt:good-averages}
 \operatorname{Re}\langle g_R,\zeta\rangle>\frac12,\quad\mbox{and}
 \quad |g_R|\le\fint_R|g|\d z<K_0.
\end{equation}

\subsection{Cones and finite absorption}\label{sub4}

Fix an aperture $0<\omega\le(4K_0)^{-1}$ depending only
on the structural data. By compactness of the unit sphere,
we can choose finitely many unit columns $\zeta$ so that
the cones
\[
 \mathcal C_\zeta=
   \bigl\{a\in\C^{1\times m}\setminus\{0\}:
       |a/|a|-\overline{\zeta}^{\,T}|\le\omega\bigr\}
\]
cover $\C^{1\times m}\setminus\{0\}$.
The number of cones depends only on $m$ and $\omega$.
For a column $v$, one has
$\overline{\zeta}^{\,T}v=\langle v,\zeta\rangle$, so this
choice of cone is compatible with the accretivity condition
\eqref{wt:good-averages}. Fix one such $\zeta$ and a cube $Q$, and use the test field
$g=\nabla_Xb_Q^\zeta$ and its stopping cubes $Q_j$.
Define
\[
 \mathsf S_Q=\mathsf T_Q\setminus\bigcup_j\mathsf T_{Q_j}.
\]
If $W_R\subset\mathsf S_Q$, then $R$ is not contained in
any stopping cube. Consequently, for
$(z,\lambda)\in W_R$, the vector
$v=\mathsf E_\lambda g(z)=g_R
 =\fint_R\nabla_X b_Q^\zeta(w)\d w$  satisfies
\eqref{wt:good-averages}.
If $\gamma_\lambda(z)\in\mathcal C_\zeta$, the cone
condition and the triangle inequality give
\[
 \begin{aligned}
 |\gamma_\lambda v|
 &\ge|\gamma_\lambda|
       \bigl(|\langle v,\zeta\rangle|-\omega|v|\bigr)\ge|\gamma_\lambda|
       \bigl(\operatorname{Re}\langle v,\zeta\rangle
                     -\omega|v|\bigr)\ge\bigl(\tfrac12-\omega K_0\bigr)|\gamma_\lambda|
 \ge\tfrac14|\gamma_\lambda|.
 \end{aligned}
\]
With
\[
 \d\nu_\zeta
 :=1_{\mathcal C_\zeta}(\gamma_\lambda)
       |\gamma_\lambda|^2\d z\,\frac{\d\lambda}{\lambda},
\]
the testing estimate \eqref{wt:testing} therefore yields
\begin{equation}\label{Tin}
 \nu_\zeta(\mathsf S_Q)
 \le16\iint_{\mathsf S_Q}
       |\gamma_\lambda\mathsf E_\lambda g|^2
       \d z\,\frac{\d\lambda}{\lambda}
 \le16J_Q^\zeta
 \le C_\varepsilon|Q|.
\end{equation}
The parameter $\varepsilon$ was fixed in terms of the structural
data in Subsection~\ref{sub2}, so $C_\varepsilon$ is structural
and independent of $Q$ and $\zeta$.

To conclude we first restrict to $\lambda_-<\lambda<\lambda_+$, with
$0<\lambda_-<\lambda_+<\infty$, and define
\[
 \mathcal C^\zeta_{\lambda_-,\lambda_+}
 =\sup_Q\frac{\nu_\zeta\bigl(\mathsf T_Q\cap
                  \{\lambda_-<\lambda<\lambda_+\}\bigr)}{|Q|},
\]
where the supremum is over all dyadic cubes.
This quantity is finite before absorption.
Indeed, for each $0<\lambda\le\ell(Q)$, the averaging-scale
cubes $R$ contained in $Q$ are descendants of $Q$ and
partition it up to a null set.
Since $\lambda\sim\ell(R)$, the one-scale estimate
\eqref{wt:gamma-one-scale} gives
\[
 \int_Q|\gamma_\lambda|^2\d z
 =\sum_R\int_R|\gamma_\lambda|^2\d z
 \lesssim\sum_R|R|
 =|Q|.
\]
Integrating over the truncated height interval, and using
$\d\nu_\zeta\le|\gamma_\lambda|^2\d z\,\d\lambda/\lambda$,
yields $\mathcal C^\zeta_{\lambda_-,\lambda_+}
 \lesssim\log(\lambda_+/\lambda_-)<\infty$. For each cube $Q$, the tent decomposition and \eqref{Tin} now give
\[
 \begin{aligned}
 \nu_\zeta\bigl(\mathsf T_Q\cap
                  \{\lambda_-<\lambda<\lambda_+\}\bigr)&\le\nu_\zeta(\mathsf S_Q)
       +\sum_j\nu_\zeta\bigl(\mathsf T_{Q_j}\cap
                  \{\lambda_-<\lambda<\lambda_+\}\bigr)\\
 &\le C_\varepsilon|Q|
       +\mathcal C^\zeta_{\lambda_-,\lambda_+}\sum_j|Q_j|\le
 \bigl(C_\varepsilon
       +(1-\theta)\mathcal C^\zeta_{\lambda_-,\lambda_+}\bigr)|Q|,
 \end{aligned}
\]
where the last step uses \eqref{eq:stopping-packing}.
Dividing by $|Q|$ and taking the supremum over $Q$, we obtain
\[
 \mathcal C^\zeta_{\lambda_-,\lambda_+}
 \le C_\varepsilon
       +(1-\theta)\mathcal C^\zeta_{\lambda_-,\lambda_+}.
\]
Since this quantity is finite and $\theta>0$, we can absorb
the last term and conclude that $\mathcal C^\zeta_{\lambda_-,\lambda_+}
 \le {C_\varepsilon}/{\theta}$, uniformly in both height cutoffs. 
 
 Finally, for each fixed $Q$, let $\lambda_-\downarrow0$ and
$\lambda_+\uparrow\infty$. Monotone convergence gives
$\nu_\zeta(\mathsf T_Q)\le(C_\varepsilon/\theta)|Q|$.
Since the finitely many cones cover every nonzero row,
\[
 \int_{\mathsf T_Q}|\gamma_\lambda(z)|^2
       \d z\,\frac{\d\lambda}{\lambda}
 \le\sum_\zeta\nu_\zeta(\mathsf T_Q)
 \lesssim|Q|.
\]
This proves Theorem~\ref{wt:carleson} for $H_A$.
The proof for $H_A^*$ uses $-T,A^*$ in place of $T,A$.
Both resolvent estimates and the unsmoothed principal-part
estimate hold for this pair. The cutoff bounds use absolute
values of transport derivatives, so the same testing and
packing proof applies.

\section{Shifted square roots and Riesz transforms}
\label{sec:consequences}

The homogeneous estimates imply a uniform inhomogeneous statement.
This also gives bounded Riesz transforms without choosing a homogeneous
completion of the energy space.

\begin{cor}\label{cor:shifted-roots}
Let $B$ denote either $H_A$ or $H_A^*$. For every $\mu>0$,
\begin{equation}\label{eq:shifted-roots}
 \mathrm D((\mu+B)^{1/2})=\V,
 \qquad
 \| (\mu+B)^{1/2}u\|_2
 \sim\sqrt\mu\|u\|_2+\|\nabla_Xu\|_2+\||T|^{1/2}u\|_2,
 \quad u\in\V.
\end{equation}
The constants depend only on $m,a_0,a_1$, independently of $\mu$.
Consequently,
\[
 \sqrt\mu(\mu+B)^{-1/2},\quad
 \nabla_X(\mu+B)^{-1/2},\quad
 |T|^{1/2}(\mu+B)^{-1/2},\quad
 D_Y^{1/3}(\mu+B)^{-1/2}
\]
are bounded on $\mathrm L^2(G)$, uniformly in $\mu>0$. The horizontal transform
takes values in $\mathrm L^2(G;\C^m)$.
\end{cor}

\begin{proof}
For a maximal accretive operator $B$, subtracting the Balakrishnan
representations of $(\mu+B)^{1/2}$ and $B^{1/2}$
\cite[Chapter~3]{HaaseFunctionalCalculus} and using the resolvent
identity gives, on $\mathrm D(B)$,
\[
 (\mu+B)^{1/2}-B^{1/2}
 =\frac\mu\pi\int_0^\infty
 t^{1/2}(t+B)^{-1}(t+\mu+B)^{-1}\d t.
\]
The integral converges in operator norm, and its norm is at most
\[
 \frac\mu\pi\int_0^\infty
 \frac{\d t}{\sqrt t(t+\mu)}=\sqrt\mu.
\]
Since $\mathrm D(B)=\mathrm D(\mu+B)$ is a core for both square roots, their domains
agree and the bounded-difference identity holds on that common domain.
The corresponding inverse-power formula gives
\[
 \|(\mu+B)^{-1/2}\|_{2\to2}
 \le\frac1\pi\int_0^\infty
 \frac{t^{-1/2}}{t+\mu}\d t=\mu^{-1/2}.
\]
It follows that
\[
 \sqrt\mu\|u\|_2\le\|(\mu+B)^{1/2}u\|_2,
 \qquad
 \|B^{1/2}u\|_2\le2\|(\mu+B)^{1/2}u\|_2.
\]
These estimates, the bounded difference and Theorem~\ref{thm:main}
prove both directions of \eqref{eq:shifted-roots}. Apply the result to
$u=(\mu+B)^{-1/2}f$ and use \eqref{eq:y-consequence} to obtain the
Riesz-transform bounds.
\end{proof}

\begin{rem}
The homogeneous transforms
$\nabla_X B^{-1/2}$, $|T|^{1/2}B^{-1/2}$ and $D_Y^{1/3}B^{-1/2}$,
initially defined on $\Ran(B^{1/2})$, extend uniquely to bounded
operators on $\mathrm L^2(G)$. Indeed, their bounds follow directly from
Theorem~\ref{thm:main} and \eqref{eq:y-consequence}. Their initial
domain is dense because $B^*$, and hence $(B^*)^{1/2}$, is injective.
Injectivity was established in Remark~\ref{rem:injective}. Thus
$\overline{\Ran(B^{1/2})}=
(\ker((B^*)^{1/2}))^\perp=\mathrm L^2(G)$.
The unbounded inverse power $B^{-1/2}$ is used only on this initial
domain.
\end{rem}

\appendix

\section{Fractional powers and the common core}
\label{app:common-core}
This appendix is devoted to the proof of  Lemma~\ref{wr:common-core}.

\begin{proof}[Proof of Lemma~\ref{wr:common-core}]
 We first prove density in $\V$. By
Lemma~\ref{lem:form-density}, it suffices to approximate
each $v\in\mathscr S(G)$ by elements of $\mathscr S_0(G)$
in the norm of $\V$.
Write $V:=\widehat v$ for the partial Fourier transform
in $(Y,t)$, and choose smooth cutoffs
$\chi_n(\eta,\tau)$ with $0\le\chi_n\le1$, supported in \(\{1/n<|(\eta,\tau)|<n\}\), and converging pointwise to one away from the origin.
Define $v_n$ by $\widehat v_n:=\chi_n V$. Since $V$ is Schwartz, each $v_n$ belongs to $\mathscr S(G)$.
Moreover, $\widehat v_n$ vanishes in a neighbourhood of
$(\eta,\tau)=0$. The Fourier identities for moments therefore
give
\[
 \int_{\R^{m+1}}Y^\beta t^\ell v_n(X,Y,t)\,\d Y\d t=0
 \qquad
 \text{for every }X,\quad
 \beta\in\mathbb N_0^m,\quad \ell\in\mathbb N_0.
\]
Multiplying by any monomial in $X$ and integrating shows
that $v_n\in\mathscr S_0(G)$. Multiplication by $\chi_n$ commutes with $\nabla_X$ and
with multiplication by $|q|^{1/2}$, where
$q=\tau-X\cdot\eta$. Dominated convergence yields
\[
 \begin{split}
 \|(\chi_n-1)V\|_2
 &+\|(\chi_n-1)\nabla_XV\|_2+\|(\chi_n-1)|q|^{1/2}V\|_2
 \longrightarrow0.
 \end{split}
\]
By partial Plancherel, this is convergence of $v_n$ to $v$
in $\V$. A diagonal choice, starting from the Schwartz
approximants supplied by Lemma~\ref{lem:form-density},
proves density of $\mathscr S_0(G)$ in $\V$.
Since $\V$ is dense in $\mathrm L^2(G)$, this also proves
density of $\mathscr S_0(G)$ in $\mathrm L^2(G)$.

We next establish invariance under complex powers.
Recall that $\mathscr S_0(G)$ is a closed subspace of
$\mathscr S(G)$, since every polynomial moment is a
continuous functional on $\mathscr S(G)$.
Convolution with $h_s$ preserves $\mathscr S_0(G)$.
Indeed, for every polynomial $Q$ and
$f\in\mathscr S_0(G)$, the polynomial group law and
invariance of Haar measure give
\[
 \int_G Q(z)(f*h_s)(z)\,\d z
 =\int_G h_s(w)
       \left(\int_G Q(y\circ w)f(y)\,\d y\right)\d w=0.
\]
The inner integral vanishes because
$y\mapsto Q(y\circ w)$ is a polynomial for each fixed $w$.  To construct negative powers on $\mathscr S_0(G)$, we
need decay of the heat semigroup in the Schwartz topology
as $s\to\infty$. Write $s=\lambda^4$ and use the
heat-kernel scaling to obtain
\[
 e^{-s\mathcal R}f(z)
 =\lambda^{-\mathfrak q}
   \int_G f(y)
   h_1\bigl((\delta_{\lambda^{-1}}y)^{-1}
                 \circ\delta_{\lambda^{-1}}z\bigr)\,\d y.
\]
We work in Euclidean coordinates and use the Schwartz
seminorms
\[
 p_{L,B}(g):=\max_{|\beta|\le B}
       \sup_z(1+|z|)^L|\partial_z^\beta g(z)|,
 \qquad L,B\in\mathbb N_0.
\]
A smooth function $g$ belongs to $\mathscr S(G)$ if and only if
\[
 p_{L,B}(g)<\infty
 \qquad\text{for every }L,B\in\mathbb N_0.
\]
Set $H(a,w):=h_1(a^{-1}\circ w)$.
Taylor's formula in $a$ gives, for each integer $k\ge1$,
\[
 H(a,w)=\sum_{|\gamma|<k}
       \frac{a^\gamma}{\gamma!}\partial_a^\gamma H(0,w)
       +\mathcal E_k(a,w).
\]
Since the group law and inversion are polynomial, the
Schwartz bounds for $h_1$ and the integral form of the
Taylor remainder imply
\[
 p_{L,B}(\mathcal E_k(a,\cdot))
 \lesssim_{L,B,k}|a|^k(1+|a|)^M
\]
for some $M=M(L,B,k)$.
Indeed, derivatives of $H$ are finite sums of derivatives
of $h_1(a^{-1}\circ w)$ with polynomial coefficients.
Writing $w=a\circ v$ bounds their weighted suprema by
a polynomial in $a$ and finitely many Schwartz seminorms
of $h_1$.

For $\lambda\ge1$, the kinetic weights $1,3,2$ give
\[
 |\delta_{\lambda^{-1}}y|\le\lambda^{-1}|y|,
 \qquad
 1+|\delta_\lambda w|\le\lambda^3(1+|w|).
\]
All Taylor terms integrate to zero against $f$, because
they are polynomials in $y$. Differentiation with respect
to $z$ introduces only additional nonpositive powers of
$\lambda$. The remainder therefore yields
\[
 p_{L,B}(e^{-\lambda^4\mathcal R}f)
 \lesssim_{L,B,k}
 \lambda^{3L-\mathfrak q-k}
       \int_G(1+|y|)^{k+M}|f(y)|\,\d y.
\]
The weighted integral is a continuous Schwartz seminorm. Since the seminorms $p_{L,B}$ generate the Schwartz topology
and increase with $L$ and $B$, every continuous seminorm $p$
on $\mathscr S(G)$ satisfies
\[
 p(g)\le C_p\,p_{L,B}(g),
 \qquad g\in\mathscr S(G),
\]
for suitable $L,B\in\mathbb N_0$. Hence, using this and choosing $k$ sufficiently large we can conclude that for every Schwartz seminorm $p$ and every $N>0$,
there is a continuous Schwartz seminorm $q_{p,N}$ such that
\begin{equation}\label{Mellin1}
 p(e^{-s\mathcal R}f)
 \lesssim_{p,N}q_{p,N}(f)s^{-N},
 \qquad s\ge1,\quad f\in\mathscr S_0(G).
\end{equation}
We next consider the behaviour as $s\downarrow0$. Applying \cite[Lemma~2.6]{FischerRuzhansky}
to the Schwartz kernel $h_1$, with integral one and
dilation parameter $s^{1/4}$, gives 
\[e^{-s\mathcal R}f\longrightarrow f\quad\mbox{in}\quad\mathscr S(G)\quad\mbox{as}\quad s\downarrow0.\] To obtain a uniform seminorm bound, write
\[
 (f*h_s)(z)=\int_G h_1(y)
       f\bigl(z\circ(\delta_{s^{1/4}}y)^{-1}\bigr)\d y.
\]
Differentiation and the polynomial group law bound each
Schwartz seminorm of this expression by a fixed Schwartz
seminorm of $f$ times a polynomial moment of $|h_1|$,
uniformly for $0<s\le1$.
Thus, for every $p$, there is a continuous Schwartz
seminorm $q_p$ such that
\begin{equation}\label{Mellin2}
 \sup_{0<s\le1}p(e^{-s\mathcal R}f)
 \lesssim_p q_p(f),
 \qquad f\in\mathscr S(G).
\end{equation}

The estimates \eqref{Mellin1} and \eqref{Mellin2} allow us
to identify negative powers through a Mellin integral.
For $f\in\mathscr S_0(G)$ and $\operatorname{Re}a>0$, set
\[
 g_a:=\frac1{\Gamma(a)}
   \int_0^\infty s^{a-1}e^{-s\mathcal R}f\,\d s.
\]
Near zero, the integral converges in every Schwartz
seminorm by \eqref{Mellin2} and
$\operatorname{Re}a>0$.
At infinity, convergence follows from \eqref{Mellin1}
with $N>\operatorname{Re}a$.
Since $\mathscr S_0(G)$ is closed in $\mathscr S(G)$
and, as shown above, $e^{-s\mathcal R}$ maps
$\mathscr S_0(G)$ into itself for every $s>0$, the integral
converges in $\mathscr S_0(G)$.
The same estimates show that $f\mapsto g_a$ is continuous
on that space. To identify $g_a$ with the spectral power, introduce $F_n:=\mathbf1_{[1/n,n]}(\mathcal R)$. On each of these compact spectral intervals, the spectral
theorem and the scalar Gamma-integral identity give
\[
 F_ng_a=\mathcal R^{-a}F_nf.
\]
Injectivity of $\mathcal R$ implies that $F_n\to I$
strongly on $\mathrm L^2(G)$. Hence
$F_nf\to f$ and $F_ng_a\to g_a$ in $\mathrm L^2(G)$.
Closedness of the spectral power $\mathcal R^{-a}$
therefore gives
\[
 f\in\mathrm D(\mathcal R^{-a}),
 \qquad
 \mathcal R^{-a}f
 =\frac1{\Gamma(a)}
   \int_0^\infty s^{a-1}e^{-s\mathcal R}f\,\d s.
\]
In particular, $\mathcal R^{-a}$ preserves
$\mathscr S_0(G)$ continuously whenever
$\operatorname{Re}a>0$.

Nonnegative integer powers of $\mathcal R$ also preserve
$\mathscr S_0(G)$ continuously. Indeed, they are differential
operators with polynomial coefficients, and integration by
parts expresses each polynomial moment of $\mathcal R^kf$
as a polynomial moment of $f$.
For arbitrary $\alpha\in\C$, choose an integer $k\ge0$
with $k>\operatorname{Re}\alpha$. The spectral calculus gives
\[
 \mathcal R^\alpha f
 =\mathcal R^k\mathcal R^{-(k-\alpha)}f,
 \qquad f\in\mathscr S_0(G).
\]
The right-hand side is well defined by the invariance just
proved. It follows that
$\mathscr S_0(G)\subset\mathrm D(\mathcal R^\alpha)$
and that every complex power maps $\mathscr S_0(G)$
continuously into itself.

The Mellin integral can also be differentiated locally
with respect to $a$ in the Schwartz topology as the
additional powers of $\log s$ are controlled by the same
endpoint estimates. Taking $k$ fixed in a neighbourhood
of any given $\alpha$ therefore shows that
$\alpha\mapsto\mathcal R^\alpha f$ is holomorphic with
values in $\mathscr S_0(G)$.
Agreement with the spectral calculus yields
\[
 \mathcal R^\alpha\mathcal R^\beta f
 =\mathcal R^{\alpha+\beta}f,
 \qquad
 \langle\mathcal R^\alpha f,g\rangle
 =\langle f,\mathcal R^{\overline\alpha}g\rangle,
 \qquad
 f,g\in\mathscr S_0(G),\quad \alpha,\beta\in\C.
\]
This proves the asserted invariance and continuity of all
complex powers on $\mathscr S_0(G)$, together with the power
and adjoint pairing identities. It remains to show that
$\mathscr S_0(G)$ is a core for $\mathcal R^{1/4}$.

To start the proof of the core assertion, the estimates \eqref{Mellin1} and \eqref{Mellin2} show that
\[
 (I+\mathcal R)^{-1}f
 =\int_0^\infty e^{-s}e^{-s\mathcal R}f\,\d s
\]
converges in $\mathscr S_0(G)$ for every
$f\in\mathscr S_0(G)$, with the identity following from
the spectral theorem. Hence $(I+\mathcal R)^{-1}$
preserves $\mathscr S_0(G)$. Let $u\in\mathrm D(\mathcal R)$ and set
$f:=(I+\mathcal R)u$.
By the density already established, choose
$f_n\in\mathscr S_0(G)$ with $f_n\to f$ in
$\mathrm L^2(G)$, and define
\[
 u_n:=(I+\mathcal R)^{-1}f_n.
\]
Resolvent invariance gives $u_n\in\mathscr S_0(G)$.
Boundedness of the resolvent then yields
\[
 u_n\longrightarrow u,
 \qquad
 \mathcal Ru_n=f_n-u_n
 \longrightarrow f-u=\mathcal Ru
 \quad\text{in }\mathrm L^2(G).
\]
Thus $\mathscr S_0(G)$ is a core for $\mathcal R$. To obtain the corresponding assertion for
$\mathcal R^{1/4}$, let
$u\in\mathrm D(\mathcal R^{1/4})$.
The spectral truncations
\[
 u^{(N)}:=\mathbf1_{[0,N]}(\mathcal R)u
\]
belong to $\mathrm D(\mathcal R)$ and converge to $u$
in the graph norm of $\mathcal R^{1/4}$.
Moreover, the spectral theorem and H\"older's inequality give
\[
 \|\mathcal R^{1/4}v\|_2
 \le\|\mathcal Rv\|_2^{1/4}\|v\|_2^{3/4},
 \qquad v\in\mathrm D(\mathcal R).
\]
Consequently, approximation in the graph norm of
$\mathcal R$ also gives approximation in the graph norm
of $\mathcal R^{1/4}$.
Approximating each $u^{(N)}$ by elements of
$\mathscr S_0(G)$ and taking a diagonal sequence proves
that $\mathscr S_0(G)$ is a core for $\mathcal R^{1/4}$.
\end{proof}

\section{The type-zero commutator estimate}
\label{app:type-zero-commutator}

Let $\mathfrak g$ denote the Lie algebra of $G$, identified
with its left-invariant vector fields. Choose the homogeneous
basis
\[
 Z_j=\partial_{X_j},\qquad
 Z_{m+1}=T,\qquad
 Z_{m+1+j}=\partial_{Y_j},
 \qquad j=1,\ldots,m.
\]
With respect to the kinetic dilations $\delta_\lambda$, these
fields have degrees
\[
 v_j=1,\qquad v_{m+1}=2,\qquad v_{m+1+j}=3,
 \qquad j=1,\ldots,m.
\]
We write $Z_\ell^L$ and $Z_\ell^R$ for the corresponding
left- and right-invariant vector fields. Explicitly, for
$f\in\mathrm C^\infty(G)$,
\[
 \begin{aligned}
 (Z_\ell^Lf)(z)
 &:=\left.\frac{\d}{\d s}
 f\bigl(z\circ\exp(sZ_\ell)\bigr)\right|_{s=0},\\
 (Z_\ell^Rf)(z)
 &:=\left.\frac{\d}{\d s}
 f\bigl(\exp(sZ_\ell)\circ z\bigr)\right|_{s=0}.
 \end{aligned}
\]

\begin{proof}[Proof of Lemma~\ref{lem:type-zero-commutator}]
We first show that the commutator has an integrable kernel
with cancellation, a positive fractional moment, and integrable
left and right derivatives. These properties give decay between
different scales and hence the square-function estimate. Throughout this proof, implicit constants may additionally
depend on finitely many homogeneous kernel seminorms of $K$.
Further dependence on a differential monomial, a decay exponent,
a Schwartz function, or a fractional moment exponent is indicated
by subscripts on $\lesssim$.

Write $k$ for the smooth restriction of $K$ to $G\setminus\{0\}$,
where $0$ denotes the group identity. With our convention
$\mathcal Cf=f*K$ and $P_\lambda f=f*p_\lambda$, put
\[
 \mathfrak h_\lambda:=p_\lambda*K-K*p_\lambda.
\]
Associativity of convolution gives $[\mathcal C,P_\lambda]f=f*\mathfrak h_\lambda$ initially for Schwartz functions. Distributional homogeneity
implies
\[
 \mathfrak h_\lambda(z)
 =\lambda^{-\mathfrak q}
       \mathfrak h(\delta_{1/\lambda}z),
 \qquad
 \mathfrak h:=p_1*K-K*p_1.
\]
This identity is initially understood in the sense of
distributions. We will show that $\mathfrak h$ is an
integrable smooth function.

To this aim, choose a smooth annular partition of unity
\[
 \psi_j(z):=\psi(\delta_{2^{-j}}z),
 \qquad
 \sum_{j\in\mathbb Z}\psi_j(z)=1
 \quad(z\ne0),
\]
where $\operatorname{supp}\psi
 \subset\{z:c_0\le\rho(z)\le C_0\}$ for fixed geometric constants $0<c_0<C_0$.
Here $\rho$ is the homogeneous quasi-norm fixed in
Subsection~\ref{subsec:kolmogorov-geometry}.
Set \[\ell_j:=2^j,\qquad K_j:=\psi_jk.\]
Since $\rho(\delta_{1/\ell_j}z)=\ell_j^{-1}\rho(z)$,
we have
\[
 \operatorname{supp}K_j
 \subset\operatorname{supp}\psi_j
 \subset\{z:c_0\ell_j\le\rho(z)\le C_0\ell_j\}.
\]
Moreover, homogeneity gives
\[
 K_j(z)=\ell_j^{-\mathfrak q}
                  K_0(\delta_{1/\ell_j}z).
\]
Thus, for every homogeneous left- or right-invariant
differential monomial $D$ of degree $v$,
\[
 \|DK_j\|_1\lesssim_D\ell_j^{-v},
 \qquad
 \|DK_j\|_\infty\lesssim_D\ell_j^{-\mathfrak q-v}.
\]
We also allow $D=I$, with $v=0$. For $D=I$, this gives a uniform $\mathrm L^1$ bound for
each annular piece $K_j$. It does not imply integrability
of the full kernel, since these bounds are not summable
over $j\in\mathbb Z$.

To relate these annular pieces to the distribution $K$, we
apply \cite[Proposition~3.2.27]{FischerRuzhanskyBook}.
Since $K$ is homogeneous of degree $-\mathfrak q$ and smooth
away from the identity, this proposition shows that its
restriction $k$ has vanishing homogeneous angular mean and that
\[
 K=\operatorname{pv}k+\widetilde c\,\delta_0
\]
for some $\widetilde c\in\C$.
Here the principal value is taken with respect to $\rho$,
\[
 \langle\operatorname{pv}k,\phi\rangle
 :=\lim_{\varepsilon\downarrow0}
       \int_{\{\rho(z)>\varepsilon\}}k(z)\phi(z)\d z,
 \qquad \phi\in\mathscr S(G).
\]
The vanishing angular mean gives, by homogeneous polar
coordinates,
\[
 \int_{\{a<\rho(z)<b\}}k(z)\d z=0,
 \qquad 0<a<b.
\]
Thus distributional homogeneity supplies the cancellation
needed at the identity.

Fix $\chi\in\mathrm C_0^\infty(G)$ equal to one near the
identity, and put $c:=\langle K,\chi\rangle$. Then
\[
 \langle K,\phi\rangle
 =c\phi(0)+\int_G
       k(z)\bigl(\phi(z)-\phi(0)\chi(z)\bigr)\d z,
 \qquad \phi\in\mathscr S(G).
\]
The integral is absolutely convergent as the expression in
parentheses is $O(\rho(z))$ near the identity and rapidly
decreasing at infinity. Equivalently, one has
\[
 K=c\delta_0+
       \sum_{j\in\mathbb Z}(K_j-c_j\delta_0)
 \quad\text{in }\mathscr S'(G),
 \qquad c_j:=\int_GK_j(z)\chi(z)\d z.
\]
Every Dirac term disappears from the commutator because
$p_1*\delta_0=\delta_0*p_1=p_1$. Hence
\begin{equation}\label{angulasum}
 \mathfrak h=\sum_{j\in\mathbb Z}\mathfrak h_j
 \quad\text{in }\mathscr S'(G),
 \qquad
 \mathfrak h_j:=p_1*K_j-K_j*p_1.
\end{equation}
Since $p_1,K_j\in\mathrm L^1(G)$, both convolutions
$p_1*K_j$ and $K_j*p_1$ are integrable.
Fubini's theorem and left invariance of Haar measure give
\[
 \int_G\mathfrak h_j(z)\d z
 =\left(\int_Gp_1\right)\left(\int_GK_j\right)
  -\left(\int_GK_j\right)\left(\int_Gp_1\right)
 =0.
\]

We next estimate the annular commutators. We use the identities relating left- and right-invariant
derivatives and the corresponding mean-value estimates,
see \cite[Corollary~3.1.30, Proposition~3.1.46 and
Remark~3.1.48(1)]{FischerRuzhanskyBook}. We record the particular consequence needed here.
For a fixed $\varphi\in\mathscr S(G)$, a differential monomial
$D$ as above, and $N\ge0$, there is an integer $M\ge0$ such that
\begin{equation}\label{w:schwartz-translation-difference}
 \begin{split}
 &\sup_{z\in G}(1+\rho(z))^N
 \left|D\bigl(\varphi(y^{-1}\circ z)-\varphi(z\circ y^{-1})\bigr)\right|\lesssim_{\varphi,D,N}\rho(y)(1+\rho(y))^M.
 \end{split}
\end{equation}
Here $D$ acts in the $z$-variable, and dependence on $\varphi$
involves only finitely many of its Schwartz seminorms.  To prove this when $\rho(y)\le1$, use
\cite[Lemma~3.1.47]{FischerRuzhanskyBook} to write
\[
 y=\exp(t_1Z_1)\circ\cdots\circ\exp(t_{2m+1}Z_{2m+1}),
 \qquad |t_\ell|\lesssim\rho(y)^{v_\ell},
\]
where $v_\ell$ is the homogeneous degree of $Z_\ell$
and the order of the basis fields is fixed.
Both translations agree at $y=0$. Applying the fundamental
theorem of calculus along the successive one-parameter
translations therefore bounds their difference in each
Schwartz seminorm by a fixed multiple of
\[
 \sum_{\ell=1}^{2m+1}|t_\ell|
 \lesssim\sum_{\ell=1}^{2m+1}\rho(y)^{v_\ell}
 \lesssim\rho(y),
\]
since $v_\ell\ge1$ and $\rho(y)\le1$. For $\rho(y)\ge1$, the polynomial group law and the
polynomial coefficients in the adjoint action imply that
the translated Schwartz seminorms grow at most polynomially
in $\rho(y)$. Increasing $M$ therefore gives
\eqref{w:schwartz-translation-difference} in this range as well.
The argument applies to both left and right derivatives.

For $\ell_j\le1$, write
\[
 \mathfrak h_j(z)
 =\int_GK_j(y)
       \bigl(p_1(z\circ y^{-1})-p_1(y^{-1}\circ z)\bigr)\d y.
\]
The two terms involve multiplication on opposite sides.
Their variation in $y$ is therefore controlled by left-
and right-invariant derivatives of $p_1$, respectively.
The translations here have size comparable to $\ell_j$.
Applying \eqref{w:schwartz-translation-difference} to $p_1$,
and using the support and $\mathrm L^1$ bounds for $K_j$, gives
\[
 |D\mathfrak h_j(z)|
 \lesssim_{D,N}\ell_j(1+\rho(z))^{-N}.
\]
Thus the small annuli gain a factor $\ell_j$ from the
difference of two nearby translates of the fixed smooth
kernel $p_1$. For $\ell_j\ge1$, use instead
\[
 \mathfrak h_j(z)
 =\int_Gp_1(y)
       \bigl(K_j(y^{-1}\circ z)-K_j(z\circ y^{-1})\bigr)\d y.
\]
Rescaling the annulus to unit size gives
\[
 \begin{split}
 \mathfrak h_j(\delta_{\ell_j}x)
 &=\ell_j^{-\mathfrak q}\int_Gp_1(y)
 \bigl[K_0((\delta_{1/\ell_j}y)^{-1}\circ x)-K_0(x\circ (\delta_{1/\ell_j}y)^{-1})\bigr]\d y.
 \end{split}
\]
A homogeneous derivative of degree $v$ contributes an
additional factor $\ell_j^{-v}$. Applying
\eqref{w:schwartz-translation-difference} to the fixed
smooth annular kernel $K_0$ therefore yields
\[
 \begin{split}
 |D\mathfrak h_j(z)|
 &\lesssim_{D,N}\ell_j^{-\mathfrak q-v}
       (1+\rho(z)/\ell_j)^{-N}\times
 \int_G|p_1(y)|\frac{\rho(y)}{\ell_j}
          \left(1+\frac{\rho(y)}{\ell_j}\right)^M\d y\\
 &\lesssim_{D,N}\ell_j^{-\mathfrak q-v-1}
       (1+\rho(z)/\ell_j)^{-N}.
 \end{split}
\]
The last step uses $\ell_j\ge1$ and the finite polynomial
moments of the Schwartz function $p_1$. This exhibits the
additional factor $\ell_j^{-1}$ gained on large annuli.

Above we have proved, for a differential monomial $D$ of degree $v$, that
\begin{equation}\label{w:annular-commutator-pointwise}
 |D\mathfrak h_j(z)|\lesssim_{D,N}
 \begin{cases}
 \ell_j(1+\rho(z))^{-N},
       &\ell_j\le1,\\
 \ell_j^{-\mathfrak q-v-1}
       (1+\rho(z)/\ell_j)^{-N},
       &\ell_j\ge1.
 \end{cases}
\end{equation}
Taking $N$ sufficiently large and integrating gives
\[
 \|D\mathfrak h_j\|_1
 \lesssim_D\min\{\ell_j,\ell_j^{-v-1}\}.
\]
In particular,
\[
 \|\mathfrak h_j\|_1
 \lesssim\min\{\ell_j,\ell_j^{-1}\}.
\]
For $0<\alpha<1$, the same estimates give
\[
 \int_G\rho(z)^\alpha|\mathfrak h_j(z)|\d z
 \lesssim_\alpha\min\{\ell_j,\ell_j^{\alpha-1}\}.
\]
Both bounds are summable over $j\in\mathbb Z$. The annular
series in \eqref{angulasum} therefore converges absolutely in $\mathrm L^1(G)$
and in $\mathrm L^1(G,\rho^\alpha\d z)$. Its distributional
limit identifies its sum with $\mathfrak h$. The derivative
estimates also justify termwise differentiation. Consequently,
\begin{equation}\label{w:commutator-integral-properties}
 \begin{gathered}
 \int_G\mathfrak h(z)\d z=0,\qquad
 \int_G(1+\rho(z))^\alpha|\mathfrak h(z)|\d z<\infty,\qquad
 Z_\ell^L\mathfrak h,\,
 Z_\ell^R\mathfrak h\in\mathrm L^1(G).
 \end{gathered}
\end{equation}
For $\rho(z)>1$, summing
\eqref{w:annular-commutator-pointwise} over
$\ell_j\le1$, $1<\ell_j\le\rho(z)$, and $\ell_j>\rho(z)$,
with $N>\mathfrak q+v+1$, gives
\begin{equation}\label{w:commutator-pointwise}
 |D\mathfrak h(z)|
 \lesssim_D(1+\rho(z))^{-\mathfrak q-v-1}.
\end{equation}
For $\rho(z)\le1$, this follows directly by summing the
uniform bounds in \eqref{w:annular-commutator-pointwise}. We now pass from these kernel properties to the square-function
estimate. Set
\[
 H_\lambda f:=f*\mathfrak h_\lambda,
 \qquad
 \mathfrak h^*(z):=\overline{\mathfrak h(z^{-1})}.
\]
The integrability of $\mathfrak h_\lambda$ makes $H_\lambda$
bounded on $\mathrm L^2(G)$. By density, it agrees with
$[\mathcal C,P_\lambda]$ on that space. Moreover,
$H_\lambda^*f=f*\mathfrak h_\lambda^*$.
Inversion preserves Haar measure and $\rho$, and exchanges
left and right derivatives up to sign. Thus
\eqref{w:commutator-integral-properties} also holds for
$\mathfrak h^*$.

We next use these kernel properties to establish almost
orthogonality between operators at different scales.
For example, the product $H_\lambda H_\sigma^*$ has
convolution kernel
$\mathfrak h_\sigma^**\mathfrak h_\lambda$.
When $\sigma\le\lambda$, cancellation of
$\mathfrak h_\sigma^*$ allows us to write this convolution
in terms of translation differences of
$\mathfrak h_\lambda$. We therefore first quantify how
the integrable left and right derivatives control such
differences in $\mathrm L^1$. The same estimates will
apply to the adjoint kernels by the preceding paragraph. Scaling gives
\[
 \|Z_\ell^L\mathfrak h_\lambda\|_1
 +\|Z_\ell^R\mathfrak h_\lambda\|_1
 \lesssim\lambda^{-v_\ell}.
\]
Expressing $y$ as an ordered product of exponentials and
integrating along the corresponding one-parameter subgroups
therefore gives, when $\rho(y)\le\lambda$,
\[
 \begin{split}
 &\|\mathfrak h_\lambda(y^{-1}\circ\cdot)-\mathfrak h_\lambda\|_1
 +\|\mathfrak h_\lambda(\cdot\circ  y^{-1})-\mathfrak h_\lambda\|_1\lesssim
 \sum_\ell\left(\frac{\rho(y)}{\lambda}\right)^{v_\ell}
 \lesssim\frac{\rho(y)}{\lambda}.
 \end{split}
\]
For $\rho(y)>\lambda$, use the uniform $\mathrm L^1$ norm
of $\mathfrak h_\lambda$. Combining the two ranges yields
\begin{equation}\label{w:commutator-translation-bound}
 \begin{split}
 &\|\mathfrak h_\lambda(y^{-1}\circ\cdot )-\mathfrak h_\lambda\|_1
 +\|\mathfrak h_\lambda(\cdot\circ  y^{-1})-\mathfrak h_\lambda\|_1\lesssim\min\{1,\rho(y)/\lambda\}.
 \end{split}
\end{equation}
The same bound holds for $\mathfrak h_\lambda^*$. Fix $0<\alpha<1$. For $0<\sigma\le\lambda$, cancellation
of $\mathfrak h_\sigma^*$ allows us to subtract the
untranslated kernel. By
\eqref{w:commutator-translation-bound},
\[
 \begin{split}
 \|\mathfrak h_\sigma^* *\mathfrak h_\lambda\|_1
 &\le\int_G|\mathfrak h_\sigma^*(y)|
       \|\mathfrak h_\lambda(y^{-1}\circ\cdot)
                         -\mathfrak h_\lambda\|_1\d y\\
 &\lesssim\int_G|\mathfrak h_\sigma^*(y)|
                      \min\{1,\rho(y)/\lambda\}\d y\lesssim\lambda^{-\alpha}
       \int_G\rho(y)^\alpha|\mathfrak h_\sigma^*(y)|\d y\lesssim_\alpha(\sigma/\lambda)^\alpha.
 \end{split}
\]
Here we used $\min\{1,t\}\le t^\alpha$ and the scaling of
the finite $\alpha$-moment. Right translation gives the
same bound for
$\|\mathfrak h_\lambda*\mathfrak h_\sigma^*\|_1$.
Thus the two convolution orders are estimated separately,
using the left and right translation bounds, respectively. For our right-convolution convention,
\[
 \begin{aligned}
 H_\lambda H_\sigma^*f
   &=f*(\mathfrak h_\sigma^* *\mathfrak h_\lambda),\qquad H_\lambda^*H_\sigma f
   &=f*(\mathfrak h_\sigma*\mathfrak h_\lambda^*).
 \end{aligned}
\]
Applying the preceding estimates also to the adjoint kernels,
using Young's inequality, and interchanging the scales gives
\[
 \|H_\lambda H_\sigma^*\|_{2\to2}
 +\|H_\lambda^*H_\sigma\|_{2\to2}
 \lesssim_\alpha
 \left(\frac{\min\{\lambda,\sigma\}}
            {\max\{\lambda,\sigma\}}\right)^\alpha,
 \qquad \lambda,\sigma>0.
\]

To deduce the square-function bound, fix
$I=[a,b]\subset(0,\infty)$ and define
\[
 \mathcal U_I\colon\mathrm L^2(G)
 \longrightarrow
 \mathrm L^2\!\bigl(I,{\lambda^{-1}}{\d\lambda};
                         \mathrm L^2(G)\bigr),
 \qquad
 (\mathcal U_If)(\lambda):=H_\lambda f.
\]
The uniform $\mathrm L^1$ bound for $\mathfrak h_\lambda$
makes $\mathcal U_I$ bounded for each such interval.
Its product with its adjoint acts by
\[
 (\mathcal U_I\mathcal U_I^*F)(\lambda)
 =\int_I H_\lambda H_\sigma^*F(\sigma)
                                      \frac{\d\sigma}{\sigma}.
\]
The preceding operator estimate and the identity
\[
 \int_0^\infty
 \left(\frac{\min\{\lambda,\sigma\}}
            {\max\{\lambda,\sigma\}}\right)^\alpha
 \frac{\d\sigma}{\sigma}
 =\frac2\alpha
\]
allow us to apply Schur's test. We obtain a bound for
$\mathcal U_I\mathcal U_I^*$, and hence for $\mathcal U_I$,
independent of $a$ and $b$. Letting $a\downarrow0$ and
$b\uparrow\infty$ gives
\begin{equation}\label{w:commutator-Haar-SF}
 \int_0^\infty
 \|H_\lambda f\|_2^2\frac{\d\lambda}{\lambda}
 \lesssim_\alpha\|f\|_2^2.
\end{equation}
Taking, for example, $\alpha=1/2$ fixes all constants.
Since $H_\lambda=[\mathcal C,P_\lambda]$, this proves
\eqref{w:commutator-SF}. \end{proof}

\noindent
{\bf Conflict of interest.}
The author states that there is no conflict of interest.

\medskip
\noindent
{\bf Data availability.}
Data sharing is not applicable to this article as no datasets
were generated or analysed during the current study.

\medskip
\noindent
{\bf Declaration of use of AI.}
The author conceived the ideas underlying
Lemma~\ref{lem:uncertainty} and
Lemma~\ref{lem:kinetic-hidden-coercivity} before writing
this manuscript. Building on these ideas, the author
planned to prove a kinetic Kato square root estimate
by adapting the parabolic arguments in
\cite{Nystrom2016,AtaeiEgertNystrom,AtaeiNystromBMO},
which in turn draw on methods from the solution of the
elliptic Kato square root problem in \cite{AHLMcT}.

The author encountered difficulties in developing suitable
smoothing operators and the associated Littlewood--Paley
estimates. During a discussion with GPT-6 Astra, developed
by OpenAI, the tool drew the author's attention to the
connection with positive Rockland operators. This theory
was unfamiliar to the author. With assistance from the tool,
the author developed the argument in
Section~\ref{sec:flat-rockland}, using the theory presented
in \cite{FischerRuzhansky,FischerRuzhanskyBook}.
The author then continued towards the proof of the main
theorem along the strategy of
\cite{Nystrom2016,AtaeiEgertNystrom}.

GPT-6 Astra was also used to help develop, revise and
complete the manuscript. This included assistance with
developing, checking and rewriting mathematical arguments,
drafting and revising text, and improving organization,
clarity and internal consistency. The author directed the
investigation and made the final mathematical and editorial
decisions. The author verified all mathematical arguments,
checked the references, and took particular care to
acknowledge prior work and place the results in their
historical context. The author acknowledges the AI
assistance described above and takes full responsibility
for the content of the manuscript.

\end{document}